\newif\ifprint
\printfalse
\documentclass[twoside,10pt]{HandbookOfModuli}

\usepackage{amsmath,amsthm}
\usepackage{amssymb}
\usepackage{latexsym}
\usepackage[utf8]{inputenc}
\usepackage{textcomp}
\usepackage{color}
\usepackage[pdftex]{graphicx}
\usepackage{titlesec}
\usepackage{xspace}

\ifprint
	\input{giovanni} 
	\usepackage[final]{microtype}
\fi
\usepackage{bm}
\renewcommand{\mathbf}[1]{\bm{#1}} 

\ifprint
	\definecolor{linkred}{rgb}{0,0,0} 
	\definecolor{linkblue}{rgb}{0,0,0} 
\else
	\definecolor{linkred}{rgb}{0.7,0.2,0.2}
	\definecolor{linkblue}{rgb}{0,0.2,0.6}
\fi

\numberwithin{equation}{section} 

\usepackage[
	hypertexnames=false,
	hyperindex,
	pagebackref,
	pdftex,
	pdftitle={Handbook of Moduli},
	pdfdisplaydoctitle,
	pdfpagemode=UseNone,
	breaklinks=true,
	extension=pdf,
	bookmarks=false,
	plainpages=false,
	colorlinks,
	linkcolor=linkblue,
	citecolor=linkred,
	urlcolor=linkred,
	pdfmenubar=true,
	pdftoolbar=true,
	pdfpagelabels,
	pdfpagelayout=TwoPage,
	pdfview=Fit,
	pdfstartview=Fit
]{hyperref}

\catcode`@=11 \baselineskip 4.5mm
\def\ps@handbook{\def\@oddhead{\hfill \leftmark \hfill\thepage }
\def\@evenhead{\thepage \hfill \rightmark \hfill}
\def\@oddfoot{}
\def\@evenfoot{}}
\def\@evenhead{}
\def\@oddfoot{}
\def\@evenfoot{\hfill\copyright\ China Higher Education Press}
\def\list#1#2{\ifnum \@listdepth >5\relax \@toodeep \else \global
\advance \@listdepth\@ne \fi \rightmargin \z@ \listparindent\z@
\itemindent\z@ \csname @list\romannumeral\the\@listdepth\endcsname
\def\@itemlabel{#1}\let\makelabel\@mklab \@nmbrlistfalse #2\relax
\@trivlist \parskip -\parsep \parindent\listparindent \advance
\linewidth -\rightmargin \advance\linewidth -\leftmargin \advance
\@totalleftmargin \leftmargin \parshape \@ne \@totalleftmargin
\linewidth \ignorespaces}
\renewcommand*\l@section{\@tocline{1}{0pt}{0em}{1em}{}}
\renewcommand*\l@subsection{\@tocline{2}{0pt}{1.5em}{2em}{}} 
\renewcommand\contentsnamefont{\bfseries\large} 
\catcode`@=12
\renewcommand{\theequation}{\thesection.\arabic{equation}}

\def\thebibliography#1{\section*{References}
\list{[\arabic{enumi}]}{\settowidth \labelwidth{[#1]} \leftmargin
\labelwidth \advance \leftmargin \labelsep \usecounter{enumi}}
\def\newblock{\hskip .11em plus .33em minus .07em} \sloppy
\clubpenalty 4000 \widowpenalty 4000 \sfcode`\.=1000 \relax}

\titleformat{\section}{\normalfont\large\bfseries}{\thesection.}{0.5em}{}[\kern0.em]
\titleformat{\subsection}{\normalfont\bfseries}{\thesubsection.}{0.3em}{}[\kern0.em]
\titleformat{\subsubsection}[runin]{\normalfont\bfseries}{\thesubsubsection.}{0.5em}{}[\kern0.5em]

\def\fofsubsubsection#1{\refstepcounter{equation}\subsubsection*{\theequation.\kern0.25em #1}}
\def\foisubsubsection#1{\refstepcounter{equation}\subsubsection*{\kern\parindent\theequation.\kern0.25em #1}}

\usepackage[all,ps,cmtip]{xy} 

\newtheorem{thm}[equation]{Theorem} 
\newtheorem{prop}[equation]{Proposition} 
\newtheorem{lem}[equation]{Lemma} 
\newtheorem{cor}[equation]{Corollary} 
\theoremstyle{remark}
\newtheorem{defn}[equation]{Definition} 
\newtheorem{rem}[equation]{Remark} 
\newtheorem{eg}[equation]{Example} 

\newcommand{\Z}{\mathbb{Z}}     
\newcommand{\N}{\mathbb{N}}     
\newcommand{\C}{\mathbb{C}}     

\newcommand{\sH}{\mathcal{H}}   
\newcommand{\sM}{\mathcal{M}}   
\newcommand{\sN}{\mathcal{N}}   
\newcommand{\sO}{\mathcal{O}}  

\theoremstyle{plain}
\newtheorem{theorem}[equation]{Theorem} 
\newtheorem{proposition}[equation]{Proposition} 
\newtheorem{lemma}[equation]{Lemma} 
\theoremstyle{remark}
\newtheorem{definition}[equation]{Definition} 
\newtheorem{remark}[equation]{Remark} 
\newtheorem{example}[equation]{Example} 
\newtheorem{question}[equation]{Question} 

\newcommand\cC{\mathcal{C}}
\newcommand\cD{\mathcal{D}}

\newcommand\cK{\mathcal{K}}
\newcommand\cL{\mathcal{L}}
\newcommand\cM{\mathcal{M}}

\newcommand\cO{\mathcal{O}}

\newcommand\cS{\mathcal{S}}

\renewcommand\AA{\mathbb{A}}

\newcommand\CC{\mathbb{C}}

\newcommand\HH{\mathbb{H}}

\newcommand\MM{\mathbb{M}}
\newcommand\NN{\mathbb{N}}
\newcommand\OO{\mathbb{O}}
\newcommand\PP{\mathbb{P}}
\newcommand\QQ{\mathbb{Q}}
\newcommand\RR{\mathbb{R}}

\newcommand\ZZ{\mathbb{Z}}

\newcommand{\ocM}{{\overline{\mathcal{M}}}}
\DeclareMathOperator{\Lift}{Lift}

\newtheorem*{philosophy}{Philosophy}

\newcommand{\til}[1]{\tilde{#1}}
\newcommand{\mf}[1]{\mathfrak{#1}}
\newcommand{\mgn}{\mathcal{M}_{g,n}}
\newcommand{\mgnb}{{\overline{\mathcal{M}}}_{g,n}}
\newcommand{\lmgnb}{{\overline{\mathcal{M}}}_{g,n}^{log}}

\newcommand{\sI}{\mathcal{I}}
\newcommand{\sS}{\mathcal{S}}
\newcommand{\sZ}{\mathcal{Z}}

\DeclareMathOperator{\Spec}{Spec}
\DeclareMathOperator{\im}{im}
\DeclareMathOperator{\coker}{coker}

\newcommand{\nc}{\newcommand}
\newcommand{\enm}{\ensuremath}

\nc{\m}{\enm}

\nc{\doc}{\begin{document}}
\nc{\docd}{\end{document}}
\nc{\mA}{\enm{\mathcal{A}}}
\nc{\mB}{\enm{\mathcal{B}}}
\nc{\mC}{\enm{\mathcal{C}}}
\nc{\mD}{\enm{\mathcal{D}}}
\nc{\mE}{\enm{\mathcal{E}}}
\nc{\mF}{\enm{\mathcal{F}}}
\nc{\mG}{\enm{\mathcal{G}}}
\nc{\mH}{\enm{\mathcal{H}}}
\nc{\mI}{\enm{\mathcal{I}}}
\nc{\mJ}{\enm{\mathcal{J}}}
\nc{\mK}{\enm{\mathcal{K}}}
\nc{\mL}{\enm{\mathcal{L}}}
\nc{\mM}{\enm{\mathcal{M}}}
\nc{\mN}{\enm{\mathcal{N}}}
\nc{\mO}{\enm{\mathcal{O}}}
\nc{\mP}{\enm{\mathcal{P}}}
\nc{\mQ}{\enm{\mathcal{Q}}}
\nc{\mR}{\enm{\mathcal{R}}}
\nc{\mS}{\enm{\mathcal{S}}}
\nc{\mT}{\enm{\mathcal{T}}}
\nc{\mU}{\enm{\mathcal{U}}}
\nc{\mV}{\enm{\mathcal{V}}}
\nc{\mW}{\enm{\mathcal{W}}}
\nc{\mX}{\enm{\mathcal{X}}}
\nc{\mY}{\enm{\mathcal{Y}}}
\nc{\mZ}{\enm{\mathcal{Z}}}

\nc{\bg}{\begin}
\nc{\ed}{\end}

\nc{\mcal}{\mathcal}

\nc{\elt}{element }

\nc{\enum}{\enumerate}
\nc{\eq}{\equation}
\nc{\etale}{\'etale }

\nc{\geom}{geometric }
\nc{\gpd}{groupoid }

\nc{\isom}{isomorphism }

\nc{\ladic}{$l$-adic}
\nc{\les}{long exact sequence }

\nc{\nbhd}{neighborhood }
\nc{\nth}{\ensuremath{n}-th }

\nc{\padic}{\ensuremath{p}-adic}
\nc{\pt}{point }

\nc{\scr}{\mathscr}
\nc{\ses}{short exact sequence }
\nc{\std}{standard }
\nc{\str}{structure }

\nc{\Thm}{Theorem }

\nc{\xet}{\enm{X_{\acute{e}t}}}

\nc{\zz}{\vspace{3ex}} 

\nc{\dc}{\DeclareMathOperator}
\dc{\Aut}{Aut}

\dc{\Coh}{Coh}
\dc{\Coker}{Coker}

\dc{\Div}{Div}

\dc{\End}{End}
\dc{\Exal}{Exal}
\dc{\Ext}{Ext}
\dc{\sExt}{} 

\dc{\GL}{GL}

\dc{\Hom}{Hom}
\dc{\sHom}{} 

\dc{\iso}{\cong}

\dc{\Ker}{Ker}

\dc{\PGL}{PGL}
\dc{\Pic}{Pic}
\dc{\Proj}{Proj}

\dc{\Tor}{Tor}

\nc{\into}{\hookrightarrow}
\nc{\onto}{\twoheadrightarrow}
\nc{\ra}{\rightarrow}

\nc{\Aone}{\enm{\mathbb{A}^1}} 
\nc{\An}{\enm{\mathbb{A}^n}} 

\nc{\Ga}{\enm{\mathbb{G}_\text{a}}} 
\nc{\Gm}{\enm{\mathbb{G}_\text{m}}} 
\nc{\Het}[1]{\enm{H^{#1}_\text{\'et}}} 

\nc{\oo}{\enm{\mathcal{O}}} 
\nc{\ox}{\enm{\mathcal{O}_X}} 
\nc{\osub}[1]{\enm{\mcal{O}_{#1}}} 

\nc{\algn}{\begin{aligned}} 
\nc{\algnd}{\end{aligned}}

\nc{\quot}[3][\Big/]{\raisebox{.6ex}{\enm{#2}}\!\!\ensuremath{#1}\!\raisebox{-1.3ex}{\ensuremath{#3}}}

\nc{\globalquotient}[2]{\left[\quot{#1}{#2}\right]}

\nc{\mor}[3]{\xymatrix{#1\ar[r]^{#2}&#3}}

\nc{\squarediagram}[8]{\xymatrix{
#1  \ar[r]^{#5}\ar[d]_{#7} &    #2  \ar[d]^{#8}  \\
#3  \ar[r]^{#6}            &    #4}}

\nc{\simplediagram}[4]{\xymatrix{
#1  \ar[r]\ar[d] &    #2  \ar[d]  \\
#3  \ar[r]       &    #4}}

\nc{\tab}[2][tbp]{\begin{table}[#1]\begin{center}\begin{tabular}{#2}}
\nc{\tabd}[1][]{\end{tabular}\caption{#1}\end{center}\end{table}} 

\theoremstyle{remark}
\newtheorem{rmk}[equation]{Remark}
\newtheorem*{conv}{Convention}

\newcommand{\be}{\begin{eqnarray*}}   
\newcommand{\ee}{\end{eqnarray*}}
\newcommand{\B}{\operatorname{B}}     
\newcommand{\Blo}{\operatorname{Blo}} 
\newcommand{\Id}{\operatorname{Id}}   
\renewcommand{\H}{\operatorname{H}}   
\renewcommand{\O}{\mathcal O}         

\theoremstyle{definition}
\newtheorem{subdefinition}[subsubsection]{Definition}
\newtheorem{subproposition}[subsubsection]{Proposition}
\newtheorem{sublemma}[subsubsection]{Lemma}
\newtheorem{subcorollary}[subsubsection]{Corollary}
\newtheorem{subtheorem}[subsubsection]{Theorem}
\newtheorem{subremark}[subsubsection]{Remark}
\newtheorem{subexample}[subsubsection]{Example}
\newtheorem{anitem}[equation]{} 
\usepackage{amsmath,amsthm,amssymb,amsxtra,amsfonts,amsbsy,mathrsfs}
\usepackage{verbatim}
\usepackage[all]{xy}

\begin{document}
\setcounter{page}{1}
%
%

\long\def\replace#1{#1}

%
%
\title{\replace{Logarithmic Geometry and Moduli}}
%
%
\author[Abramovich]{\replace{Dan Abramovich}}
\address{\replace{
Brown University\\
Box 1917\\
Providence, RI 02912\\
U.S.A.}}
\email[Abarmovich]{\replace{abrmovic@math.brown.edu}}

\author[Chen]{\replace{Qile Chen}}
\email[Chen]{\replace{q.chen@math.brown.edu}}

\author[Gillam]{\replace{Danny Gillam}}
\email[Gillam]{\replace{wdgillam@math.brown.edu}}

\author[Huang]{\replace{Yuhao Huang}}
\address{\replace{University of California, Berkeley}}
\email[Huang]{\replace{huang@math.berkeley.edu}}

\author[Olsson]{\replace{Martin Olsson}}
\email[Olsson]{\replace{molsson@math.berkeley.edu}}

\author[Satriano]{\replace{Matthew Satriano}}
\email[Satriano]{\replace{satriano@math.berkeley.edu}}

\author[Sun]{\replace{Shenghao Sun}}
\email[Sun]{\replace{shenghao@math.berkeley.edu}}
%
%

\begin{abstract}
	
\replace{We discuss the role played by logarithmic structures in the theory of moduli}. 

\end{abstract}  

\maketitle
\thispagestyle{empty}
%
%

\tableofcontents

\section{Introduction}\label{intro}

\subsection*{Logarithmic structures in algebraic geometry}

It can be said that Logarithmic Geometry is concerned with a method of
finding and
using ``hidden smoothness'' in singular varieties. The original
insight comes from consideration of de Rham cohomology. Since
singular varieties naturally occur ``at the boundary'' of many moduli
problems, logarithmic geometry was soon applied in the theory of
moduli. 

Foundations for this theory were first given by Kazuya Kato in
\cite{KKato}, following ideas of Fontaine and Illusie. The main body
of work on logarithmic geometry  has been concerned with deep
applications in the  cohomological study 
of $p$-adic and arithmetic schemes. This gave the theory an aura of
``yet another extremely complicated theory''.
The treatments of the theory are however quite accessible.
We hope to convince the reader here that the theory is simple
enough and useful enough to be considered by anybody interested in
moduli of singular varieties, indeed enough to be included in a
Handbook of Moduli.

\subsection*{Normal crossings and logarithmic smoothness} So what is
the original insight? Let $X$ be a
nonsingular complex variety, $S$ a curve with a point $s$ and $f:X \to
S$ a dominant morphism smooth away from $s$, in such a way that
$f^{-1}s=X_s=Y_1\cup\ldots\cup Y_m$ is a reduced simple normal crossings
divisor. Then of course $\Omega_{X/S}=\Omega_X/f^*\Omega_S$ fails to be locally free at the
singular points of $f$. But consider instead the sheaves
$\Omega_X(\log(X_s))$ of differential forms with at most logarithmic poles
along the $Y_i$, and similarly $\Omega_S(\log(s))$. Then  there is an
injective sheaf homomorphism $f^*\Omega_S(\log(s))\to
\Omega_X(\log(X_s))$, and {\em the quotient sheaf $\Omega_X(\log(X_s))/\Omega_S(\log(s))$
  is locally free}. 

So in terms of logarithmic forms, {\em the morphism $f$ is as good as a
smooth morphism.}

There is much more to be said: first, this $\Omega_X(\log(X_s))/\Omega_S(\log(s))$ can be
extended to a logarithmic de Rham complex, and its hypercohomology,
while not recovering the cohomology of the singular fibers, does give
rise to the limiting Hodge structure. So it is evidently worth
considering. 

Second, the picture is quite a
bit more
general, and can be applied to all toric and toroidal maps between
toric varieties or toroidal embeddings (with a
little caveat about the characteristic of the residue fields). So
there is some flexibility in choosing $X\to S$.

\subsection*{The search for a structure}

Since we are considering moduli, then as soon as we consider $X
\to S$ as above we must also consider the normal crossings fiber $X_s
\to \{s\}$.  But what structure should we put on this variety? The
notion of differentials with logarithmic poles along $X_s$ is not
in itself intrinsic to $X_s$. Also
the
normal crossings variety $X_s$ is not in itself toric or 
toroidal, so a new structure is needed to incorporate it into the
picture.

One is tempted to consider varieties which are assembled from nice
variety by some sort of gluing, as normal crossings varieties are. But
already normal crossings varieties do not give a satisfactory
answer in general, because their deformation spaces have ``bad''
components. Here is a classical example: consider a smooth projective  variety $Z$
such that $Pic^0(Z)$  
is nontrivial. Let $L$ be a line bundle on $Z$ and set
$Y=\PP(\cO\oplus L)$, with zero section $Z\subset Y$.  Let $X$ be the
blowing up of $Z\times 0\subset Y \times \AA^1$. We have a flat morphism $f:X
\to \AA^1$ with fiber  $X_0 = f^{-1}(0) \simeq Y\cup Y$, where the two copies
of $Y$ are glued with the zero section of one attached to the $\infty$
section of the other. 

So clearly $X_0$ is a normal crossings variety with a nice smoothing
to a copy of $Y$. But there are other deformations: the variety $Y\cup
Y$ also deforms to $Y \cup Y'$ where $Y' = \PP(\cO\oplus L')$ and $L'$
a deformation of the line bundle $L$. And it is not hard to see that $Y \cup Y'$ does
not have a smoothing. Ideally one really does not want to see this deformation
$Y\cup Y'$ in the picture - and ideally $X_0$ should have a natural
structure whose deformation space excludes $Y\cup Y'$ automatically.

Such a structure was proposed by Friedman in \cite{Friedman}, where
the notion of {\em d-semistable varieties} was introduced. This structure
is somewhat subtle, and while it solves the issue in this case, it is
not quite as flexible as one could wish. As we will see in Section
\ref{D-SS}, logarithmic structures subsume d-semistability and
do provide an appropriate flexibility.

\subsection*{Organization of this chapter} 
In this chapter we briefly describe logarithmic structures and
indicate where they can be useful in the study of moduli spaces.   
Section \ref{Qile1} gives the basic definitions of  logarithmic
structures, and section \ref{Qile2} discusses logarithmic
differentials and log smooth deformations, which are important in
considering moduli spaces.

 Section \ref{Satriano}  
gives the first example where logarithmic geometry fits well with
moduli spaces: the moduli space of stable curves is the moduli space
of log smooth curves. The issue of d-semistability does not arise
since a nodal curve is automatically d-semistable. So the theory for
curves is simple.   
Turning to higher dimensions, Section \ref{D-SS} shows how
d-semistability can be described using logarithmic structures. 

If one is to enlarge algebraic geometry to include logarithmic
structures, the task of generalizing the techniques of algebraic
geometry to logarithmic structure can certainly seem daunting. In
section \ref{LogStacks} we show how to encode logarithmic structure in
terms of certain algebraic stacks. This allows us to reduce various
constructions to the case of algebraic stacks. (One can argue that the
theory of stacks is not simple either, but at least in the theory of
moduli they have come to be accepted, with some exceptions
\cite{Kollar-companion}.)  

In section \ref{LogDef} we make use of logarithmic stacks to describe
the complexes which govern deformations and obstructions for
logarithmic structures even in the non-smooth case. This comes in
handy later. For instance, even when studying moduli of log-smooth
schemes, the moduli spaces tend to be singular, and their cotangent
complexes are a necessary ingredients in constructing virtual
fundamental classes. 

Section \ref{Rounding} describes a beautiful construction, similar to
polar coordinates, in which families of complex log smooth varieties
give rise canonically to families of topological
manifolds. Differential geometers have used polar coordinates on nodal
curves to ``make space" for monodromy to act by Dehn twists. Rounding
(using Ogus's terminology) is a magnificent way to generalize this. 

The immediate implications of logarithmic structures for De Rham
cohomology and Hodge structures are described in  Section
\ref{DeRham}. 

We conclude by describing three applications, where logarithmic
structures serve as the proverbial ``magic powder" (term suggested by
Kato and Ogus) to clarify or remove unwanted behavior from moduli
spaces. 

Section \ref{MainComp} describes a number of cases where the
main irreducible component of a moduli space can be separated from
other ``unwanted" components by sprinkling the objects with a bit of
logarithmic structures. 

In Section \ref{Roots} we introduce twisted curves, a central object of orbifold stable maps, and show how logarithmic structures  give a palatable way to construct the moduli stack of twisted curves. 

Section \ref{Stablemaps} gives background for the work of B. Kim,
in which Jun Li's moduli space of relative stable maps, with its
obstruction theory and virtual fundamental class, is beautifully
simplified using logarithmic structures. 

\subsection*{Notation} Following the lead of Ogus \cite{Ogus}, we try whenever possible to denote a logarithmic scheme by a regular letter (such as $X$) and the underlying scheme by $\underline X$.  When this is impossible we write $X$ for the underlying scheme and $(X,\cM_X)$ for a logarithmic scheme over it.

\subsection*{Acknowledgements} This chapter originated from lectures given by Olsson at the  School and Workshop on Aspects of Moduli, June 15-28, 2008 at the De Giorgi Center at the Scuola Normale Superiore in Pisa, Italy.  The material was revisited and expanded in our seminar during the Algebraic Geometry program at MSRI, 2009. We thank the De Giorgi Center, MSRI, their  staff and program organizers  for providing these opportunities. Thanks are due to Arthur Ogus and Phillip Griffiths, who lectured on two topics at the MSRI seminar. While no new material is intended here, we acknowledge that research by Abramovich, Gillam and Olsson is supported by the NSF.  Finally thanks are due to Gavril Farkas and Ian Morrison for the invitation to write this chapter.


\section{Definitions and basic properties}\label{Qile1}
In this  section we introduce the basic definitions of logarithmic geometry in the sense of \cite{KKato}. Good introductions are given in \cite{KKato} and \cite{Ogus}. Further technique is developed in \cite{Gabber-Ramero}. 

\subsection*{Logarithmic structures} 

The basic definitions are as follows:

\begin{defn} 
A \emph{monoid} is a commutative semi-group with a unit. A morphism of monoids is required to preserve the unit element. We use $Mon$ to denote the category of Monoids.
\end{defn}

\begin{defn}\label{Def:log-str} 
Let $\underline{X}$ be a scheme. A \emph{pre-logarithmic structure} on $\underline{X}$ is a sheaf of
monoids $\sM_{X}$ on the \'etale site $\underline{X}_{\acute{et}}$ combined with a
morphism of sheaves of monoids: $\alpha
: \sM_{X} \longrightarrow \sO_{\underline{X}}$, called the structure morphism, where
we view $\sO_{\underline{X}}$ as a monoid under multiplication. A pre-log
structure is called a \emph{log structure} if
$\alpha^{-1}(\sO_{\underline{X}}^{*})\cong \sO^{*}_{\underline{X}}$ via $\alpha$. The
pair $(\underline{X},\sM_{X})$ is called a \emph{log scheme}, and will be denoted by $X$.
\end{defn}

Note that, given a log structure $\sM_{X}$ on $\underline{X}$,  we can view $\sO^{*}$ as a subsheaf $\sM_{X}$.

\begin{defn}
\label{Def:chara}
Given a log scheme $X$, the quotient sheaf $\overline{\sM}_{X}=\sM_{X}/\sO_{\underline{X}}^{*}$ is called the {\em characteristic of the log structure} $\sM_{X}$.
\end{defn}

\begin{defn} 
Let  $\sM$ and $\sN$ be pre-log structures on $\underline{X}$. A {\em morphism}
between them is a morphism $\sM \rightarrow \sN$ of sheaves of monoids
which is compatible with the structure morphisms. 
\end{defn}

How should one think of such a beast? There are two extreme cases: 

\begin{enumerate} 
\item If an element $m\in \sM$ has $\alpha(m)=x \neq 0$, one often thinks
of $m$ as some sort of partial data of a ``branch of the logarithm of
$x$''. Evidently no data is added if $x$ is invertible, but some is
added otherwise. In particular, we will see later that $m$ permits us
to take the logarithmic 
differential $dx/x$ of $x$.

\item If $\alpha(m)=0$ it is often the case that it $m$ comes by
restricting the log structure of an ambient space, and serves as the
``ghost'' of a logarithmic cotangent vector coming from that space. So
the log structure ``remembers'' deformations that are lost when
looking at the underlying scheme.
\end{enumerate}

\subsection*{The log structure associated to a pre-log structure} 
We have a natural inclusion $$i :  (\text{log structures on } 
\underline{X}) \hookrightarrow (\text{pre-log structures on } 
\underline{X})$$ by viewing a log structure as a pre-log structure. We now construct a left adjoint.

Let  $\alpha : \sM \rightarrow \sO_{\underline{X}}$ be a pre-log structure on $\underline{X}$. We define the  \emph{associated  log structure} $\sM^{a}$ to be the push-out of 
\[\xymatrix{
\alpha^{-1}(\sO_{\underline{X}}^{*}) \ar[d] \ar[r] & \sM \\
\sO_{\underline{X}}^{*}
}
\]
in the category of sheaves of monoids on $\underline{X}_{\acute{et}}$, endowed with $$\sM^{a} \rightarrow \sO_{\underline{X}}  \ \ (a,b)\mapsto \alpha(a)b \ (a \in \sM, b \in \sO_{\underline{X}}^{*}).$$

In this way, we obtain a functor $a$ : (pre-log structures on $\underline{X}$) $\rightarrow$ (log structures on $\underline{X}$). From the universal property of push-out, any morphism of pre-log structure from a pre-log structure $\sM$ to a log structure on $\underline{X}$ factor through $\sM^{a}$ uniquely.

\begin{lem}\cite[1.1.5]{Ogus}
The functor $a$ is left adjoint to $i$.
\end{lem}

\begin{eg}
The category of log structures on $\underline{X}$ has an initial object, called the trivial log structure, given by the inclusion $\sO_{\underline{X}}^{*} \hookrightarrow \sO_{\underline{X}}$.  It also has a final object, given by the identity map $\sO_{\underline{X}} \rightarrow \sO_{\underline{X}}$. Trivial log structures are quite useful as they make the category of schemes into a full subcategory of the category of log schemes (see Definition \ref{Def:log-mor}). The final object is rarely used since it is not fine, see definition \ref{Def:fine}.
\end{eg}

\begin{eg}\label{NClog}
Let $\underline{X}$ be a regular scheme, and $D \subset \underline{X}$  a
divisor. We can define a log structure $\sM$ on $\underline{X}$ associated to the
divisor $D$ as 
$$\sM(U)\ =\ \big\{\,g\in \sO_{\underline{X}}(U)\, :\ g|_{U\setminus D} \in \cO_{\underline{X}}^*(U\setminus D)\,\big\}\ \subset\ \cO_{\underline{X}}(U).$$ 

The case where $D$ is a normal
crossings divisor is special - we will see later that it is {\em log smooth}.

Note that the concept of normal crossing is
local in the \'etale topology. This is one reason we use
the \'etale topology instead of the Zariski topology. 
\end{eg}

\begin{eg}\label{Afflog}
Let $P$ be a monoid, $R$ a ring, and denote by $R[P]$  the monoid algebra. 
Let $\underline{X} = \Spec R[P]$. Then $\underline{X}$ has a
canonical log structure associated to the canonical map $P \rightarrow
R[P]$. We denote by $\Spec\ (P\rightarrow R[P])$ the log
scheme with underlying $\underline{X}$, and the canonical log structure. 
\end{eg}

\subsection*{The inverse image and the category of log schemes} 

Let $f:\underline{X}\rightarrow \underline{Y}$ be a morphism of schemes. Given a log structure
$\sM_{Y}$ on $\underline{Y}$, we can define a log structure on $\underline{X}$, called the
invese image of $\sM_{Y}$, to be the log structure associated to the
pre-log structure $f^{-1}(\sM_{Y})\rightarrow
f^{-1}(\sO_{\underline{Y}})\rightarrow \sO_{\underline{X}}$. This is usually denoted by
$f^{*}(\sM_{Y})$. Using the inverse image of log structures, we can
give the following definition. 

\begin{defn}\label{Def:log-mor} 
A \emph{morphism of log schemes} $X\rightarrow Y$
consists of a morphism of underlying schemes $f:\underline{X}\rightarrow \underline{Y}$, and a
morphism  $f^{\flat}: f^{*}\sM_{Y}\rightarrow \sM_{X}$ of log
structures on $\underline{X}$.

We denote by $LSch$ the category of 
log schemes.  
\end{defn}

\begin{eg}
In Example \ref{Afflog}, the log structure on $\Spec\ (P\rightarrow
R[P])$ can be viewed as the inverse image of the log structure on
$\Spec\ (P\rightarrow \Z[P])$ via the canonical map
$\Spec\ (R[P]) \rightarrow \Spec\ (\Z[P])$. 
\end{eg}

\begin{eg}\label{Logpt} 
Let $k$ be a field, $\underline{Y}$=Spec $k[x_{1},\cdots,x_{n}]$,
$D$=$V(x_{1}\cdots x_{r})$. Note that $D$ is a normal crossing divisor
in $\underline{Y}$. By example \ref{NClog}, we have a log structure $\sM_{Y}$ on
$\underline{Y}$ associated to the divisor $D$. In fact, $\sM_{Y}$ can be viewed as
a subsheaf of $\sO_{\underline{Y}}$ generated by $\sO_{\underline{Y}}^{*}$ and
$\{x_{1},\cdots,x_{r}\}$. 

Consider the inclusion $j: p$ = Spec $k \hookrightarrow \underline{Y}$ sending the
point to the origin of $\underline{Y}$. Then $j^{*}\sM_{Y} = k^{*} \oplus \N^{r}$,
and the structure map $j^{*}\sM\longrightarrow \sO_{\underline{X}}$ is given by
$(a,n_{1},\cdots,n_{r})\mapsto a\cdot 0^{n_{1}+\cdots+n_{r}}$, where
we define $0^{0} = 1$ and $0^{n}=0$ if $n\neq 0$. Such point with the
log structure above is call a {\em logarithmic point}; when $r = 1$ we call
it the {\em standard logarithmic point}. 
\end{eg}

\subsection*{Charts of log structures} 

\begin{defn}\label{Def:chart} 
Let $X$ be a log scheme, and $P$ a monoid. A {\em chart} for $\sM_{X}$
is a morphism $P\rightarrow \Gamma(X,\sM_{X})$, such that the induced
map of log strucutres $P^{a}\to \sM_{X}$ is an isomorphism, where
$P^{a}$ is the log structure associated to the pre-log structure given
by
$P\rightarrow \Gamma(X,\sM_{X})\rightarrow \Gamma(X,\sO_{\underline{X}})$.  
\end{defn}

In fact, a chart of $\sM_{X}$ is equivalent to a morphism $$f:X \rightarrow\Spec\ (P\rightarrow\Z[P]),$$ such that $f^{\flat}$ is an isomorphism. In general, we have the following:

\begin{lem}\cite[1.1.9]{Ogus} The morphism
$$Hom_{LSch}(X,\Spec\ (P\rightarrow \Z[P]))\ \to \ Hom_{Mon}(P,\Gamma(X,\sM_{X}))$$
associating to $f$ the composition
 $$\xymatrix{P \ar[r]&\Gamma(\underline{X}, P_X)\ar[r]^{\Gamma(f^\flat)} & \Gamma(\underline{X}, \cM_X)}$$
 is an isomorphism.
\end{lem}

We can also consider  charts for log morphisms.

\begin{defn}\label{Def:chart-mor} 
Let $f:X\rightarrow Y$ be a morphism of log schemes. A {\em chart} for $f$ is a triple $(P_{X}\rightarrow\sM_{X},Q_{Y}\rightarrow\sM_{Y},Q\rightarrow P)$ where $P_{X}$ and $Q_{Y}$ are the constant sheaves associated to the monoids $P$ and $Q$, which satisfy the following conditions:
 \begin{enumerate}
   \item $P_{X}\rightarrow\sM_{X}$ and $Q_{Y}\rightarrow\sM_{Y}$ are charts of $\sM_{X}$ and $\sM_{Y}$;
   \item the morphism of monoids $Q\rightarrow P$ makes the following diagram commutative:
\[
\xymatrix{
Q_{X} \ar[r]\ar[d] & P_{X} \ar[d]\\
f^{*}\sM_{Y} \ar[r]& \sM_{X}.
}
\]
 \end{enumerate}
\end{defn}

\subsection*{Fine log structures} 
Arbitrary log structures are too wild to manipulate; they are roughly
analogous to arbitrary ringed spaces: both notions are useful for
general constructions, but a narrower, more geometric category is
desirable. In Definition \ref{Def:fine} below we introduce the notion
of fine log structures.  Continuing the analogy above, these are 
well-behaved log structures
analogous to noetherian 
schemes, in the sense that you can do geometry on them.   

Given a monoid $P$, we can associate a group 
\begin{equation*}{
P^{gp}:=\{(a,b) |(a,b)\sim(c,d) \ \mbox{if} \ \exists s \in
P \ \mbox{such that} \ s+a+d=s+b+c\}.} 
\end{equation*}
 Note that any morphism from $P$ to an abelian group factors through $P^{gp}$ uniquely.

\begin{defn}\label{Def:integral} 
$P$ is called {\em integral} if $P\rightarrow P^{gp}$ is injective. It is called \emph{saturated} if it is integral and for any $p\in P^{gp}$, if $n\cdot p \in P$ for some positive integer $n$ then $p\in P$.
\end{defn}

\begin{defn}\label{Def:fine} 
A log scheme $X$ is said to be {\em fine}, if \'etale locally there is a chart $P\rightarrow \sM_{X}$ with $P$ a finitely generated integral monoid. If moreover $P$ can be choosen to be saturated, then $X$ is called a {\em fine and saturated} (or {\em fs}) log structure. This is equivalent to saying that for every geometric point $\bar{x} \to \underline{X}$ the monoid $\ocM_{\bar{x}, X}$ is saturated as in Definition \ref{Def:integral}. Finally if $P\simeq \NN^k$ we say that the log strcture is {\em locally free}.
\end{defn}

In the following, we will focus on fine log schemes.

\section{Differentials, smoothness, and log smooth deformations}\label{Qile2}

The main reference in this section is \cite{KKato}.

\subsection*{Logarithmic differentials} 

In \cite{EGA_IV} Grothendieck defines a derivation as the difference
of infinitesimal liftings of a  section. We can do the same thing with
logarithmic schemes. First, 
we need a concept of infinitesimal extension, which requires the
following definition.  

\begin{defn} 
A morphism $f: X\rightarrow X$ of log schemes is called {\em strict} if $f^{\flat}: f^{*}\sM_{Y}\rightarrow \sM_{X}$ is an isomorphism. It is called a {\em strict closed immersion} \footnote{The term used in \cite{KKato} is {\em an exact closed immersion}.} if it is strict and the underlying map $\underline{X}\rightarrow \underline{Y}$ is a closed immersion in the usual sense.
\end{defn}

Let us consider a commutative diagram of solid arrows of log schemes:

\[
\xymatrix{
T_{0} \ar[r]^{\phi} \ar[d]^{j}_{J}& X \ar[d]^{f}\\
T_{1} \ar[r]_{\psi} \ar@<1ex>@{.>}[ur]^{g_{1}} \ar@<-1ex>@{.>}[ur]_{g_{2}} & Y
}
\]
where $j$ is a strict closed immersion defined by an ideal $J$ with 
$J^{2}=0$. Note that $T_{0}$ and $T_{1}$ have the same underlying topological space, and isomorphic \'etale sites. Then we have the following commutative diagram of sheaves of algebras:

\[
\xymatrix{
\sO_{\underline{T}_{0}} & \phi^{-1}\sO_{\underline{X}} \ar[l] \ar@<-1ex>@{.>}[dl]_{g_{1}^{\#}} \ar@<1ex>@{.>}[dl]^{g_{2}^{\#}}\\
\sO_{\underline{T}_{1}} \ar[u] & \psi^{-1}\sO_{\underline{Y}} \ar[u] \ar[l]
}
\]
Then $g_{1}^{\#}-g_{2}^{\#}$ is a derivation $\partial_{g_{1}-g_{2}}: \phi^{-1}\sO_{\underline{X}}\rightarrow J$ in the usual sense. We also have a commutative diagram given by the log structures:

\[
\xymatrix{
\sM_{T_{0}} & \phi^{-1}\sM_{X} \ar[l] \ar@<-1ex>@{.>}[dl]_{g_{1}^{\flat}} \ar@<1ex>@{.>}[dl]^{g_{2}^{\flat}}\\
\sM_{T_{1}} \ar[u] & \psi^{-1}\sM_{Y} \ar[u] \ar[l]
}
\]

Note that we have an ``exact sequence" of mutiplicative
monoids $$\textbf{1}\rightarrow (1+J)\rightarrow
\sM_{T_{1}}\rightarrow \sM_{T_{0}}\rightarrow \textbf{1},$$ by which we mean that the group $1+J$ acts freely on $\cM_{T_1}$ with quotient $\cM_{T_0}$. Hence we
obtain a morphism $D_{g_{1}-g_{2}}:\phi^{-1}\sM_{X}\rightarrow J$ such
that for every $m \in \phi^{-1}\sM_{X}$ we have
$(g^{\flat}_{1}-g^{\flat}_{2})(m) = 1 + D_{g_{1}-g_{2}}(m).$  It is not 
hard to check that it is a monoid homomorphism:
$D_{g_{1}-g_{2}}(m\cdot n)=D_{g_{1}-g_{2}}(m)+D_{g_{1}-g_{2}}(n)$ for
any $m,n \in \phi^{-1}(\sM_{X})$. By the definition of log structures,
we also have  
\begin{enumerate}
 \item $\alpha(m)D_{g_{1}-g_{2}}m = \partial_{g_{1}-g_{2}}(\alpha(m))$, $\forall m\in \phi^{-1}\sM_{X}$;
 \item $D_{g_{1}-g_{2}}|_{\psi^{-1}\sM_{Y}}=0$.
\end{enumerate}

\begin{rem}
 \begin{enumerate}
       \item Since the log structure contains all the invertible
         elements in the structure sheaf, the map $D_{g_{1}-g_{2}}$
         determines $\partial_{g_{1}-g_{2}}$. 
       \item The above properties show that $D_{g_{1}-g_{2}}$ behaves
         like ``$d\log$". This is one of the reasons for the name
         ``logarithmic structure''. 
 \end{enumerate}
\end{rem}

Summarizing the above discussion gives the following definitions:
\begin{defn}\cite{Ogus-crystals},\cite[Definition 1.1.1]{Ogus} \label{LogDer} 
Consider a morphism $f: X \rightarrow Y$ of fine log schemes. Let $I$ be 
an $\sO_{\underline{X}}$-module. A {\em log derivation of $X$ over $Y$ to 
$I$} is a pair $(\partial, D)$ where $\partial \in 
\cD er_{\underline{Y}}(\underline{X},I)$ and $D:\sM_{X}\rightarrow I$ is an additive map such that the following conditions hold:
 \begin{enumerate}
   \item $D(ab)=D(a)+D(b)$ for $a,b\in\sM_{X}$;
   \item $\alpha(a)D(a)=\partial(\alpha(a))$, for $a\in\sM_{X}$.
   \item $D(a)=0$, for $a\in f^{-1}\sM_{Y}$.
 \end{enumerate}
The sheaf $\cD er_{Y}(X,I)$ of log derivations of $X$ over $Y$ to $I$ is the sheaf of germs of pairs $(\partial, D)$. The sheaf $\cD er_{Y}(X,\sO_{\underline{X}})$ is usually denoted by $T_{X/Y}$, and is called {\em the logarithmic tangent sheaf of $X$ over $Y$}. 
\end{defn} 

As an analogue of differentials of usual schemes, we have the following result:

\begin{prop}\label{prop-logdiff}\cite[IV.1.1.6]{Ogus} 
of Log differentials
There exists an $\sO_{\underline{X}}$-module $\Omega_{X/Y}^{1}$ with a universal derivations $(\partial,D)\in \cD er_{Y}(X,\Omega_{X/Y}^{1})$, such that for any $\sO_{\underline{X}}$-module $I$, the canonical map
\[\sH om_{\sO_{\underline{X}}}(\Omega_{X/Y}^{1},I)\rightarrow \cD er_{Y}(X,I), \ \ \ u\mapsto (u\circ \partial,u\circ D)\]
is an isomorphism of $\sO_{\underline{X}}$-modules. In fact, we have the following construction: \[\Omega_{X/Y}^{1}=\Omega_{\underline{X}/\underline{Y}}\oplus(\sO_{\underline{X}}\otimes_{\Z}\sM^{gp}_{X})/\mathcal{K}\]
where $\mathcal{K}$ is the $\sO_{\underline{X}}$-module generated by local 
sections of the following forms:
 \begin{enumerate}
   \item $(d\alpha(a),0)-(0,\alpha(a)\otimes a)$ with $a\in\sM_{X}$;
   \item $(0,1\otimes a)$ with $a\in Im(f^{-1}(\sM_{Y})\rightarrow\sM_{X})$.
 \end{enumerate}
The universal derivation $(\partial,D)$ is given by $\partial : \sO_{\underline{X}}\stackrel{d}\rightarrow\Omega_{\underline{X}/\underline{Y}}\rightarrow\Omega_{X/Y}^{1}$ and $D: \sM_{X}\rightarrow \sO_{\underline{X}}\otimes_{\Z}\sM_{X}^{gp}\rightarrow\Omega_{X/Y}^{1}$.
\end{prop}

\begin{defn} \label{LogDiff} 
Given a morphism $f:X\to Y$ of log schemes, the $\sO_{\underline{X}}$-module $\Omega_{X/Y}^{1}$ is called {\em the sheaf of logarithmic differentials}. Sometimes we use the short notation $\Omega_{f}^{1}$ for $\Omega_{X/Y}^{1}$.

Note that $\sH om(\Omega_{X/Y}^{1},\sO_{\underline{X}})\cong T_{X/Y}$.
\end{defn}

\begin{rem} 
If we consider only fine log structures, and assume that $\underline{Y}$ is locally noetherian and $\underline{X}$ locally of finite type over $\underline{Y}$, then both $\cD er_{Y}(X,I)$ and $\Omega_{X/Y}^{1}$ in the definitions above are coherent sheaves. The proof of this can be found in \cite[IV.1.1]{Ogus}
\end{rem}

\begin{eg}\label{NCsmooth}
Consider $R=k[x_{1},\cdots,x_{n}]/(x_{1}\cdots x_{r})$, where $k$ is a field. Denote $\underline{X} = \mbox{Spec}\ R$. Let $\sM_{X}$ be the log structure on $\underline{X}$ given by $\N^{r}\rightarrow R$, $e_{i}\mapsto x_{i}$, where $e_{i}$ is the standard generator of the monoid $\N^{r}$. Let $Y=$\ \mbox{Spec}\ $(\N\rightarrow k)$ be the logarithmic point described in \ref{Logpt}. Now we can define a morphism $f:X\rightarrow Y$ by the following diagram:
\[
\xymatrix{
\N^{r} \ar[r] & R \\
\N \ar[r] \ar[u]^{\Delta} & k \ar[u]
}
\]
where $\Delta : e \mapsto e_{1}+\cdots + e_{r}$, and $e$ is the standard generator of $\N$. Then it is easy to see that $\cD er_{Y}(X,\sO_{\underline{X}})$ is a free $\sO_{\underline{X}}$-module generated by 
\[x_{1}\frac{\partial}{\partial x_{1}},\cdots,x_{r}\frac{\partial}{\partial x_{r}},\frac{\partial}{\partial x_{r+1}},\cdots,\frac{\partial}{\partial x_{n}},\] 
with a relation $x_{1}\frac{\partial}{\partial x_{1}}+\cdots+x_{r}\frac{\partial}{\partial x_{r}}=0$. The sheaf $\Omega_{f}^{1}$ is a free $\sO_{X}$-module generated by the logarithmic differentials: 
\[\frac{dx_{1}}{x_{1}},\cdots,\frac{dx_{r}}{x_{r}},dx_{r+1},\cdots,dx_{n}\] 
with a relation $\frac{dx_{1}}{x_{1}}+\cdots+\frac{dx_{r}}{x_{r}}=0$.
\end{eg}

\begin{eg}
Let $h: Q \rightarrow P$ be a morphism of fine monoids. Denote $X=\ $Spec\ $(P\rightarrow\Z[P])$ and $Y=$ Spec\ $(Q\rightarrow\Z[Q])$. Then we have a morphism $f: X\rightarrow Y$ induced by $h$. A direct calculation shows that $\Omega_{f}^{1}=\sO_{\underline{X}}\otimes \textbf{Cok} (h^{gp})$. This can also be seen from the universal property of the sheaf of logarithmic differentials.
\end{eg}

\subsection*{Logarithmic Smoothness}

Let us go back to the following cartesian diagram of log schemes:

\begin{equation}\label{diag:smooth}
\xymatrix{
T_{0} \ar[r]^{\phi} \ar[d]^{j}_{J}& X \ar[d]^{f}\\
T_{1} \ar[r]_{\psi}  & Y
}
\end{equation}
where $j$ is a strict closed immersion defined by $J$ with $J^{2}=0$. As in the usual case, we can define log smoothness by the infinitesimal lifting property.

\begin{defn}\label{defn:logsmooth}
A morphism $f:X\rightarrow Y$ of fine log schemes is called {\em log smooth} (resp. {\em \'etale}) if the underlying morphism $\underline{X}\rightarrow \underline{Y}$ is locally of finite presentation and for any commutative diagram (\ref{diag:smooth}), \'etale locally on $T_{1}$  there exists a (resp. there exists a unique) morphism $g:T_{1}\rightarrow X$ such that $\phi=g\circ j$ and $\psi=f\circ g$.
\end{defn}

We have the following useful criterion for smoothness from \cite[Theorem 3.5]{KKato}.

\begin{thm}\label{KatoStrThm} 
(K.Kato) Let $f:X\rightarrow Y$ be a morphism of fine log schemes. Assume we have a chart $Q\rightarrow \sM_{Y}$, where $Q$ is a finitely generated integral monoid. Then the following are equivalent:
 \begin{enumerate}
  \item $f$ is log smooth (resp. log \'etale);
  \item \'etale locally on $X$, there exists a chart $(P_{X}\rightarrow \sM_{X},Q_{Y}\rightarrow \sM_{Y},Q\rightarrow P)$ extending the chart $Q_{Y}\rightarrow\sM_{Y}$, satisfying the following properties.
    \begin{enumerate}
      \item The kernel and the torsion part of the cokernel (resp. the kernel and the cokernel) of $Q^{gp}\rightarrow P^{gp}$ are finite groups of order invertible on $X$.
      \item The induced morphism from $\underline{X}\rightarrow \underline{Y}\times_{Spec \ \Z[Q]} Spec \ \Z[P]$ is \'etale in the classical sense.
    \end{enumerate}
 \end{enumerate}
\end{thm}

\begin{rem}
 \begin{enumerate}
   \item We can require $Q^{gp}\rightarrow P^{gp}$ in (a) to be injective, and replace the requirement of  $\underline{X}\rightarrow \underline{Y}\times_{Spec\ \Z[Q]} Spec\ \Z[P]$ be \'etale in (b) by requiring it to be smooth without changing the conclusion of the theorem \ref{KatoStrThm}. 
   \item The arrow in (b) shows that a log smooth morphism is ``locally toric'' relative to the base. Consider the case $Y$ is a log scheme with underlying space given by Spec $\C$ with the trivial log structure, and $X=$ Spec$(P\rightarrow \C[P])$ where $P$ is a fine, saturated and torsion free monoid. Then $\underline{X}$ is a toric variety with the action of Spec $\C[P^{gp}]$. According to the theorem, $X$ is log smooth relative to $Y$, though the underlying space might be singular. These singularities are called toric singularities in \cite{ToricSing}. This is closely related to the classical notion of toroidal embeddings, see \cite{KKMS}.
 \end{enumerate}
\end{rem}

\begin{eg} 
Using the theorem, we can check directly that the morphism $f$ in example \ref{NCsmooth} is log smooth, but the underlying map has normal crossing singularities. We will see later that one of the major advantages of log structures is in dealing with the normal crossing singularities.
\end{eg}

Let $X\stackrel{f}{\rightarrow}Y\stackrel{g}{\rightarrow}Z$ be morphisms of fine log schemes. Consider the sheaves of log differentials $\Omega_{g}^{1}$ and $\Omega_{g\circ f}^{1}$, with their universal derivations $(\partial_{g},D_{g})$ and $(\partial_{g\circ f},D_{g\circ f})$ respectively. We have a canonical map $f^{*}\Omega_{g}^{1}\rightarrow \Omega_{g\circ f}^{1}$ induced by 
\[f^{*}(\partial_{g}u)\mapsto \partial_{g\circ f}f^{*}(u) \mbox{\ \ and\ \ }f^{*}(D_{g}v) \mapsto D_{g\circ f}f^{\flat}(v),\]
where $u\in \sO_{Y}$ and $v\in \sM_{Y}$. Denote by $(\partial_{f},D_{f})$ the universal derivation associated to $\Omega_{f}^{1}$. Similarly, we have a canonical map $\Omega_{g\circ f}^{1}\to \Omega_{f}^{1}$ induced by 
\[\partial_{g\circ f}u'\mapsto \partial_{f}u' \mbox{\ \ and\ \ }D_{g\circ f}v' \mapsto D_{f}v',\]
where $u'\in \sO_{X}$ and $v'\in \sM_{X}$. The following proposition shows that log differentials behave like usual differentials, especially for log smooth morphisms.

\begin{prop}\label{prop:logsmoothdif} 
 \begin{enumerate}
  \item The sequence $f^{*}\Omega_{g}^{1} \rightarrow \Omega_{g\circ f}^{1}\rightarrow \Omega_{f}^{1}\rightarrow 0$ is exact.
  \item If $f$ is log smooth, then $\Omega_{f}^{1}$ is a locally free $\sO_{X}$-module, and we have the following exact sequence: $0\rightarrow f^{*}\Omega_{g}^{1} \rightarrow \Omega_{g\circ f}^{1}\rightarrow \Omega_{f}^{1}\rightarrow 0$.
  \item If $g\circ f$ is log smooth and the sequence in (2) is exact and splits locally, then $f$ is log smooth.
 \end{enumerate}
\end{prop}

A proof can be found in \cite[Chapter IV]{Ogus}.

\subsection*{Logarithmic smooth deformation}

Having discussed log smoothness, a natural thing to do is to develop log smooth deformations. In many cases, we would require this to be a flat deformation for the underlying space. Unfortunately log smoothness does not imply flatness, so we need the following definition.

\begin{defn}\label{defn:integralmonmap} 
A map of fine monoids $h:Q\rightarrow P$ is called {\em integral} if the induced map on monoid algebra $\Z[Q]\rightarrow\Z[P]$ is flat.
\end{defn}

\begin{defn}\label{defn:integralmap} 
A morphism $f:X\rightarrow Y$ of fine log schemes is called \emph{integral} if for every geometric point $\bar{x} \in X$, the map of characteristic monoids $h:f^{-1}(\overline{\sM}_{Y})_{\bar{x}}\rightarrow\overline{\sM}_{X,\bar{x}}$ is integral.
\end{defn}

\begin{rem}
 \begin{enumerate}
  \item If $f$ is integral, then \'etale locally we have a chart $(P_{X}\rightarrow\sM_{X},Q_{Y}\rightarrow\sM_{Y},Q\stackrel{h}{\rightarrow}P)$ such that $h$ is integral.
  \item If $h:Q\rightarrow P$ is integral of integral monoids, then for any integral monoid $Q^{'}$, the push-out of $P\leftarrow Q\rightarrow Q^{'}$ in the category of monoids is integral. Thus integral morphisms are stable under base change by integral log schemes.
  \item Given a morphism $h: Q\rightarrow P$ of integral monoids, there is an explicit criterion, which looks complicated, but sometimes is useful for checking integrality of $h$ directly: if $a_{1},a_{2}\in Q$, $b_{1},b_{2}\in P$ and $h(a_{1})b_{1}=h(a_{2})b_{2}$, then there exist $a_{3},a_{4}\in Q$ and $b\in P$ such that $b_{1}=h(a_{3})b,b_{2}=h(a_{4})b$ and $a_{1}a_{3}=a_{2}a_{4}$. This comes essentially from the equational criterion for flatness.  
 \end{enumerate}
\end{rem}

Now we have the following fact from \cite[4.5]{KKato}.

\begin{prop}\label{prop:logsmooth-flat} 
If $f$ is a log smooth and integral morphism of fine log schemes, then $\underline{f}$ the underlying map is flat in the usual sense.
\end{prop}

Now let us consider the following deformation problem. We are given a log smooth integral morphism $f_{0}:X_{0}\rightarrow B_{0}$ of fine log schemes, and a strict closed immersion $j:B_{0}\rightarrow B$ defined by an ideal $J$ with $J^{2}=0$. We want to find a log smooth lifting $f:X \rightarrow B$ fitting in the following cartesian diagram:
\[
\xymatrix{
X_{0} \ar@{^{(}->}[r] \ar[d] & X \ar[d]\\
B_{0} \ar@{^{(}->}[r] & B.
}
\]

\begin{rem}
Since $f_{0}$ is integral, and $\sM_{X}/(1+J) \cong \sM_{X_{0}}$, it is not hard to show that the lifting $f$ is automatically integral and hence flat.
\end{rem}

We have the following theorem for log smooth deformations.

\begin{thm}\cite[3.14]{KKato}\label{thm:basiclogdef}
With the notation as above, we have:
 \begin{enumerate}
  \item There is a canonical obstruction $\eta\in 
H^{2}(\underline{X_0},T_{X_{0}/B_{0}}\otimes J)$ such that $\eta=0$ if and 
only if there exists a log smooth lifting.
  \item If $\eta=0$, then the set of log smooth deformations form a torsor 
under $H^{1}(\underline{X_0},T_{X_{0}/B_{0}}\otimes J)$.
  \item The automorphism group of any deformation is given by 
$H^{0}(\underline{X_0},T_{X_{0}/B_{0}}\otimes J)$.
 \end{enumerate}
\end{thm}

The theorem can be proved in a manner similar to the case of usual deformation theory as in \cite[Expos\'e 3]{SGA_I}. Another proof using the logarithmic cotangent complex can be found in \cite[Thm 5.6]{LogCot}, which we will discuss later.

\section{Log smooth curves and their moduli}\label{Satriano}

Now that we have reviewed some of the foundations, we can discuss the
first application to the study of moduli spaces: F. Kato's interpretation
of $\mgnb$ as a moduli space for log curves.  The general 
philosophy is expressed by F. Kato in the introduction to \cite{FKato}:
\begin{philosophy}
Since log smoothness includes some degenerating objects like semistable reductions, etc., 
the moduli space of log smooth objects should be already compactied, once its existence 
has been established.
\end{philosophy}
Along the lines of this philosophy, to compactify $\mgn$, we want to
introduce a notion of \emph{log curve} which extends  
the notion of smooth curve.  Following F. Kato, we do so after some
preliminaries. 

\subsection*{Relative characteristic sheaves}
Recall from Definition \ref{Def:chara} that the characteristic $\overline{\sM}_X$ of a log scheme $X$ is defined as 
$\sM_X/\sO^*_{\underline{X}}$.  In the study of log curves, the following relative notion of characteristic plays 
an important role.
\begin{definition}
\label{def:char}
Given a morphism $f:X\rightarrow Y$ of log schemes, the \emph{relative characteristic} $\overline{\sM}_{X/Y}$ is 
defined as the quotient $\sM_X/\im(f^*\sM_Y\rightarrow\sM_X)$ in the category of integral monoids.
\end{definition}
\begin{example}
\label{ex:char}
Let $f:X\rightarrow Y$ be the morphism from Example \ref{NCsmooth}.  Then the relative characteristic 
$\overline{\sM}_{X/Y}$ is the cokernel in the category of integral monoids of the diagonal map 
$\Delta:\N\rightarrow \N^2$, which is $\Z$.
\end{example}
\begin{lemma}[{\cite[Lemma 1.6]{FKato}}]
If $f:X\rightarrow Y$ is an integral morphism of fine log schemes, then $\overline{\sM}_{X/Y,\bar{x}}=0$ 
if and only if $f$ is strict in an \'etale neighborhood of $x$.
\end{lemma}
As the following example illustrates, the integrality assumption on $f$ is necessary.
\begin{example}
\label{ex:intrel}
Let $P$ be the monoid on three generators $x$, $y$, and $z$ subject to the relation $x+y=2z$.  We have an injection
\[
i:P\longrightarrow \N^2
\]
sending $x$ to $(2,0)$, $y$ to $(1,1)$, and $z$ to $(0,2)$.  Let $X=\Spec k[\N^2]$ and $Y=\Spec k[P]$ with their canonical 
log structures.  Then it follows from \cite[Prop 3.4]{KKato} that the morphism of log schemes $f:X\rightarrow Y$ 
induced by $i$ is log \'etale, but $\underline{f}$ is not flat, and hence, by \cite[Cor 4.5]{KKato}, $f$ is not 
integral.  It is easy to check that
\[
\sM_{X/Y,0}=\N^2/P\simeq\Z/2,
\]
so $\overline{\sM}_{X/Y,0}=0$, but $f$ is not strict.
\end{example}

\subsection*{Log curves}
\begin{definition}
\label{def:logcurve}
A \emph{log curve} is a log smooth integral morphism $f:X\rightarrow S$ of fs log schemes such that the 
geometric fibers of $\underline{f}$ 
are reduced connected $1$-dimensional schemes.
\end{definition}
We require $f$ to be integral so that, by \cite[Cor 4.5]{KKato}, $\underline{f}$ is flat.  The reason for the fs 
assumption is to avoid cusps or worse singularities, as the following example shows.
\begin{example}
\label{ex:cusp}
If $X=\Spec k[\N-\{1\}]$ is given its canonical log structure and $S=\Spec k$ is given the trivial log structure, then 
$X\rightarrow S$ is log smooth and integral; however, $\underline{X}=\Spec k[x,y]/(y^2-x^3)$ has a cusp.
\end{example}
It is a remarkable fact that by endowing our curves with log structures as in Definition \ref{def:logcurve}, this is 
enough to control the singularities of the curve.
\begin{theorem}[{\cite[Thm 1.3]{FKato}}]
\label{thm:node}
If $k$ is a seperably closed field and $f:X\rightarrow S$ is a log curve with $\underline{S}=\Spec k$, then 
$\underline{X}$ has at worst nodal singularities.  Moreover, if $r_1,\dots,r_\ell$ are the nodes of $\underline{X}$, 
then there exist smooth points $s_1,\dots,s_n$ of $\underline{X}$ such that
\[
\overline{\sM}_{X/S}=\Z_{r_1}\oplus\dots\oplus\Z_{r_\ell}\oplus\N_{s_1}\oplus\dots\N_{s_n};
\]
here $M_x$ denotes the skyscraper sheaf for a monoid $M$ supported at a point $x\in \underline{X}$.
\end{theorem}
The reader should think of the $s_i$ in the above theorem as marked points.  So we can already see how $n$-pointed curves 
emerge naturally from the log geometry perspective.
\begin{example}
\label{ex:logcurve}
Consider the closed subscheme $\underline{X}$ of $\mathbb{P}^2_k\times_{k}\mathbb{A}^1_k$ defined by $xz=ty$, where $t$ is the 
coordinate of $\mathbb{A}^1_k$ and $x$, $y$, and $z$ are the coordinates of $\mathbb{P}^2_k$.  Then $\underline{X}$ has a natural 
log structure $\sM_X$.  For example, on the locus where $z$ is invertible, $\underline{X}$ is given by $\Spec k[P_z]$ with $P_z$ 
a monoid on five generators $a$, $b$, $c$, $c'$, $u$ subject to the relations $c+c'=0$ and 
$a+c=b+u$; here $\sM_X$ is given by the canonical log structure associated to $P_z$.  Then the projection
\[
X\longrightarrow \mathbb{A}^1_k
\]
is a log curve, where $X$ is given the log structure above and $\mathbb{A}^1_k$ 
is given the log structure defined by the divisor $t=0$.  We see that every fiber above $t\neq0$ is 
isomorphic to $\mathbb{P}^1_k$ with log structure given by the divisor at 0 and $\infty$; the fiber above $t=0$ is nodal.  
The $n$ in Theorem \ref{thm:node} is equal to 2 for all geometric fibers.
\end{example}
Since our goal is to give a log geometric 
description of $\mgnb$, we would like to express the stability condition purely in terms of log geometry.  The following 
proposition provides the key.
\begin{proposition}[{\cite[Prop 1.13]{FKato}}]
\label{prop:dualizing}
With notation as in Theorem \ref{thm:node}, there is a natural isomorphism
\[
\Omega^1_{X/S}\longrightarrow\omega_{\underline{X}}(s_1+\dots+s_n),
\]
where $\omega_{\underline{X}}$ is the dualizing sheaf of $\underline{X}$.
\end{proposition}
We therefore make the following definition.
\begin{definition}
\label{def:stable}
Let $f:X\rightarrow S$ be a log curve and for all geometric points $\bar{t}$ of $\underline{S}$, let $\ell(\bar{t})$ 
and $n(\bar{t})$ be such that
\[
\overline{\sM}_{X_{\bar{t}}/\bar{t}}=\Z_{r_1}\oplus\dots\oplus\Z_{r_{\ell(\bar{t})}}\oplus\N_{s_1}\oplus\dots\N_{s_{n(\bar{t})}}.
\]
We say $f$ is \emph{of type $(g,n)$} if $\underline{f}$ is proper, $\underline{X}$ has genus $g$, and $n(\bar{t})=n$ for all $\bar{t}$.  We say 
$f$ is \emph{stable of type of $(g,n)$} if it of type $(g,n)$ and 
\[
H^0(\underline{X}_{\bar{t}},T_{X_{\bar{t}}/\bar{t}})=0
\]
for all geometric points $\bar{t}$ of $\underline{S}$.
\end{definition}
It is, in fact, true (\cite[Prop 1.7]{FKato}) that if $f:X\rightarrow S$ is a log curve of type $(g,n)$, 
then the $s_i$ in each geometric fiber fit together to yield $n$ sections $\sigma_i$ of $\underline{f}$.  It follows 
then that every stable log curve of type $(g,n)$ is an $n$-pointed stable curve of genus $g$ in the classical sense.

\subsection*{Log structures on stable curves}
Having now shown that every log curve is naturally a pointed nodal curve, we shift gears and ask the following 
question: given a stable genus $g$ curve $\underline{f}:\underline{X}\rightarrow\underline{S}$ with $n$ marked points, how 
many log structures can we put on $\underline{X}$ and $\underline{S}$ so that the associated morphism of log schemes 
is a log curve with relative characteristic supported on our given $n$ marked points?  We begin by sketching the 
construction of a canonical such log structure.

Lemmas 2.1 and 2.2 of \cite{FKato} show that to endow $\underline{X}$ and $\underline{S}$ with log structures as 
desired, it is enough to consider the case when $\underline{S}=\Spec A$ and $A$ is strict Henselian.  For every node 
$r_i$ of the closed fiber of $\underline{f}$, we can find an \'etale neighborhood $U_i$ of the points 
specializing to $r_i$ and a diagram
\[
\xymatrix{
U_i\ar[d]\ar[r]^-{\psi_i} & \Spec A[x,y,t]/(xy-t)\ar[d]^{\pi}\\
S\ar[r]^-{\varphi_i} & \Spec A[t]
}
\]
which is cartesian.  Let $t_i\in A$ be the image of $t$ under the morphism induced by $\varphi_i$.  
Endowing $\Spec A[t]$ with the log structure associated to the morphism $\N\rightarrow A[t]$ sending 1 to $t$, and 
$\Spec A[x,y,t]/(xy-t)$ with the log structure associated to $\N^2\rightarrow A[x,y,t]/(xy-t)$ sending $e_1$ (resp. $e_2$) 
to $x$ (resp. $y$), we see that $(\pi,\Delta)$ 
is a morphism of log schemes, where $\Delta:\N\rightarrow\N^2$ is the diagonal map.  Pulling back these log structures 
under $\varphi_i$ (resp. $\psi_i$), we obtain log structures $\mathcal{L}_i$ (resp. $\sM'_i$) on $S$ (resp. $U_i$).  
Away from the points specializing to $r_i$, we define a log structure $\sM''_i$ as the pullback of $\mathcal{L}_i$.  
The log structures $\sM'_i$ and $\sM''_i$ glue to yield a log structure $\sM_i$ on $X$.  Let $\mathcal{N}$ be the 
log structure on $X$ associated to the divisor defined by the marked points.  We let
\[
\sM_X=\sM_1\oplus_{\sO_X^*}\dots\oplus_{\sO_X^*}\sM_{\ell}\oplus_{\sO_X^*}\mathcal{N}
\]
and 
\[
\sM_S=\mathcal{L}_1\oplus_{\sO_S^*}\dots\oplus_{\sO_S^*}\mathcal{L}_{\ell}.
\]
It is not difficult to see that with these definitions, we have endowed $f$ with the structure of a log curve.

Moreover, a detailed analysis of the proof of Theorem \ref{thm:node} shows that this log structure we have just 
constructed is ``minimal'' among all possible log structures giving $\underline{X}/\underline{S}$ the structure of 
a log curve (see 1.8 and Thm 2.3 of \cite{FKato}):
\begin{theorem}
\label{thm:basic}
Let $\underline{X}/\underline{S}$ be a stable genus $g$ curve with $n$ marked points and let $X/S$ be the log curve 
obtained by endowing $\underline{X}/\underline{S}$ with the canonical log structure above.  If $X'/S'$ is a log curve 
and $\underline{a}:\underline{S}'\rightarrow\underline{S}$ and 
$\underline{b}:\underline{X}'\rightarrow\underline{X}$ are morphisms 
such that $\underline{X}'\simeq\underline{X}\times_{\underline{S}}\underline{S}'$ and such that the divisors of marked 
points in $\underline{X}'$ are sent scheme-theoretically to the divisors of marked points in $\underline{X}$, 
there are unique morphisms $a$ and $b$ of log schemes extending the morphisms $\underline{a}$ and $\underline{b}$ 
above such that 
\[
\xymatrix{
X'\ar[r]^-{b}\ar[d] & X\ar[d]\\
S'\ar[r]^-{a} & S
}
\]
is cartesian in the category of fs log schemes.
\end{theorem}
\begin{definition}
\label{def:basic}
A log curve $X/S$ is called \emph{basic} if it satisfies the universal property in Theorem \ref{thm:basic}.
\end{definition}

\subsection*{Moduli}
Before discussing moduli of log smooth curves, we begin with some generalities about log structures and stacks.  
Note that the definition of a log structure is not particular to schemes; indeed, Definition \ref{Def:log-str} 
makes sense for any ringed topos.  We can therefore define log structures on the \'etale site of a Deligne-Mumford 
stack or the lisse-\'etale site of an Artin stack.  The notions of fine and fs log structures carry over to this 
setting as well, so one can speak of fine (or fs) log algebraic stacks.

There is another equivalent way to talk about log structures on stacks; namely, if $\mf{X}$ is a stack over the 
category of schemes, and $\sM_{\mf{X}}$ is a log structure on $\mf{X}$, then $\mf{X}$ can naturally be viewed as 
a stack over the category of log schemes.  For concreteness, say $\mf{X}$ is a stack over the 
category of schemes with the \'etale topology.  Then we obtain a category $\tilde{\mf{X}}$ fibered over the category 
of log schemes by defining $\tilde{\mf{X}}(S,\sM_S)$ to be the category whose objects are morphisms 
\[
f:(S,\sM_S)\rightarrow(\mf{X},\sM_{\mf{X}})
\]
of log stacks and whose morphisms $g$ from $f\in\tilde{\mf{X}}(S,\sM_S)$ to $f'\in\tilde{\mf{X}}(S',\sM_{S'})$ are 
given by diagrams
\[
\xymatrix{
(S,\sM_S)\ar[rr]^g\ar[dr]_f & & (S',\sM_{S'})\ar[dl]^{f'}\\
 & (\mf{X},\sM_{\mf{X}}). &
}
\]
One checks that $\tilde{\mf{X}}$ is a stack over the category of log schemes where coverings are given by surjective strict 
\'etale morphisms.

Conversely, given any stack $\mathcal{Y}$ over the category of log schemes with the strict \'etale topology, we obtain 
a log stack $(\mathcal{Y}',\sM_{\mathcal{Y}'})$ over the category of schemes with the \'etale topology by letting 
$\mathcal{Y}'(S)$ be the category of pairs $(\sM_S,\xi)$ where $\xi$ is an object of $\mathcal{Y}(S,\sM_S)$.  The 
log structure $\sM_{\mathcal{Y}'}$ is then defined by the following property: if $f:S\rightarrow\mathcal{Y}$ is 
a morphism which corresponds to the pair $(\sM_S,\xi)$, then $f^*\sM_{\mathcal{Y}'}=\sM_S$.

It is however important to note that these {\em two procedures are not
inverse to each other}. If we start with a logarithmic stack
$(\mf{X},\cM_{\mf{X}})$, and take the result  $(\tilde{\mf{X}})'$ of
the composite operation, we do not get $\mf{X}$ but rather the stack
$LOG_{\mf{X}}$ described in section \ref{LogStacks}. In order to
recover the stack $\mf{X}$ over schemes from a stack $\tilde{\mf{X}}$
over log schemes, it is necessary to distinguish objects similar to
the basic log curves of Definition \ref{def:basic}. Rather than launch
into a premature categorical discussion, let us see how this works for
log curves. The issue is revisited in sections \ref{Roots} (with
Olsson's terminology of {\em special log structures})
and \ref{Stablemaps} (with Kim's terminology of {\em minimal log
structures}).

Let $\lmgnb$ be the stack over the category of fs log schemes with the strict \'etale topology where 
$\lmgnb(S,\sM_S)$ is the category of stable log curves of type $(g,n)$ over $(S,\sM_S)$.  
Let $\mgn^{bas}$ be the substack of basic stable log curves of type $(g,n)$.  
By the above discussion, we obtain a log stack over the category of schemes with the \'etale topology, which 
we again denote by $(\mgn^{bas},\sM_{\mgn^{bas}})$.  

Note that the Deligne-Mumford compactification $\mgnb$ carries a natural log structure 
$\sM_{\mgnb}$ coming from the simple normal crossing divisor at the boundary.  It follows from \cite[\S2]{DM} that 
$\sM_{\mgnb}$ can also be described as the log structure which assigns to each stable curve 
$\underline{X}/\underline{S}$ the basic log structure obtained on $\underline{S}$.  
The discussion following Definition  \ref{def:logcurve} above shows that we have a natural morphism 
\[
F:(\mgn^{bas},\sM_{\mgn^{bas}})\longrightarrow(\mgnb,\sM_{\mgnb})
\]
and from Theorem \ref{thm:basic} we see 
\begin{theorem}[{\cite[\S4]{FKato}}]
The morphism $F$ is an equivalence of log stacks.
\end{theorem}

\begin{anitem} {\bf Back to the big picture}
\label{Curves-big-picture}
\end{anitem}
We end by mentioning a type of converse to the philosophic principal mentioned at the beginning of this 
section.  We have seen that since log smoothness includes degenerate objects, log geometry can naturally lead to 
compactifications; however, it is also generally true that we do not end up with ``too many'' degenerate objects.
\begin{philosophy}
\emph{
Log geometry controls degenerations.
}
\end{philosophy}
In higher dimensions, compactifications tend to have unwanted extra components.  Log geometry helps 
to cut down on these components.  Let us give some inkling of an idea as to why this should be true.  Suppose $\mf{X}$ is 
an algebraic stack which is irreducible.  Suppose we can find a proper algebraic stack $\bar{\mf{X}}$ with a fine log structure 
and an open immersion $i:\mf{X}\rightarrow\bar{\mf{X}}$ such that $\mf{X}$ is the trivial locus of $\bar{\mf{X}}$.

As we now explain, log geometry provides us with a good method of trying to show that $\bar{\mf{X}}$ is irreducible 
as well.  If $k$ is separably closed and $x:\Spec k\rightarrow \bar{\mf{X}}$ is a morphism, then pulling back the log 
structure on $\bar{\mf{X}}$ endows $\Spec k$ with a fine log structure.  By \cite[Lemma 2.10]{KKato}, it follows that this 
is the log structure associated to a morphism of monoids $P\rightarrow k$, where $P$ is fine.  Hence, 
we have a strict closed immersion of log schemes $j:\Spec k\rightarrow \Spec k[[P]]$, where $\Spec k[[P]]$ is given its 
canonical log structure.  Note that the generic point $\Spec K$ of $\Spec k[[P]]$ carries the trivial log structure.  
Therefore, if $x$ factors as a morphism of log stacks through $j$, then we automatically obtain a commutative 
diagram
\[
\xymatrix{
\mf{X}\ar[rr]^i & & \bar{\mf{X}}\\
\Spec K\ar[r]\ar[u]^y & \Spec k[[P]]\ar[ur] & \Spec k\ar[u]_x\ar[l]^-{j}
}
\]
and hence $x$ is the specialization of the point $y$ of $\mf{X}$.  We see then that the log structure on $\Spec k$ obtained 
from $\bar{\mf{X}}$ somehow serves as a compass telling us which way to look in order to find a family degenerating to 
our given point $x$ of $\bar{\mf{X}}$.

\section{D-semistability and log structures}\label{D-SS}

The main references here are \cite{Friedman,FKato-def,Olsson2003}.

\bg{conv}
Throughout this section, every scheme is over an algebraically closed
 field. The notation $X$ will be reserved for a {\em normal crossing variety}.  
 By this we mean a variety for which every closed point $x\in X$ has
 an \'etale neighborhood $x \to U$ which admits an \'etale map
$U \to \Spec\dfrac{k[x_1,\cdots,x_n]}{(x_1\cdots x_r)}\ra X$ such that
 $x$ maps to the  point with coordinates  
$(0,\cdots,0)$. By a \std\nbhd of $x\in X$, we mean an \'etale \nbhd
 of $x$ as above. The notation $\sM_k$ denotes the log \str 
of the \std log \pt $\NN\ra k, n\mapsto 0^n$, see example \ref{Logpt}.
\ed{conv}

\subsection*{Introduction}

To study the geometry of a normal crossing variety $X$, e.g. 
to study the deformation theory of $X$, one would like to ask the
following questions: 

\begin{question}\label{1} {Can we embed $X$ into another variety $i:X\ra\scr{X}$
as a normal crossing divisor?}   
\end{question}
\begin{question}
\label{2}{Can we find a semi-stable smoothing of $X$, i.e. embed
$X\ra\Spec{k}$ in a flat family over a curve $\scr{X}\ra C$,  
in such a way that  there exists a diagram 
\[
{
   \xymatrix{
  X \ar[d] \ar[r]          & \scr{X} \ar[d]  & \ar[l] \scr{X}^*=\scr{X}\backslash X\ar[d]_{f^*}\\
  \Spec{k}  \ar[r]^{0}         & C               & \ar[l] C^*=C\backslash{0}
  }
}
\]\newline
 where $\scr{X}$ is smooth, the squares are cartesian and $f^*$ is 
smooth?}
\end{question}

\zz
The answer to these questions are not always yes, because the
existence of such maps would imply the existence of certain log
structures on $X$, which in turn would imply intrinsic condition on $X$,
 so their existence is not guaranteed.

\bg{eg}
[log structure of embedding type]\label{emb}
If we can find an embedding as in Question (\ref{1}) above, then $i^*(j_* \mathcal{O}^{\times}
_{\scr{X}^*})\ra\ox$ defines a log \str $\sM_X$ on $X$, which \etale locally has a chart $\mathbb{N}^r\ra\Spec\dfrac{k[x_1,\cdots,x_n]}{(x_1\cdots
x_r)}$ sending the \elt $e_i$ in the \std basis of $\NN^r$ to $x_i$.
This is called a {\em log \str of embedding type.} 
\end{eg}


\bg{eg}[log \str of semi-stable type]\label{sm} If we can find a
semi-stable smoothing as in Question (\ref{2}) above, then 
what we have in this case is not only a log \str of embedding type on
$X$, but also a morphism (of sheaf of monoids) $f^\flat:
f^*\sM_k\ra \sM_X$, which makes $X$ a log smooth variety over the
standard log point ($\Spec
k, \sM_k$). \'Etale locally a chart for the log \str on X can be put
in the form $$\big(\Spec\dfrac{k[x_1,\cdots,x_n]}{(x_1\cdots 
x_r)}, \NN^r, e_i\mapsto x_i, i=1,\cdots,r\big).$$ Modulo the units in
the monoid, the morphsim of quotient monoids induced by
$f^\flat:f^*\sM_k\ra \sM_X$ is just the diagonal
$\Delta:\NN\ra\NN^r$. Such a pair $(\sM_X,f^\flat:f^*\sM_k\ra \sM_X)$
is called a {\em log \str of semi-stable type on $X$.} 
\end{eg}

\bg{rmk}
The log structure $\sM_k$ on $\Spec k$ can be defined by the pullback of the log \str $(\sM_0=j_* \mathcal{O}^{\times}_{C^*}\ra\oo_C)$ on $C$, i.e. the log \str defined by the divisor $0$ of $C$. We have an \isom $\overline{\sM_k}:=\sM_k/k^\times\iso (\sM_0/\oo_C^\times)_0\iso\NN$, where the second \isom assigns each function to its vanishing order at the point $0$ in $\NN$. This gives a geometric interpretation of the \std log point.
\end{rmk}

Concerning the existence of such log \str on the normal crossing variety $X$, we have the following theorems (\cite{FKato-def}, Sec.11):

\begin{thm}\label{embed}
Let X be a normal crossing variety over the spectrum of
an algebraically closed field, then X can be equipped with
a log structure of embedding type iff there exists a line bundle $\mathcal{L}$ on $X$ such that
\begin{equation*}
\mathcal{E}xt_X^1(\Omega^1_X,\mathcal{O}_X)|_D\iso\mathcal{L}|_D
\end{equation*}
where D is the non-smooth locus of X.
\end{thm}

\bg{rmk}
It is not hard to see $\mathcal{E}xt_X^1(\Omega^1_X,\mathcal{O}_X)|_D$ is a line bundle on $D$.
\ed{rmk}

\bg{defn} (see \cite[Def. 1.13 and Prop. 2.3]{Friedman})
Let $X$ and $D$ as before. If
$\mathcal{E}xt_X^1(\Omega^1_X,\mathcal{O}_X)|_D$ is a trivial on $D$,
then we say that $X$ is \emph{d-semistable}.  
\ed{defn}

\begin{thm}[d-semistability]\label{d-semi}
Let $X$ be a normal crossing variety over the spectrum of
an algebraically closed field, then $X$ can be equipped with
a log structure over the standard log point, such that the structure
morphism is log smooth if and only if $X$ is d-semistable.
\end{thm}

Generalization of these theorems can be found in
\cite{Olsson2003}, section 3.

\bg{cor}
If $X$ has a semistable smoothing as in Question \ref{2}, then $X$ is d-semistable.
\ed{cor}
In fact, we can put a log structure on $X$ such that it is log smooth over the standard log point by Example \ref{sm},
 so we can apply Theorem \ref{d-semi}

\begin{rmk}
Being d-semistable is not equivalent to having a semi-stably
smoothing. See \cite[Section 3]{MR705035} for counterexamples.
\end{rmk}

\begin{eg}[a normal crossing variety that is not
d-semistable \cite{Friedman}]\label{tetra} 
Let $X$ be the subvariety of $\mathbf{P}^3$ defined by the product
of 4 linear equations $f=L_1L_2L_3L_4=0$. It is a normal crossing variety
provided the four planes has no points in common. Then $D$ is
defined by the homogeneous ideal $(L_iL_j|1\leq i<j\leq 4)$, and it
is not hard to calculate that
$\mathcal{E}xt_X^1(\Omega^1_X,\mathcal{O}_X)|_D\cong
\mathcal{O}_D(4)$, which is not trivial. So $X$ is not d-semistable.
\end{eg}

\bg{cor}
The $X$ in Example \ref{tetra} does not admit a  semi-stably smoothing.
\ed{cor}

\begin{rmk}
If we put $X$ in the $1-$dimensional family $\scr{X}$ defined by
 $f+tg=0$ with parameter $t$, 
 where $g$ is a smooth quartic, then $X$ is the fiber over $t=0$ and for generic $g$, $\scr{X}$
 over $(t\neq 0)$ is smooth. But the whole space of this family is not smooth: in fact, 
for a generic $g$, this family is singular at the 24 points of $D\cap\{g=0\}$. However, 
a single blowing up at such a point gives a $\mathbf{P}^1\times\mathbf{P}^1$. Contract 
along either ruling gives back a family with parameter $t$ which has one less singularity. 
If we do this process to all 24 points, we will get a family which is a semi-stable smoothing 
of $\tilde{X}$, the blowing-up of $X$ at those $24$ points. $\tilde{X}$ is d-semistable. 
For details, see \cite[Rem. 1.14]{Friedman}.
\end{rmk}

\subsection*{Refined analysis of the existence of log structures}
Continuing with the notation $X,D$ above, let us analyze the
situation. We want to break the job of finding a  
suitable log structure over the standard log point into 2 steps:

\begin{enumerate}
\item Put a log structure of embedding type on $X$.
\item See if it is possible to make that log structure semi-stable.
\end{enumerate}

We will see that two related obstructions arise naturally, where the
vanishing of the first corresponds to the first step, and the
vanishing of the second, which means precisely being d-semistable,
allows us to do the second step.

Since \'etale locally a log \str of embedding type always exists, let us 
consider the stack $\mcal{G}$, which to each $U\in \xet$ associates the 
\gpd of log \str of embedding type on $U$. 
Using Artin's approximation theorem one can show that any two elements of $U$ are 
locally 
isomorphic, which means $\mcal{G}$ is a gerbe. Since
$\Aut{\mcal{G}}\iso\mcal{K}$, where $\mcal{K}$ is the kernel of the
restriction map $\ox^\times\ra\oo_D^\times$, we have: 

\bg{prop}
There is an obstruction $\eta$ in $H^2(\xet,\mcal{K})$ 
whose vanishing is equivalent to the existence of a log \str of embedding type on $X$.
\ed{prop}

For the calculation of this obstruction, we state the following result (See \cite{FKato-def}, Sec.11
and \cite{Olsson2003}, Sec.3):

\bg{prop}\label{gerbe}
In the \les of cohomology associated to the \ses 
$1\ra\mcal{K}\ra\ox^\times\ra\oo_D^\times\ra1$, the line bundle $\mathcal{E}xt_X^1(\Omega^1_X,\mathcal{O}_X)|_D$ maps to $-\eta\in H^2(\xet,\mcal{K}).$ 
\end{prop}

Combining these results with the exactness of the \les, we get \Thm\ref{embed}.

\bg{rmk}
\label{gr}
General theory tells us if $\eta=0$, then the set of all log \str of embedding type on $X$ is naturally a torsor under $H^1(\xet,\mcal K)$. In general the set of all log \str of embedding type is only a pseudo torsor under this group.
\ed{rmk}

Suppose $\eta=0$, then we can put a log structure of embedding type on $X$, which maps to $(\Spec k,k^\times)$. To go one step further, i.e., to make it a log \str of semistable type, we need a morphism of monoids $f^*\sM_k\ra \sM_X$, such that in an \std \nbhd of $x\in X$, a chart of the morphism is given by the diagonal $\Delta:\NN\ra\NN^r$, where $r$ is the number of irrducible components passing through $x$.
 
Since $\sM_k\iso\NN\oplus k^\times$ (non-canonically), and the image of $k^\times\subset \sM_k$ is determined by the underlying morphism of schemes. To give the morphism wanted from $\sM_k$ to $\sM_X$, we only have to specify the image of an \elt in $\sM_k$ having vanishing order 1.

Now the question becomes a lifting problem for morphism of sheaves in monoids:
\[
   \xymatrix{
                               & \sM_X \ar[d]^{\beta}  \\
  {\NN}\iso f^{-1}\overline{\sM_k} \ar@{..>}[ur] \ar[r]^(.65){ \overline{f^\flat}}         & \overline{\sM_X}   
  }
\]
where \'etale locally $\overline{f^\flat}$ is the diagonal
$\Delta$. To lift $\overline{f^\flat}$ it is equivalent to lift the
element $\overline{f^\flat}(1)$. 

Consider the sheaf of all the possible local liftings of $\overline{f^\flat}(1)$,
$T=\beta^{-1}(\overline{f^\flat}(1))$, then $T$ is a torsor under $\ox^\times$. To
find a lifting of $\overline{f^\flat}(1)$, it is equivalent to find a global section
of $T$, i.e. a trivialization of $T$.

It seems like we have got an obstruction of finding a log \str on $X$ 
which is semi-stable in $H^1(\xet,\ox)=\Pic(X)$. This is, however, not 
quite true. In fact, what we got is for each log scheme $X$ with a log 
\str of embedding type, the obstruction of making it semi-stable. And our 
original question (on the existence of log \str of semi-stable type) 
allows some ambiguity of choosing the log \str of embedding type $\sM_X$ 
on $X$. As we said in remark \ref{gr}, in this case the set of all log 
\str of embedding type on $X$ is an $H^1(\xet,\mcal K)$-torsor, which 
implies:

\bg{prop}
If $\eta=0$, i.e. there exists a log \str of embedding type on $X$. Then there is an obstruction for finding a log \str of semi-stable type on $X$, $\eta'\in H^1(\xet,\ox^\times)/H^1(\xet,\mcal K)$, whose vanishing is equivalent to the existence of such a log structure.  
\end{prop}

By the \les of cohomology, $H^1(\xet,\ox^\times)/H^1(\xet,\mcal K)$ embeds into $H^1(D_{\acute et},\oo_D^\times)\iso\Pic(D)$. For the calculation of $\eta'$ as an \elt of $\Pic(D)$, we state the following proposition (See \cite{Olsson2003}).

\bg{prop}
We have $-\eta'=[\mathcal{E}xt_X^1(\Omega^1_X,\mathcal{O}_X)|_D]\in \Pic(D)$.
\ed{prop}

Combining these two proposition, we get \Thm \ref{d-semi}.

\section{Stacks of logarithmic structures}\label{LogStacks}
The main reference of this section is \cite{LogStack}.

\subsection*{A motivating example}
Before introducing the stack $LOG_S$ classifying fine log 
structures on schemes over a fine log scheme $S,$ constructed 
by Olsson in \cite{LogStack}, let us look at an example, which will give 
the local covers of $LOG_S.$

\begin{definition}\label{DF-log-str}
Let $\underline{X}$ be a scheme and $r\ge1$ an integer. A 
\emph{Deligne-Faltings log 
structure of rank $r$ on $\underline{X}$} (abbreviated as a 
\emph{DF log structure of rank $r$}) is the following date:

$\bullet$ a sequence $L_1,\cdots,L_r$ of line bundles on 
$\underline{X},$ and 

$\bullet$ a morphism $s_i:L_i\rightarrow\sO_X$ of line bundles, 
for each $i.$ 
\end{definition}

Consider the following three categories fibered in groupoids 
over the category of schemes:

\begin{enumerate}
\item \quad the category of triples $(\underline{X}, 
L,s:L\to\sO_{\underline{X}})$ consisting of a scheme 
$\underline{X}$ and a DF log structure of rank 1 on 
$\underline{X};$

\item \quad the category of pairs $(X,\beta:\mathbb 
N\to\overline{\sM}_X)$ consisting of a fine log scheme $X$ and a 
morphism of sheaves of monoids $\beta$ that \'etale locally 
lifts to a chart: $\widetilde{\beta}:\mathbb N\to\sM_X;$

\item \quad the quotient stack $[\mathbb A^1/\mathbb G_m],$ 
where the quotient is formed with respect to the multiplication 
action of $\mathbb G_m$ on $\mathbb A^1.$
\end{enumerate}

\begin{lemma}(\cite{KKato}, complement 1). 
These three categories fibered in groupoids are equivalent.
\end{lemma}

Let us sketch the proof. Given a DF log structure 
$(L,s:L\to\sO_{\underline{X}})$ of rank 1 on $\underline{X},$ 
define a sheaf of monoids $\sM'$ on $\underline{X}$ to be 
$$
\coprod_{n\ge0}\underline{\text{Isom}}(\sO_{\underline 
{X}},L^{\otimes n}),
$$
the sheafification of the presheaf that takes $U$ to 
$\coprod_{n\ge0}\text{Isom}(\sO_U,(L|_U)^{\otimes n}).$ 
It comes with a natural morphism of sheaves of monoids 
$\sM'\to\mathbb N,$ where the monoid structure on $\sM'$ is induced by 
$$
(n,a:\sO\to L^{\otimes n})\cdot(m,b:\sO\to 
L^{\otimes m})=(n+m,a\otimes b).
$$
The map $s:L\to\sO_{\underline{X}}$ induces a morphism 
$$
\underline{\text{Isom}}(\sO_{\underline{X}},L^{\otimes 
n})\overset{\otimes s}{\longrightarrow}\underline{\text{Hom}} 
(\sO_{\underline{X}},\sO_{\underline{X}})=\sO_{\underline{X}}
$$
of sheaves, hence giving a pre-log structure on $\sM':$
$$
\sM'\to\sO_{\underline{X}}.
$$
We take $\sM_X$ to be the log structure associated to this 
pre-log structure $\sM'.$ Note that $\sM'/\sO_X^*\cong\mathbb 
N,$ and we define $\beta:\mathbb N\to\overline{\sM}_X$ to be the 
composite
$$
\beta:\mathbb N\cong \sM'/\sO_X^*\to \sM_X/\sO_X^*.
$$
Locally the line bundle $L$ is trivial, and one can choose a 
trivialization of $L,$ which gives trivializations of all 
$L^{\otimes n}.$ Sending $n\in\mathbb N$ to this 
trivialization defines a section $\mathbb N\to \sM',$ and hence 
a section $\widetilde{\beta}:\mathbb N\to \sM'\to \sM_X.$ One can 
check that this is a chart. 
  
Conversely, given a fine log structure $(\sM_X,\sM_X\overset{ 
\alpha}{\to}\sO_X)$ on $\underline{X}$ with a morphism  
$\beta:\mathbb N\to\overline{\sM}_X$ that \'etale locally lifts 
to a chart $\widetilde{\beta}:\mathbb N\to \sM_X,$ we have a 
section $\beta(1)$ of $\overline{\sM}_X,$ and its inverse image 
under $\pi:\sM_X\to\overline{\sM}_X$ is an $\sO_X^*$-torsor, which 
corresponds to a line bundle $L.$ The composition 
$$
\pi^{-1}(\beta(1))\subset \sM_X\overset{\alpha}{\longrightarrow} 
\sO_X
$$
gives a morphism of line bundles $s:L\to\sO_X.$ 

Giving a morphism $\underline{X}\to[\mathbb A^1/\mathbb G_m]$ is 
equivalent to giving a $\mathbb G_m$-torsor (namely a line 
bundle $L$) with a $\mathbb G_m$-equivariant morphism to 
$\mathbb A^1:$
$$
\xymatrix@C=.7cm{
Y \ar[r] \ar[d] & \mathbb A^1 \ar[d] \\ 
X \ar[r] & [\mathbb A^1/\mathbb G_m].}
$$
This diagram is equivalent to the following one
$$
\xymatrix@C=.7cm{
Y \ar[r]^-s \ar[d] & \mathbb A^1_X \ar[d] \\
X \ar[r] & [\mathbb A^1/\mathbb G_m]_X,} 
$$
and the top arrow is $\mathbb G_m$-equivariant, 
namely a morphism of line bundles $s:L\to\sO_X.$ This finishes the proof.

In fact, in the three fibered categories, one can replace $\mathbb N$ by 
$\mathbb N^r,$ and rank 1 DF-log structure by rank $r$ DF-log 
structure, and replace $[\mathbb A^1/\mathbb G_m]$ by 
$[\mathbb A^r/\mathbb G_m^r]$ (which is equivalent to $[\mathbb 
A^1/\mathbb G_m]^r$), and they are still equivalent. 

More generally, let $P$ be a fine monoid and $S$ a scheme, 
and let $S[P]$ be the product $S\times_{\text{Spec }\mathbb 
Z}\text{Spec }\mathbb Z[P],$ which has a fine log structure 
coming from the chart $P\to\mathbb Z[P]$ (\ref{Afflog}). For an affine 
$S$-scheme $\text{Spec }R,$ the set of $R$-points $S[P](R)$ is the set of 
monoid homomorphisms $Hom_{\text{mon}}(P,R),$ where $R$ is
regarded as a multiplicative monoid. Let $P^{\text{gp}}$ be the group 
associated to $P.$ For any affine $S$-scheme $\text{Spec }R,$ the group 
$Hom_{\text{mon}}(P^{\text{gp}},R)=Hom_{\text{gp}}(P^{\text{gp}},R^*)$ 
acts on the set $Hom_{\text{mon}}(P,R)$ by pointwise multiplication. This 
induces an action of the $S$-group scheme $S[P^{\text{gp}}]$ on $S[P].$
When $S=\text{Spec }k$ for a field $k$ and $P$ is saturated and 
torsion-free, the $k$-group variety $S[P^{\text{gp}}]$ is a torus, and 
$S[P]$ is a toric variety with respect to this torus action. 

We have the following.

\begin{lemma}\label{toric-stack}(\cite{LogStack}, 5.14, 5.15)
The following two categories fibered in groupoids over the category of 
$S$-schemes are equivalent:

\begin{enumerate}
\item \quad the category of pairs $(X,\beta:P\to\overline{\sM}_X)$ 
consisting of a fine log scheme $X$ with a morphism $\underline{X}\to S$ 
and a morphism of sheaves of monoids $\beta$ that fppf locally lifts to a 
chart: $\widetilde{\beta}:P\to\sM_X;$

\item \quad the quotient stack $\cS_P:=[S[P]/S[P^{\emph{gp}}]].$
\end{enumerate}

If in addition $P$ is fs, then one can replace ``fppf" by ``\'etale".
\end{lemma}

\begin{anitem}\label{functoriality-toric-stack}
In fact, the action of $S[P^{\text{gp}}]$ on $S[P]$ extends to an action 
on the log structure on $S[P],$ and so this log structure descends to a 
log structure $\sM_{\cS_P}$ on the stack $\cS_P$ (cf. Section 
\ref{Satriano} for the definition of log structures on stacks), and there 
is a natural 
morphism $\pi_P:P\to\overline{\sM}_{\cS_P}$ of sheaves of monoids that 
fppf locally lifts to a chart. This is the universal pair 
$(\sM_{\cS_P},\pi_P)$ on $\cS_P$ that induces the equivalence in 
(\ref{toric-stack}) above. 

Moreover, for a morphism $h:Q\to P$ of fine monoids, the induced morphism 
$$
S[h]:S[P]\to S[Q]
$$
is compatible with the actions of $S[P^{\text{gp}}],S[Q^{\text{gp}}]$ and 
the homomorphism $S[h^{\text{gp}}]:S[P^{\text{gp}}]\to S[Q^{\text{gp}}],$ 
hence it descends to a morphism 
$$
S(h):\cS_P\to\cS_Q 
$$
of $S$-stacks. The map $h:Q\to P,$ regarded as a morphism of constant 
sheaves, induces a morphism $S(h)^*\sM_{\cS_Q}\to\sM_{\cS_P}$ of log 
structures, making $S(h)$ into a morphism of $S$-log stacks. 
\end{anitem}

\subsection*{The stack of log structures}
Now we can discuss the stack $LOG_S$ parameterizing fine log 
structures.

Let $S$ be a fine log scheme. Define $LOG_S$ to be the 
category with 

$\bullet$ objects: morphisms $X\to S$ of fine log 
schemes, and 

$\bullet$ morphisms: strict morphisms $X\to Y$ over 
$S.$

With the functor $(X\to S)\mapsto(\underline{X}\to\underline{ 
S})$ from $LOG_S\to\text{Sch}_{\underline{S}},$ this defines a 
fibered category over $\underline{S}.$ One of the main results 
in \cite{LogStack} is the following.

\begin{theorem}(\cite{LogStack}, 1.1)
$LOG_S$ is an algebraic stack locally of finite presentation 
over $\underline{S}.$
\end{theorem}

Here for an algebraic stack we use a slightly different 
definition from \cite[4.1]{LMB}. Namely the first axiom there 
that the diagonal is representable, separated and 
quasi-compact, is replaced by that the diagonal is 
representable and of finite presentation. In fact the stack 
$LOG_S$ is not quasi-separated \cite[3.17]{LogStack}. 

Here are two basic properties of $LOG_S.$

\begin{proposition}(\cite{LogStack}, 3.19, 3.20)
(1). The natural map $i_S:\underline{S}\to LOG_S$ corresponding to the 
identity morphism $S\to S$ is an open immersion;

(2). The 2-functor 
$$
S\mapsto LOG_S:\{\text{fine log 
schemes}\}\to\{\text{algebraic stacks}\}
$$
preserves fiber product. More precisely, if 
$$
\xymatrix{
X' \ar[r] \ar[d] & X \ar[d] \\
S' \ar[r] & S}
$$
is a Cartesian square of fine log schemes, then the induced 
diagram
$$
\xymatrix{
LOG_{X'} \ar[r] \ar[d] & LOG_X \ar[d] \\
LOG_{S'} \ar[r] & LOG_S}
$$
is a 2-Cartesian square of algebraic stacks.
\end{proposition}

\subsection*{What is $LOG_S$ good for?}
One can use this stack $LOG_S$ to reinterprete many concepts 
in log geometry. Note that for a morphism $f:X\to S$ of fine 
log schemes, the induced morphism $LOG(f):LOG_X\to LOG_S$ is 
faithful, hence representable.

\begin{definition}\label{def-LOG(f)}
Let $P$ be a property of representable morphisms of algebraic stacks. Then 
we say that $f:X\to S$ \emph{has property $LOG(P)$} if $LOG(f):LOG_X\to 
LOG_S$ has property $P.$ We say that $f$ \emph{has property weak $LOG(P)$} 
if the map $\underline{X}\to LOG_S$ corresponding 
to the given morphism $f:X\to S$ has property $P.$
\end{definition}

\begin{anitem}
\textbf{Caution:} The diagram 
$$
\xymatrix{
\underline{X} \ar[r]^-{i_X} \ar[d]_-{\underline{f}} & LOG_X 
\ar[d]^-{LOG(f)} \\
\underline{S} \ar[r]_-{i_S} & LOG_S}
$$
does not necessarily commute. It commutes if and only if $f$ 
is strict. In \cite[Section 2]{LogCot} a device is introduced in order to fix this issue, using stacks of diagrams of logarithmic structures.
\end{anitem}

Recall from (\ref{defn:logsmooth}) the notion of log smoothness and log 
\'etaleness.

\begin{theorem}\label{Th-LOG-sm}
For a morphism $f:X\to S$ of fine log schemes, $f$ is LOG 
smooth (resp. LOG \'etale) if and only if $f$ is log smooth 
(resp. log \'etale), if and only if $f$ is weakly LOG smooth 
(resp. weakly LOG \'etale).
\end{theorem} 

This is part of \cite[4.6]{LogStack}.

\begin{anitem}
Another application is the following. Consider the deformation problem for 
a log smooth integral morphism $f_0:X_0\to B_0$ of fine log schemes 
and a strict square-zero thickening $B_0\to B$ defined by an ideal 
$J\subset\sO_{\underline{B}}.$ In (\ref{thm:basiclogdef}) we gave the 
relation between this deformation problem and the cohomology 
groups of the log tangent bundle $T_{X_0/B_0}.$ The stack $LOG_S$ provides 
another way of thinking of this problem. 

By (\ref{Th-LOG-sm}), the log smooth morphism $f_0:X_0\to B_0$ induces a 
representable smooth morphism $\underline{X_0}\to LOG_{B_0},$ denoted 
$\cL_{f_0},$ and the deformation problem 
$$
\xymatrix{
X_0 \ar@{^{(}..>}[r] \ar[d]_-{f_0} & X \ar@{..>}[d]^-f \\
B_0 \ar@{^{(}->}[r] & B}
$$
is equivalent to the following 
$$
\xymatrix{
\underline{X_0} \ar@{..>}[r] \ar[d]_-{\cL_{f_0}} & \underline{X} 
\ar@{..>}[d]^-{\cL_f} \\ LOG_{B_0} \ar[r] & LOG_B.}
$$
The solution to this deformation problem is the cohomology groups of 
the \textit{ordinary tangent bundle} $T_{\underline{X_0}/LOG_{B_0}},$ 
therefore, theorem (\ref{thm:basiclogdef}) holds with $T_{X_0/B_0}$ 
replaced by $T_{\underline{X_0}/LOG_{B_0}}.$ In fact, we have 
$\Omega^1_{X_0/B_0}\cong\Omega_{\underline{ 
X_0}/LOG_{B_0}}$ (cf. (\cite{LogCot}, 3.8)). 

See section \ref{LogDef} for the general deformation theory of log 
schemes. 

\end{anitem}

\subsection*{Local structure of $LOG_S.$}\label{section-loc-stru}

For a fine log scheme $S,$ the relation between the quotient stacks 
$\cS_P$ and $LOG_S$ is that, the stack $LOG_S$ can be covered by the 
relative versions of the $\cS_P$'s. 

Let $u:U\to S$ be a strict morphism of fine log schemes, such that the 
underlying morphism $\underline{u}$ is \'etale. We will just say that $u$ 
is an \'etale strict morphism, if there is no confusion. Let 
$\beta:Q\to\sM_U$ be a chart, and let $h:Q\to P$ be a morphism of fine 
monoids. The chart $\beta$ induces a strict morphism 
$U\to\underline{S}[Q],$ which we also denote by $\beta.$ 

Let $\cS_P$ be the quotient stack 
$[\underline{S}[P]/\underline{S}[P^{\text{gp}}]]$ with the natural fine 
log structure $\sM_{\cS_P}$ in (\ref{functoriality-toric-stack}), and let 
$\underline{\cS_P}$ 
be the underlying stack. Consider the 2-commutative diagram 
$$
\xymatrix@C=.6cm @R=.5cm{
\underline{\cS_P}\times_{\underline{\cS_Q}}\underline{U} \ar[r]^-{pr_2} 
\ar[d]_-{pr_1} & \underline{U} \ar[r]^-{\underline{\beta}} & \underline
{\underline{S}[Q]} \ar[d]^-{\underline{\pi}} \\ 
\underline{\cS_P} \ar[rr]_-{\underline{\underline{S}(h)}} &&
\underline{\cS_Q}.}
$$
Let $Z$ be the $\underline{U}$-log stack with underlying stack 
$\underline{Z}=\underline{\cS_P}\times_{\underline{\cS_Q}}\underline{U}$ 
and the inverse image log structure $\sM_Z=pr_1^*\sM_{\cS_P}.$ Applying 
$pr_1^*$ to the morphism of log structures  
$\underline{S}(h)^*\sM_{\cS_Q}\to  
\sM_{\cS_P} 
$
(cf. (\ref{functoriality-toric-stack})), noting that 
$\pi\circ\beta:U\to\cS_Q$ is strict, we obtain a morphism 
$
pr_2^*\sM_U\to 
\sM_Z, 
$
making $pr_2$ into a morphism of log stacks $pr_2:Z\to U.$ This gives a 
morphism 
$$
\underline{Z}\to LOG_U\to LOG_S.
$$

\begin{proposition}\label{Prop-loc-stru}(\cite{LogStack}, 5.25)
For any fine log scheme $S,$ the natural morphism 
$$
\coprod_{(U,\beta,h)}\underline{\cS_P}\times_{\underline{\cS_Q}}
\underline{U}\to LOG_S
$$
is a representable \'etale surjection, where the disjoint 
union is taken over the isomorphism classes of all triples $(U,\beta,h)$ 
consisting of an \'etale strict morphism $U\to S,$ a chart 
$Q\overset{\beta}{\to}\sM_U,$ and a morphism $h:Q\to P$ of fine monoids, 
for some fine monoids $P$ and $Q.$
\end{proposition}

\section{Log deformation theory in general}\label{LogDef}


The main reference here is \cite{LogCot}.

As is well  known, the general deformation theory of schemes and
morphisms of schemes is not as easy
as in the smooth case. To understand deformation theory of general
morphisms, one has to use the full power of the cotangent
complex, see \cite{Illusie1,Illusie2}. In log geometry, one can
generalize it to get a reasonable theory of logarithmic cotangent
complex.

This log cotangent complex will be compatible with the usual cotangent
complex when the morphism in question is strict and is also compatible
with log smooth deformation theory for log smooth
morphisms. Basically, this is an application of the deformation theory
of representable morphisms to algebraic stacks (\cite{defStacks}) to
the classifying morphisms from the underlying scheme of \m X to the
stack \m{LOG_Y} (\cite{LogStack}, see also section \ref{LogStacks} of
this chapter.) 

\begin{conv}
We will focus on the category of fine log schemes. For a log scheme \m X, \underline{\m X} means the underlying scheme of \m X.
\end{conv}

\begin{rmk}We will work with the category \m{D'(\underline{X}_{\acute{e}t})} and similar categories, and one can talk about distinguished triangles and Ext's in these categories. For relevant definitions, see \cite{LogCot} and \cite{defStacks}.
\end{rmk}

Our presentation here follows \cite{LogCot}. An alternative approach
to the contangent complex due to Gabber is also explained
in \cite[Section 8]{LogCot}. 

\subsection*{The Log Cotangent Complex}
In Section \ref{LogStacks}, an Artin stack \m{LOG_Y} is defined for a
log scheme $Y$. It has the property that morphisms of log
schemes \m{X\ra Y} are equivalent to  morphisms \m{\underline{X}\ra
LOG_Y}. Thus one may 
interpret deformations of a morphism of log schemes \m{X\ra Y} as 
deformations of the associated representable
morphism   \m{\underline{X}\ra
LOG_Y}.  In \cite{defStacks}, the deformation theory of
representable morphisms of stacks was studied in detail. As an application of this
theory, one makes the following definition: 

\begin{defn}
For a morphism of log schemes \m{f:X\ra Y}, the {\em logarithemic
cotangent complex of \m{f}} is the complex \m{L_f=L_{\underline{X}/LOG_Y}}, where
the right hand side is the cotangent complex of the
morphism \m{\underline{X}\ra LOG_Y} defined in
Section \ref{LogStacks}. 
\end{defn}

\begin{rmk}
One should think about \m{L_f} as an object of the category \m{D'_{qcoh}(\underline{X}_{\acute{e}t})}. In the above definition, the right hand side is an object of the category \m{D'_{qcoh}(\underline{X}_{lis-\acute{e}t})}. As the restriction functor
\m{
D'_{qcoh}(\underline{X}_{lis-\acute{e}t})\ra D'_{qcoh}(\underline{X}_{\acute{e}t})
}
is an equivalence of categories, no information of the cotangent complex would lost.
\end{rmk}

\subsection*{Basic Properties}

For every morphism of fine log schemes \m{f:X\ra Y} the log cotangent complex is a projective system
\[
L_f=(\cdots\ra L_f^{\geq -n-1}\ra L_f^{\geq -n}\ra\cdots\ra L_f^{\geq0})
\]
where each \m{L_f^{\geq -n}} is an essentially constant ind-object in \m{D^{[-n,0]}(\ox)} (The derived category of \m{\ox}-modules supported in  \m{[-n,0]}).

The log cotangent complex \m{L_f} has the following properties:
\begin{enumerate}
\item For any \m{n\geq 0}, the natural map \m{\tau_{\geq -n}L_f^{\geq -n-1}\ra L_f^{\geq -n}} is an isomorphism.
\item If \m f is strict, then the system \m{(\tau_{\geq -n}L_{f'})} represents \m{L_f}, where \m{L_{f'}} is the usual cotangent complex of the underlying morphism of schemes \m{f'}.
\item If \m{f:X\ra Y} is log smooth, then the sheaf of log differentials \m{\Omega^1_{X/Y}} represents \m{L_f}.\label{logsm}
\item \label{split}If \[\squarediagram{X'}{X}{Y'}{Y}{a}{b}{g}{f}\] is a commutative diagram of fine log schemes, then there is a natural map \[
a^*L_f\ra L_g
\]
which is an isomorphism if the square above is cartesian and \m f is log flat. Furthermore, if the composite \m{X'\ra Y'\ra Y} satisfies the condition (T) below, then the map
\[
g^*L_b\oplus a^*L_f\ra L_{bg}
\]
is also an isomorphism.
\item \label{triangle}Given a composite
\[
\xymatrix{X\ar[r]^{f}&Y\ar[r]^{g}&Z}
\]
satisfying condition (T) below, there is a natural map
\[
L_f\ra f^*L_g[1]
\]
making the resulting triangle
\begin{equation}
f^*L_g\ra L_{gf}\ra L_f \ra f^*L_g[1]
\label{tri}
\end{equation}
distinguished.
\end{enumerate}

\begin{rmk}
In (\ref{split}) and (\ref{triangle}) above, \m{f^*, a^*, g^*} should
be understood in the derived sense. 
\end{rmk}

\begin{rmk}
One might hope for a theory of log cotangent complex in which every
triangle (\ref{tri}) is distinguished. This is unfortunately not the
case - an example due to W. Bauer is given in \cite[Section 7]{LogCot}. 
\end{rmk}

On the other hand, Gabber has shown (see \cite[Section 8]{LogCot})
that if one loosens the requirement \ref{logsm}, then one can obtain a
theory of log cotangent complexes for which one has a distinguished
triangle (\ref{tri}) for all
composites \xymatrix{X\ar[r]^{f}&Y\ar[r]^{g}&Z}.

The Condition (T) mentioned above is the following:

There exists a family of commutative diagrams
\[
\xymatrix
{
X_i \ar[r]^-{\pi_{X_i}}\ar[dr] &X\times_Y Y_i\ar[d]\ar[r]&X\times_Z Z_i\ar[d]\ar[r]&X\ar[d]_f\\
& Y_i \ar[r]^-{\pi_{Y_i}}\ar[dr] &Y\times_Z Z_i\ar[r]\ar[d]&Y\ar[d]_g\\
&&Z_i\ar[r]^-{\pi_{Z_i}}& Z
}
\]
such that
\begin{enumerate}
\item The underlying schemes of $X_i, Y_i, Z_i$ are all affine.
\item The \m{\pi}'s are all strict, and their underlying morphisms are flat and locally of finite presentation.
\item The underlying family of morphisms of schemes of \{\m{X_i\ra X}\} is jointly surjective.
\item There exists charts 
\[\beta_{X_i}:Q_{X_i}\ra \mM_{X_i},\beta_{Y_i}:Q_{Y_i}\ra \mM_{Y_i},\beta_{Z_i}:Q_{Z_i}\ra \mM_{Z_i}\]
and injective maps
\[
Q_{Z_i}\ra Q_{Y_i}\ra Q_{X_i}
\]
compatible with the morphisms \m{f_i, g_i} and
\[
\Tor^j_{\mO_{Z_i}\otimes_{\ZZ[Q_{Z_i}]}\ZZ[Q_{Y_i}]}(\mO_{Z_i}\otimes_{\ZZ[Q_{Z_i}]}\ZZ[Q_{X_i}], \mO_{Y_i}[G])=0 \text{ for all } j>0.
\]
Here \m{G:= \Coker(Q_{Z_i}^{gp}\ra Q_{Y_i}^{gp})} and \m{\mO_{Y_i}[G]} is viewed as an \m{\mO_{Z_i}\otimes_{\ZZ[Q_{Z_i}]}\ZZ[Q_{Y_i}]}-algebra via the map
\[
\mO_{Z_i}\otimes_{\ZZ[Q_{Z_i}]}\ZZ[Q_{Y_i}]\ra\mO_{Y_i}[G], t\otimes e_q\mapsto g_i^*(t)\beta_{Q_{Y_i}}(q)\cdot \bar{q}
\]
where \m{\bar{q}} denotes the image of \m q in \m G.
\end{enumerate}

\subsection*{Deformation Theory of Log Schemes in General}

In this section, we explain the relation between the log cotangent complex 
and deformation theory of log schemes. Let \m{f:X\ra Y} be a morphism of 
fine log schemes and let \m I be a quasi-coherent sheaf on \underline{\m 
X}. Define a \m{Y}-extension of \m{X} by \m{I} to be a commutative diagram 
of log schemes
\[
\xymatrix
{
X\ar[r]^{j}\ar[d]_f&X'\ar[dl]^{f'}\\
Y
}
\]
where \m{j} is an strict closed immersion 
 defined by a 
square-zero ideal, together with an isomorphism \m{\epsilon_j: I\cong 
\Ker(\mO_{X'}\ra\mO_X)}. The set of \m{Y}-extensions of $X$ by \m I forms, in a 
natural way, a
category \m{\underline{\Exal}_Y(X,I)}. Let \m{\Exal_Y(X,I)} be the set
of isomorphism classes of this category. 

There is a tautological equivalence of categories (see \cite[Problem
1]{defStacks} for the meaning of the right hand side): 
\[
\underline{\Exal}_Y(X,I)\cong\underline{\Exal}_{LOG_Y}(\underline{X},I).
\]
Hence, by \cite[Theorem 1.1]{defStacks} and our definition of \m{L_f},
we obtain the following result: 

\begin{thm} (\cite[Theorem 5.2]{LogCot})
There is a natural bijection
\[
\Exal_Y(X,I)\cong\Ext^1(L_f,I).
\]
\end{thm}

It is precisely the theorem above that guarantees that general deformation theory is controlled by our logarithmic cotangent complex.

\begin{defn}
Let \m{j_0:Y_0\into Y} be an strict closed immersion of fine log
schemes defined by a square-zero ideal \m{I\subset\mO_Y}, and
let \m{f_0:X_0\ra Y_0} be a LOG flat morphism
(Definition \ref{def-LOG(f)}). A {\itshape log flat deformation
of \m{X_0} 
to \m Y} is a cartesian square 
\[
\squarediagram{X_0}{X}{Y_0}{Y}{j}{j_0}{f_0}{f}
\]
with \m{f} LOG flat. 
\end{defn}

To give a log flat deformation as above is equivalent to give a 2-commutative diagram
\[
\squarediagram{\underline{X_0}}{\underline{X}}{LOG_{Y_0}}{LOG_Y}{j}{j_0}{\mL_{f_0}}{\mL_{f}} 
\]
with \m{\mL_f} flat. Thus from (\cite{defStacks}, 1.4) we obtain the following:

\begin{thm}
Let \m J denote the ideal of \m{LOG_{Y_0}} in \m{LOG_{Y}}. Then
\begin{enumerate}
\item{There exists a canonical class \m{o\in\Ext^2(L_{f_0},\mL_{f_0}^*J)} whose vanishing is equivalent to the existence of a log flat deformation of \m{X_0} to \m Y.}
\item{If \m{o=0}, then the set of isomorphism classes of log flat deformations of \m{X_0} to \m Y is naturally a torsor under \m{\Ext^1(L_{f_0},\mL^*_{f_0}J)}.}
\item{The automorphism group of any log flat deformation of \m{X_0} to \m{Y} is canonically isomorphic to \m{\Ext^0(L_{f_0},\mL^*_{f_0}J))}.}
\end{enumerate}
\end{thm}
This theorem gives an answer to the question of general deformation
theory of log schemes. 

\subsection*{Deformations of morphisms} As in the case of usual
schemes, once one understands deformations of log schemes, one obtains
a solution to the related problem of deformations of a log morphisms.

We are given a commutative diagram of solid arrows
\begin{equation}
\begin{xy}
(0,0)*+{X_0}="a"; (40,0)*+{X}="b";%
(20,-16)*+{Y_0}="c"; (60,-16)*+{Y}="d";%
(20,-30)*+{Z_0}="e"; (60,-30)*+{Z}="f";%
{\ar^-{i} "a";"b"};
{\ar^(0.3){j} "c";"d"};
{\ar^-{k} "e";"f"};
{\ar^-{f_0} "a";"c"};
{\ar_(0.7){h_0} "a";"e"};
{\ar@{-->}^{f} "b";"d"};
{\ar|>>>>>>>>>>>>>>{\hole}_(0.7){h} "b";"f"};
{\ar^{g_0} "c";"e"};
{\ar^{g} "d";"f"};
\end{xy}
\label{defmor}
\end{equation}
where $i,j,k$ are strict closed immersion defined by square-zero ideal sheaves $I,J,K$ living on $X,Y,Z$ respectively. The question is to find a dotted arrow \m f fitting in the diagram. To nail down \m f, we need some more data.

The morphisms \m h and \m g induce morphisms \m{w:h_0^*K\ra I} and \m{v: g_0^*K\ra J}. Assume given a morphism \m{u:f_0^*J\ra I} such that the composite\[
h_0^*K=f_0^*g_0^*K\ra f_0^*J\ra I
\]
is equal to \m w.

What we want to find is an \m f fitting in the diagram (\ref{defmor}) such that the morphism \m{f_0^*J\ra I} induced by \m f is equal to \m u. This can also be solved by the logarithmic cotangent complex.

\begin{thm}
In the situation above, assume in addition that \m u induces a map \m{\mL^*_{f_0}J'\ra I}, where \m{J'} is the ideal of \m{LOG_{Y_0}} in \m{LOG_{Y}}, then there is a canonical class \m{o\in\Ext^1(f^*_0L_{Y_0/Z_0},I)} which vanishes iff there exists \m f fitting in the diagram (\ref{defmor}) such that the morphism \m{f_0^*J\ra I} induced by \m f is equal to \m u. If \m{o=0}, then the set of such maps \m f is a torsor under the group \m{\Ext^0(f^*_0L_{Y_0/Z_0},I)}.
\end{thm}

For a proof of this theorem, see \cite{LogCot}, Theorem 5.9.


\section{Rounding}\label{Rounding}

The main reference for this section is \cite{Kato-Nakayama}. We have
benefitted from a lecture by A. Ogus, who reported on results in
\cite{Nakayama-Ogus}. 

\subsection*{What is rounding?}

The process of ``rounding", in its most basic form, produces a manifold \emph{with corners} from a \emph{smooth} analytic space with a normal crossings divisor. So the corners are not rounded but rather the opposite: they are created.  
 On the other hand, these corners are rather round and shapely. Also, anybody who has seen the construction under any name and hears the name ``rounding" immediately knows what this is about. Evidently then,  even though ``rounding" might be  something of a misnomer, it is a very good name.
The origin of the name seems to be in work of Kajiwara, Nakayama and Ogus \cite{Kajiwara-Nakayama,Nakayama-Ogus}.

In various moduli problems, the rounding of the moduli space often has
a more natural topological interpretation than the moduli space
itself.  A good  example is the Deligne-Mumford moduli space
$\overline{\mathcal{M}}_{g,n}$ of marked nodal curves, whose boundary
is a normal crossings divisor.  The ``interior" $\mathcal{M}_{g,n}$
can be described topologically as a quotient of Teichmuller space by
the appropriate mapping class group.  There is a natural
generalization of Teichmuller space  involving $2$-manifolds decorated
with circles, due to Harvey \cite{Harvey}; the analogous quotient
yields a topological description of the \emph{rounding} of
$\overline{\mathcal{M}}_{g,n}$, rather than the moduli space
$\overline{\mathcal{M}}_{g,n}$ itself. 

Another example is that of twisted curves, discussed in section
\ref{Roots}. A twisted curve is an algebraic stack with a log
structure, so it is a bit exotic. But its rounding is a good old
topological space. A similar  
example occurs in current work of one of us (Gillam): the relative
Hilbert stack of a marked Riemann surface can be defined algebraically
using Jun Li style expansions, but it is not representable (except in
some trivial cases). However its rounding is a topological space (even
a manifold with corners) which is relatively easy to describe.  A
similar phenomenon occurs in many moduli problems involving
expansions, discussed briefly below.   

The topological preeminence of the rounding of moduli spaces might
ultimately be traced back to the preference in topology for operations
involving real codimension one subspaces (e.g.\ connected sum of
manifolds) as opposed to the algebro-geometric preference for complex
codimension one operations (e.g.\ pushout of two smooth varieties
along a common divisor).  From this point of view, one might think of
log geometry as an attempt to speak algebraically about various 
``real codimension one" phenomena. 

\subsection*{The oriented real blowup}
The most general rounding operation is the \emph{Kato--Nakayama
  logarithmic space} associated to a log analytic space
\cite{Kato-Nakayama}.  In the basic example of a smooth analytic space
with log structure from a normal crossings divisor, the Kato--Nakayama
space can be described in terms of \emph{oriented real blowup}, which
is a relatively simple rounding operation that can be described as
follows. 

 Suppose $X$ is a topological space, $\pi : L \to X$ is a
complex line bundle, and $s : X \to L$ is a section of $\pi$.  Locally
on $X$ we can choose a trivialization $(\pi, \phi) : L \to X \times
\CC$ and consider the subspace  
\be \B_{L,s,\phi} X & := & \bigg\{\, l \in L\  :\ \  |\phi(l)|
\cdot(\phi s \pi)(l)\, = \ \phi(l) \cdot|(\phi s \pi)(l)|\, \bigg\}
\ee 
 of $L$.  A continuous function $u : X \to \CC^*$ yields a new
 trivialization $(\pi, u \cdot \phi)$, where $(u \cdot \phi)(l) := (u
 \pi)(l) \phi(l)$ .  The key observation is that $\B_{L,s,\phi} X=
 \B_{L,s,u \cdot \phi} X$, so the subspace $\B_{L,s, \phi} X$ is
 independent of the choice of $\phi$, hence one can define a subspace
 $\B_{L,s} X \subseteq L$ by defining it locally on $X$ using a
 trivialization, then gluing the locally defined subspaces.  From the
 local picture using a trivialization, it is clear that the subspace
 $\B_{L,s} X$ contains the zero section and $L|_{Z(s)}$ (where $Z(s)
 \subseteq X$ is the zero locus of $s$) and is invariant under the
 $\RR_{>0}$ action on $L$ inherited from the full $\CC^*$ scaling
 action.  We let $\B_{L,s}^* X$ be the complement of the zero section
 in $\B_{L,s}$ and we call \be \Blo_{L,s} X & := & (\B_{L,s}^* X ) /
 \RR_{>0} \ee the \emph{oriented real blowup} of $X$ along $(L,s)$.   

The space $\Blo_{L,s} X$ is a closed subspace of the oriented circle bundle $S^1 L := L^* / \RR_{>0}$ associated to $L$ and is, in particular, proper over $X$.  The projection $\tau : \Blo_{L,s} X \to X$ is an isomorphism away from $Z(s)$ and $\tau^{-1}(Z(s))$ is oriented circle bundle $S^1 L|_{Z(s)}$.  The spaces $\B_{L,s} X$ and $\Blo_{L,s} X$ are natural under pulling back line bundles and sections.

If $X$ is an analytic space and $D \subseteq X$ is a Cartier divisor, then $D$ determines a line bundle $\O_X(D)$ together with a section $s$ whose zero locus is $D$.  In this situation, we will write $\B_D X$, $\Blo_D X$, etc.\ and speak of the \emph{oriented real blowup of} $X$ \emph{along} $D$.  The space $\Blo_D X$ inherits a differentiable structure from its inclusion in $S^1 \O_X(D)$.

The basic example to keep in mind is the oriented real blowup $\Blo_0 \CC$ of the complex plane $\CC$ at the origin.  The origin is the zero locus of the identity map $\Id : \CC \to \CC$, hence \be \Blo_0 \CC & = & \big\{ (z,Z) \in \CC \times \CC^*\ :\ |z|Z = z |Z| \big\}\  /\  \RR_{>0} \\ & = & \big\{ (z,Z) \in \CC \times S^1\ :\ |z| Z = z \big\} \\ & \cong & \RR_{\geq 0} \times S^1 , \ee where the last isomorphism from $\RR_{\geq 0} \times S^1$ is given by $(\lambda, Z) \mapsto (\lambda Z , Z)$.  Evidently $\Blo_0 \CC$ is a half-infinite annulus whose boundary $S^1$ is the exceptional locus of $\tau : \Blo_0 \CC \to \CC$ (the fiber over the origin).

\subsection*{The Kato--Nakayama space} 
Let $(X,\mathcal{M}_X)$ be a fine and saturated logarithmic analytic space. 
\begin{definition}\label{def-Xlog}[\protect{\cite[1.2]{Kato-Nakayama}}]
We define its {\em canonical rounding}, or {\em Kato-Nakayama space}, denoted $X^{\log}$, as the space whose points are pairs $(x,F)$ where $x \in X$ 
and $F :  \mathcal{M}_{X,x} \to S^1$ is a monoid homomorphism satisfying  $F(u) = u(x) / |u(x)|$ for every $u \in \O_{X,x}^* \subseteq \mathcal{M}_{X,x}$. 
\end{definition}
This space has a natural topology. Let us describe the topology in the
special case where $\cM_X$ is the canonical log structure associated
to a Cartier divisor $D \subseteq X$, see Example \ref{NClog}.  
Locally on $X$ we can find $f_1,\dots,f_n \in \mathcal{M}_X(X)$ which, 
together with the units, generate $\mathcal{M}_X$.  The map \be X^{\log} & 
\to & X \times (S^1)^n \\ (x,F) & \mapsto & (x, F(f_{1,x}), \dots, F(f_{n,x})) \ee is then easily seen to be a monomorphism onto a closed subset of $X \times (S^1)^n$, so we give $X^{\log}$ the subspace topology so that this is a closed embedding.  Since one can check easily that this topology does not depend on the choice of generators $f_1,\dots,f_n$, the locally defined topologies glue to a topology on $X^{\log}$ making the projection $\tau : X^{\log} \to X$ given by $\tau(x,F) := x$ a proper map.

\begin{anitem} {\bf Relating Kato--Nakayama spaces to oriented blowups}
\label{Xlog-blowup}
\end{anitem}

There is a morphism \be \phi : X^{\log} & \to & \Blo_D X \\ (x,F) & 
\mapsto & (x, f \mapsto F(\overline{f})) \ee of topological spaces over 
$X$ which requires a little explanation.  Here $f \in S^1 \O_X(-D)|_x$ is in the circle bundle associated to the fiber $\O_X(-D)|_x \cong \CC$ and $\overline{f} \in \O_{X,x}(-D)$ is a lifting of $f$ to the stalk (one shows that any such $\overline{f}$ is actually in $\mathcal{M}_{X,x}$ and that $F(\overline{f})$ does not depend on this choice of lifting $\overline{f}$).  If we use the identification $$S^1 \O_X(D) = \Hom_{S^1}( S^1 \O_X(-D), S^1), $$ then we can think of $\phi$ as a map from $X^{\log}$ to $S^1 \O_X(D)$; one then shows that this $\phi$ is continuous and that $\phi$ factors through $\Blo_D X \subseteq S^1 \O_X(D)$.

When $X$ is a smooth analytic space and $D$ is a smooth divisor, the map $\phi$ is easily seen to be an isomorphism since one can reduce to the case $(X,D) = (\CC,0)$ on formal grounds.  Slightly more generally, if $X$ is smooth, but $D$ is only a normal crossings divisor, then locally we can write $D$ as a union of smooth divisors $D_1,\dots,D_i$ which look like the first $i$ coordinate hyperplanes in $\CC^n$ ($n = \dim X$), and we can define a variant of the oriented real blowup \be \Blo_D' X & := & (\Blo_{D_1} X ) \times_X \cdots \times_X (\Blo_{D_i} X)  \ee and a map $\phi : X^{\log}  \to  \Blo_D' X$.  In this local picture, the log structure $\mathcal{M}_X$ is the direct sum (in the category of log structures) of the log structures $\mathcal{M}_X^j$ from $D_1,\dots,D_i$, the associated Kato--Nakayama space $X^{\log}$ is the fibered product over $X$ of the $(X,\mathcal{M}_X^j)^{\log}$, and $\phi$ is just the fibered product over $X$ of the previously constructed isomorphisms $\phi_j : (X,\mathcal{M}_X^j)^{\log} \to \Blo_{D_j} X$.  The locally defined variants can be glued to define a global variant $\Blo_D' X$ of the oriented real blowup and an isomorphism $\phi : X^{\log} \cong \Blo_D' X$ of topological spaces over $X$.  

\begin{anitem} {\bf Topology, cohomology, and the Kato--Nakayama space}
\label{top-Xlog}
\end{anitem}
Locally, if $X = \CC^n$ and $D$ is the union of the first $i$ coordinate hyperplanes, then $D$ is the zero locus of $(z_1,\dots,z_n) \mapsto z_1 \cdots z_i \in \CC$ and we have \be \Blo_D X & = & \{ (z_1,\dots,z_n,Z) \in \CC^n \times S^1 : | z_1 \cdots z_i | Z = z_1 \cdots z_i \} \\ \Blo_D' X & = & \{ (\overline{z}, \overline{Z}) \in \CC^n \times (S^1)^i: |z_j| Z_j = z_j \; {\rm for} \; j=1,\dots, i \}. \ee In the general normal crossings divisor situation, the fiber of $\tau : X^{\log} \to X$ over a point $x \in X$ is naturally identified with $$S^1 N_{D_1/X}|_x \times \cdots \times S^1 N_{D_i / X}|_x,$$ where $D_1,\dots,D_i$ are the branches of $D$ containing $x$.  When $\dim X = n$, a point $y \in \tau^{-1}(x)$ has a neighborhood diffeomorphic to a neighborhood of the origin in $\RR_{\geq 0}^i \times \RR^{2n-i}$.  (Note $i \leq n$, so the depth of the corners in a Kato--Nakayama space is somewhat constrained.)  Recall that the \emph{topology} near the origin only depends on whether $i > 0$, but the differentiable structure depends on the actual value of $i$.  In particular, the topological boundary of the manifold $X^{\log}$ is given by $\tau^{-1}(D)$, and this manifold with boundary is homotopy equivalent to its interior, so $\H^*(X^{\log}) = \H^*(X \setminus D).$

The topology of morphisms of Kato--Nakayama spaces is also very nice, as shown by the following beautiful general result, see \cite[Theorem 0.3]{Nakayama-Ogus}:

\begin{theorem}
Let $f: X \to Y$ be a log smooth and integral morphism of fine log analytic spaces. Then the associated map $X^{\log} \to Y^{\log}$ is a topological submersion.
\end{theorem}

In fact, Nakayama and Ogus prove a more general result, replacing integrality by K. Kato's notion of {\em exact morphisms}, see \cite[Definition 4.6]{KKato}. 

The fact that the topology is nice suggests that one expects nice cohomological implications. This is indeed the original motivation leading Kato and Nakayama to define $X^{\log}$, see \cite[Theorem 0.2 (1)]{Kato-Nakayama}:

\begin{theorem}
Let $(X,M)$ be a fine and saturated log scheme with $X$ of finite type over $\CC$. 
Let $F$ be a constructible sheaf on the log-\'etale site 
$X_{\emph{log-\'et}}$, and let $F^{\log}$ be its pullback to the 
topological space $X^{\log}$. Then for all $q\in \ZZ$ we have 
$$ H^q(X_{\text{log-\'et}}, F) = H^q(X^{\log}, F^{\log}).$$
\end{theorem}

Very strong results hold true for de Rham cohomology. In fact the 
Kato--Nakayama space is a model for the log de Rham cohomology of $X$ in 
the sense that \be \H^*(X^{\log}, \CC) & = & \HH^*(X, \land^{\bullet} 
\Omega_{(X,M)}) \ee under mild assumptions on $X$. We discuss this in  
Section \ref{DeRham} below.

\subsection*{Kato--Nakayama spaces of expanded pairs}
Given a pair $(X,D)$ consisting of a smooth variety $X$ over $\CC$ with a smooth divisor $D \subseteq X$, the notion of an \emph{expanded pair} $t : \mathcal{X} \to B$ over a base $B$ arises in various relative curve counting theories.  The fiber of $t$ over a point $b \in B$ always looks like \be X[n]_0 & = & X \coprod_D \Delta_1 \coprod_D \cdots \coprod_D \Delta_n \ee (for an appropriate $n$), where $\Delta_i =  \PP( N_{D/X} \oplus \O_D)$ is a $\PP^1$ bundle over $X$.  Both $\mathcal{X}$ and $B$ have natural log structures making $t$ a log smooth map of log schemes.  The fiber of $t^{\log} : \mathcal{X}^{\log} \to B^{\log}$ over a point $c \in \tau_B^{-1}(b)$ looks like $$ X^{\log} \coprod_{c_1 : S^1 N_{D/X} \cong S^1 N_{D/ \Delta_1}} \Delta_1^{\log} \cdots \coprod_{ c_n : S^1 N_{D/\Delta_{n-1}} \cong S^1 N_{D / \Delta_n}} \Delta_n^{\log} , $$ where the choice of $c \in \tau^{-1}(b) \cong (S^1)^n$ determines the choice of orientation reversing $S^1$ bundle isomorphisms $c_1,\dots,c_n$.  Here each $\Delta_j$ has the log structure from the two copies of $D$, and $\Delta_j^{\log}$ is a cylinder bundle over $D$ (better: an $I$-bundle over $D^{\log}$).

The action of $(\CC^*)^n$ on $X[n]_0$ given by scaling the fibers of the $\PP^1$ bundles $\Delta_i$ is an action by isomorphisms of log schemes, so it lifts to an action on Kato--Nakayama spaces.  This lifted action is nontrivial on $B^{\log}$ as the $(S^1)^n$ factor of $(\CC^*)^n$ acts simply transitively on $\tau_B^{-1}(b) \cong (S^1)^n$.  In the usual moduli problems involving expansions, the isotropy group of a point $b$ involves elements of $(\CC^*)^n$ such that the induced action on $X[n]_0$ respects a map from a curve to $X[n]_0$, a subscheme of $X[n]_0,$ etc., and this isotropy is usually required to be finite to have a good moduli problem.  Since $G \cap \RR_{>0} = \{ \Id \}$ for any finite subgroup $G$ of $\CC^*$, the Kato--Nakayama space of the moduli problem is often representable even if the moduli problem itself is not.  This is always the case for moduli problems involving, say, quotients of sheaves on $X[n]_0$ pulled back from $X$, since these quotients themselves have no automorphisms and the only isotropy comes from the subgroup of $(\CC^*)^n$ preserving the quotient.

\section{Log de Rham and Hodge structures}\label{DeRham}
The main references of this section are 
\cite{Kato-Usui-book,Kato-Nakayama,K-M-N}. This section owes much to a 
lecture by Phillip Griffiths \cite{Gri}.

\subsection*{Moduli spaces of polarized Hodge structures.}
We assume the reader to be familiar with some basic concepts of Hodge 
theory. First of all, we briefly summarize the classical theory of 
the moduli spaces of polarized Hodge structures. 

\begin{anitem} {\bf The moduli space $M_h=\Gamma\backslash D_h.$}
\label{subsection1.1}
Let $n$ be an integer, and let $h$ be a sequence of non-negative 
integers $(h^{n,0},h^{n-1,0},\cdots,h^{0,n})$ 
satisfying $h^{p,q}=h^{q,p},$ called the \textit{Hodge 
numbers}. Let $H_{\mathbb Z}$ be a free abelian group of 
rank $\sum h^{p,n-p},$ with a non-degenerate bilinear 
form $Q:H_{\mathbb Z}\otimes H_{\mathbb Z}\to\mathbb Z,$ 
which is symmetric (resp. anti-symmetric) if $n$ is even 
(resp. odd). Let $G_{\mathbb Z}$ be the group functor 
$\text{Aut}(H_{\mathbb Z},Q)$ on commutative rings, sending a ring $R$ 
to the group of automorphisms on the free $R$-module 
$H_R:=H_{\mathbb Z}\otimes R$ preserving the bilinear form 
$Q.$ It is an affine group scheme of finite type over $\mathbb Z$ (which 
is clear if we write down the matrix representing the bilinear form $Q$ 
with respect to some basis of $H_{\mathbb Z}$). Let 
$\Gamma$ be an arithmetic subgroup of $G_{\mathbb 
Z}(\mathbb Z).$ 
\end{anitem}

The set of Hodge structures of weight $n$ on $H_{\mathbb 
R}$ with prescribed Hodge numbers $h,$ such that $Q$ induces a 
\textit{polarization} on $H_{\mathbb R}$ (i.e. it induces a 
morphism $H_{\mathbb R}\otimes H_{\mathbb R}\to\mathbb R(-n)$ of 
Hodge structures, and the bilinear form $Q_C(u,v):=Q(u,Cv),$ 
where $C$ is the Weil operator, is symmetric and positive 
definite), is parameterized by the homogeneous space 
$D_h=G_{\mathbb R}/K,$ where $K$ is the stabilizer group of a 
fixed polarized Hodge structure $F_0$ on $H_{\mathbb R}.$ See for instance 
(\cite{Del1}) for these concepts.  

This homogeneous space $D=D_h=G_{\mathbb R}/K$ has a 
complex structure defined as follows. It is clear that 
$Q:H_{\mathbb R}\otimes H_{\mathbb R}\to\mathbb R(-n)$ is a 
morphism of Hodge structures if and only if 
$Q(F^p,F^{n-p+1})=0$ for all $p.$ Let $f^p=\sum_{r\ge 
p}h^{r,n-r},$ and let $D^{\vee},$ the \textit{compact dual 
of $D,$} be the subspace of the product of the 
Grassmannians $\prod_p\text{Gr}(f^p,H_{\mathbb R})$ 
consisting of all flags $F^{\bullet}:$ 
$$
\cdots\subset F^{p+1}\subset F^p\subset\cdots 
$$
such that $Q(F^p,F^{n-p+1})=0.$ Then $D^{\vee}\simeq G_{\mathbb 
C}/P,$ where $P$ is a parabolic subgroup preserving a 
fixed flag. This gives $D^{\vee}$ a complex structure. We 
see that $D\subset D^{\vee}$ is the locus of flags 
satisfying 

(i) $F^p\cap\overline{F}^{n-p+1}=0$ (so that $F^p\oplus 
\overline{F}^{n-p+1}\cong H_{\mathbb C})$ for all $p,$ and 

(ii) $Q(\overline{u},Cu)>0$ for all $u\ne0$ in $H_{\mathbb C}.$

They are both open conditions, so $D\subset D^{\vee}$ is an 
open complex submanifold. The group $\Gamma$ acts on $D_h$ 
properly discontinuously, and the quotient $M_h=\Gamma 
\backslash D_h$ is the moduli space of 
\textit{$\Gamma$-equivalence classes of $Q$-polarized Hodge 
structures on $H_{\mathbb C}$ with Hodge type $h.$} See 
(\cite{Kato-Usui-book}, 0.3.6, 0.3.7). 

\subsubsection*{Variations of Hodge structures.}

\begin{definition}\label{defn-VHS}(\cite{Del1}, 3.11, 3.12)
Let $S$ be a complex manifold. A \emph{variation of Hodge 
structures $\mathscr H$ of weight $n$ on $S$} is given by 

$\bullet$ a local system $\mathscr H_{\mathbb Z}$ of free 
abelian groups of finite rank on $S;$ 

$\bullet$ a finite decreasing filtration 
$F^{\bullet}\mathscr H_{\mathscr O}$ of the vector bundle 
$\mathscr H_{\mathscr O}:=\mathscr H_{\mathbb Z}\otimes_{ 
\mathbb Z}\mathscr O_S$ by holomorphic sub-bundles, 

\noindent such that the following conditions are satisfied: 

1) (Griffiths transversality) the natural flat connection 
$\nabla=d\otimes\text{id}_{\mathscr H_{\mathscr 
O}}:\mathscr H_{\mathscr O}\to\Omega^1_S\otimes\mathscr 
H_{\mathscr O}$ takes $F^p\mathscr H_{\mathscr O}$ into 
$\Omega^1_S\otimes F^{p-1}\mathscr H_{\mathscr O},$ for 
every $p;$ 

2) for each point $s\in S,$ the fiber $F^{\bullet}(s)$ over 
$s$ is a Hodge structure of weight $n.$ 

A \emph{polarization} of the variation of Hodge structures 
$\mathscr H$ is a locally constant bilinear form 
$$
\mathscr Q:\mathscr H_{\mathbb Z}\otimes\mathscr H_{\mathbb 
Z}\to\underline{\mathbb Z}
$$ 
such that on each fiber over $s\in S,$ it induces a 
polarization of the fiber Hodge structure. 
\end{definition}

Suppose we have a polarized family of Hodge structures 
$(\mathscr H,\mathscr Q:\mathscr H_{\mathbb Z}\otimes\mathscr 
H_{\mathbb Z}\to\mathbb Z)$ of weight $n$ on $S$ (not necessarily a 
variation of Hodge structures), and a global 
section of the sheaf 
$$
\Gamma\backslash\underline{\text{Isom}}((\mathscr 
H_{\mathbb Z},\mathscr Q),(H_{\mathbb Z},Q)),
$$
where $H_{\mathbb Z}$ is regarded as a constant sheaf on $S.$ 
If the monodromy group of this family of Hodge structures on $S$ is 
contained in $\Gamma,$ i.e. the image of the homomorphism 
$$
\pi_1(S)\to G_{\mathbb Z}(\mathbb Z)
$$
is contained in $\Gamma,$ then there is a well-defined function 
$$
\varphi:S\to M_h 
$$
inducing this family of Hodge structures. This map is 
locally liftable to $D_h.$ 

If $f:X\to S$ is 
a projective smooth morphism between quasi-projective complex 
algebraic manifolds, with a relative hyperplane section 
$\eta\in H^0(S,R^2f_*\mathbb Z),$ then the family of the 
primitive parts $P^n(X_s,\mathbb Z)$ of the cohomology groups 
$H^n(X_s,\mathbb Z)$ modulo torsion form a polarized variation 
of Hodge structures of weight $n$ on $S,$ and it induces a map 
$$
\varphi:S\to M_h
$$ 
if $\Gamma$ contains the monodromy group. The map $\varphi$ is called the 
\textit{period map.} To be precise, the family of $H^n(X_s,\mathbb C)$'s 
are the stalks of $R^nf_*(f^{-1}\mathscr O_S),$ and the Hodge filtration 
on $R^nf_*(f^{-1}\mathscr O_S)$ is given by the degenerate 
spectral sequence 
$$
E_1^{pq}=R^qf_*\Omega_{X/S}^p\Longrightarrow 
R^{p+q}f_*(f^{-1}\mathscr O_S),
$$
which is induced from the natural quasi-isomorphism (the relative 
holomorphic Poincar\'e lemma) 
$$
f^{-1}\mathscr O_S\longrightarrow\Omega_{X/S}^{\bullet}.
$$ 
Since $\eta$ is a global section, 
the primitive part form a family of sub-Hodge structures on $S.$ Griffiths 
proved that the period map is holomorphic, and the family of Hodge 
structures is a variation of Hodge structures. See 
(\cite{Del1}, 3) for more detail. 

For instance, one can take $S$ to be the moduli space $A_g$ (resp. $M_g$) 
of principally polarized abelian varieties of dimension $g$ (resp. 
projective smooth curves of genus $g$), and take $X$ to be the universal 
family of such objects over $S,$ and $n=1.$ In this case, the local system 
$\mathscr H_{\mathbb Z}=R^1f_*\underline{\mathbb Z}$ is torsion-free and 
is equal to its primitive part, and the filtration $F^{\bullet}\mathscr 
H_{\mathscr O}$ is given by $F^0\mathscr H_{\mathscr O}=\mathscr 
H_{\mathscr O},\ F^1\mathscr H_{\mathscr O}=R^0f_*\Omega_{X/S}$ and 
$F^2\mathscr H_{\mathscr O}=0.$ This variation of Hodge structures (and 
its tensor powers) are of great arithmetic interest. 

\subsection*{Logarithmic Hodge structures.}

Consider the following situation. Let $f:X\to S$ be a family 
of projective manifolds, where $S$ is the complement of a 
normal crossing divisor $D$ in a compact manifold 
$\overline{S},$ and suppose one can extend the family $f$ to a 
family $\overline{f}:\overline{X}\to\overline{S}$ which is log 
smooth (here $\overline{S}$ has the log structure induced by 
the divisor $D$ (\ref{NClog})). Then one can ask if it is possible to 
enlarge the moduli space $M_h$ to some $\overline{M}_h$ so that the period 
map extends to $\overline{\varphi}:\overline{S}\to\overline{M}_h.$ 

To study the degenerations of Hodge structures, Kato and 
Usui introduced the notion of logarithmic Hodge structures 
\cite{Kato-Usui-book}. 

\subsubsection*{Log de Rham complex.}

Let $(X,\alpha:M_X\to\mathscr O_X)$ be an fs log analytic 
space over $\mathbb C,$ and let $X^{\log}$ be the Kato-Nakayama space 
(Definition \ref{def-Xlog}), with $\tau=\tau_X:X^{\log}\to X$ the natural 
proper map. Over the open set $X^*\subset X$ where the log structure is 
trivial, the map $\tau$ is a homeomorphism, and the section 
$j^{\log}:X^*\hookrightarrow X^{\log}$ is a homotopy 
equivalence. For $x\in X,$ the fiber $\tau^{-1}(x)$ is a compact torus 
$(S^1)^m,$ where $m$ is the rank of $\overline{M}_{X,x}^{\text{gp}}.$ For 
instance, let $\Delta$ be the open unit disk $\{|z|<1\}$ equipped with the 
log structure induced by the center $\{z=0\}$ (cf. \ref{NClog}). Then 
$\Delta^{\log}$ is homeomorphic to $[0,1)\times S^1.$ See 
sections \ref{Xlog-blowup} and \ref{top-Xlog} for more examples. 

One can define a sheaf of rings $\mathscr O_{X^{\log}}$ on 
$X^{\log}.$ Roughly speaking, this is the subsheaf of rings of 
$j^{\log}_*\mathscr O_{X^*}$ on $X^{\log}$ generated over 
$\tau^{-1}\mathscr O_X$ by $``\log(q)",$ for all $q\in 
M^{\text{gp}}_X.$ See (\cite{Kato-Nakayama}, Section 3) for the precise 
definition. 

For example, if $x\in X$ and $y\in\tau^{-1}(x),$ and the free 
abelian group $\overline{M}_{X,x}^{\text{gp}}$ has rank $m$ 
and is generated by $f_1,\cdots,f_m\in M_{X,x}^{\text{gp}},$ 
then the stalk $\mathscr O_{X^{\log},y}$ is isomorphic to the 
polynomial ring $\mathscr O_{X,x}[\log(f_1),\cdots,\log(f_m)].$ This shows 
that in general, $(X^{\log},\mathscr O_{X^{\log}})$ is not a locally 
ringed space. 

For a morphism $f:X\to Y$ of fs log analytic spaces, one can define the 
\textit{sheaf of relative log differentials} $\Omega_{X/Y}^1$ in the same 
way as Definition \ref{LogDiff}, namely it is the sheaf representing the 
functor of $Y$-log derivations (Definition \ref{LogDer}), where we use the 
sheaf $\mathscr O_X$ of holomorphic functions as the structure sheaf. The 
explicit description in Proposition \ref{prop-logdiff} still applies. If 
$(\partial,D)$ is the universal $Y$-log derivation of $X$ to 
$\Omega_{X/Y}^1,$ the morphism $D:M_X\to\Omega_{X/Y}^1$ is also written as 
$d\log,$ and it can be extended by linearity to $M_X^{\text{gp}}.$ In the 
explicit description 
$$
\Omega_{X/Y}^1=(\Omega_{\underline{X}/\underline{Y}}\oplus\mathscr 
O_{\underline{X}}\otimes_{\mathbb Z}M_X^{\text{gp}})/\cK,
$$
$d\log(a)$ is the image of $0\oplus(1\otimes a),$ for a local section $a$ 
of $M_X^{\text{gp}}.$ 

The sheaf $\Omega_{X/Y}^1$ is an analytic coherent $\mathscr O_X$-module. 
For an integer $r\ge1,$ let $\Omega_{X/Y}^r$ be the $r$-th exterior power 
of $\Omega_{X/Y}^1.$ The derivation $\partial:\mathscr 
O_{\underline{X}}\to\Omega_{X/Y}^1$ can be prolonged to a complex 
$\Omega^{\bullet}_{X/Y}:$
$$
\xymatrix@C=.6cm{
\mathscr O_{\underline{X}} \ar[r]^-{\partial} & \Omega_{X/Y}^1 \ar[r]^-d & 
\Omega_{X/Y}^2 \ar[r]^-d & \cdots \ar[r] & \Omega_{X/Y}^r \ar[r] & \cdots}
$$
by imposing that $d(d\log(a))=0$ for $a\in M_X^{\text{gp}}.$ This is a 
complex of $\underline{f}^{-1}\mathscr O_{\underline{Y}}$-modules, called 
the \textit{relative log de Rham complex on $X$ with respect to $f.$}

For any sheaf $F$ of $\mathscr O_X$-modules, define 
$$
\tau^*F:=\tau^{-1}F\otimes_{\tau^{-1}\mathscr O_X}\mathscr O_{X^{\log}}
$$
as a sheaf on $X^{\log}.$ For an integer $r\ge1,$ define 
$$
\Omega_{X^{\log}/Y^{\log}}^r:=\tau^*\Omega_{X/Y}^r.
$$
The structure sheaf $\mathscr O_{X^{\log}}$ comes with a natural 
derivation $d:\mathscr O_{X^{\log}}\to\Omega_{X^{\log}/Y^{\log}}^1$ 
(see (\cite{Kato-Nakayama}, 3.5)), 
which can be prolonged to a complex $\Omega^{\bullet}_{X^{\log}/Y^{\log}}$ 
of $(f^{\log})^{-1}\mathscr O_{Y^{\log}}$-modules 
$$
\xymatrix@C=.6cm{
\mathscr O_{X^{\log}} \ar[r]^-d & \Omega_{X^{\log}/Y^{\log}}^1 \ar[r]^-d &
\Omega_{X^{\log}/Y^{\log}}^2 \ar[r]^-d & \cdots \ar[r] & 
\Omega_{X^{\log}/Y^{\log}}^r \ar[r] & \cdots,}
$$
which can be called the \textit{relative log de Rham complex on 
$X^{\log}$ with respect to $f.$} 

When $Y$ is a point with trivial log structure, we denote 
$\Omega^{\bullet}_{X/Y}$ (resp. $\Omega^{\bullet}_{X^{\log}/Y^{\log}}$) by 
$\Omega^{\bullet}_X$ (resp. $\Omega^{\bullet}_{X^{\log}}$), and call it 
the \textit{absolute log de Rham complex on $X$} (resp. $X^{\log}$).

Under mild conditions, F. Kato proved the relative log Poincar\'e lemma 
and the logarithmic analogue of the de Rham theorem. We state in the 
following a weaker version. See \cite{FKato-Poincare} for the more general 
version. We say an fs log analytic space is \textit{log smooth} if it is 
so over a point with trivial log structure $(pt,\mathbb C^*).$ 

\begin{theorem}\label{log-Poincare-lemma}(\cite{FKato-Poincare}, 3.4.2, 
3.2.5)
Let $f:X\to Y$ be a log smooth morphism of fs log analytic spaces, with 
$Y$ log smooth. Assume that the 
induced morphism $f^{-1}\overline{M}_Y\to\overline{M}_X$ is injective, 
that the stalk $\overline{M}_{X/Y,x}$ of the relative characteristic 
(\ref{def:char}) is torsion-free for every $x\in X,$ and that $f$ is exact 
(\cite{KKato}, 4.6). Then there is a natural quasi-isomorphism 
$$
(f^{\log})^{-1}\mathscr 
O_{Y^{\log}}\longrightarrow\Omega^{\bullet}_{X^{\log}/Y^{\log}}.
$$
\end{theorem}

\begin{cor}\label{log-de-Rham-thm}(\cite{FKato-Poincare}, 4.1.5)
Let $f:X\to Y$ be as in Theorem \ref{log-Poincare-lemma} above, and 
assume in addition that $\underline{f}$ is proper. Then there is a natural 
quasi-isomorphism 
$$
\tau_Y^*Rf_*\Omega^{\bullet}_{X/Y}\longrightarrow Rf^{\log}_*\underline 
{\mathbb Z}\otimes_{\mathbb Z}\mathscr O_{Y^{\log}}.
$$
\end{cor}

These results apply in particular to semi-stable degenerations. Also, the 
absolute log Poincar\'e lemma was proved earlier by Kato and Nakayama 
(\cite{Kato-Nakayama}, Theorem 3.8).

\subsubsection*{Log variations of polarized Hodge structures.}

For $y\in X^{\log}$ and $x=\tau(y)\in X,$ let $\text{sp}(y)$ be the set of 
all ring homomorphisms $s:\mathscr O_{X^{\log},y}\to\mathbb C$ that extend 
the evaluation map $\text{ev}_x:\mathscr O_{X,x}\to\mathbb C.$ Since 
$\mathscr O_{X^{\log},y}$ is isomorphic to the polynomial ring over 
$\mathscr O_{X,x}$ generated by log of a basis for $\overline{M}_{X,x},$ 
if we fix an $s_0\in\text{sp}(y),$ then we have a bijection: 
$$
s\mapsto(f\mapsto s(\log(f))-s_0(\log(f))):\text{sp}(y) 
\overset{\sim}{\longrightarrow}Hom_{\text{group}}(\overline 
{M}_{X,x}^{\text{gp}},\mathbb C), 
$$
where $\mathbb C$ is viewed as an additive group. 

\begin{definition}\label{LVPHS}
Let $X$ be an fs log analytic space. A \emph{log variation of 
polarized Hodge structures of weight $n$} (abbreviated as \emph{LVPHS}) on 
$X$ is given by 

$\bullet$ a local system of free abelian groups of finite rank 
$\mathscr H_{\mathbb Z}$ on $X^{\log},$ 

$\bullet$ a bilinear form $\mathscr Q:\mathscr H_{\mathbb 
Z}\otimes \mathscr H_{\mathbb Z}\to\underline{\mathbb Z},$ 

$\bullet$ a finite decreasing filtration $F^{\bullet}\mathscr 
H_{\mathscr O}$ of $\mathscr H_{\mathscr O}:=\mathscr 
H_{\mathbb Z}\otimes_{\mathbb Z}\mathscr O_{X^{\log}}$ by $\mathscr 
O_{X^{\log}}$-submodules,

\noindent such that the following conditions are satisfied: 

1) there exist a locally free $\mathscr O_X$-module $\mathscr 
E$ and a finite decreasing filtration $F^{\bullet}\mathscr E$ 
by $\mathscr O_X$-submodules, such that $\text{Gr}_p(\mathscr 
E)$ is locally free and 
$$
F^p\mathscr H_{\mathscr O}\cong\tau^*F^p\mathscr E
$$
for each $p;$

2) for $y\in X^{\log}$ and $x=\tau(y)\in X,$ let $s\in\text 
{sp}(y)$ and let $f_1,\cdots,f_r\in M_{X,x}-\mathscr 
O_{X,x}^*$ generate the monoid $\overline{M}_{X,x}.$ If the real numbers 
$|\exp(s(\log(f_i)))|$ are sufficient small for all $i,$ 
then $(\mathscr H_{\mathbb Z,y},\mathscr Q,F^{\bullet}(s))$ is 
a polarized Hodge structure of weight $n;$
 
3) (Griffiths transversality) the connection 
$\nabla=d\otimes\text{id}_{\mathscr H_{\mathbb 
Z}}:\mathscr H_{\mathscr O}\to\Omega_{X^{\log}}^1\otimes_{\mathscr 
O_{X^{\log}}}\mathscr H_{\mathscr O}$ takes $F^p\mathscr 
H_{\mathscr O}$ into $\Omega_{X^{\log}}^1\otimes 
F^{p-1}\mathscr H_{\mathscr O}.$
\end{definition}

Here $F^{\bullet}(s),$ the \textit{specialization of $F^{\bullet}$ at 
$s$}, is the decreasing filtration of $\mathscr H_{\mathbb 
C,y}:=\mathbb C\otimes_{\mathbb Z}\mathscr H_{\mathbb Z,y}$ 
defined by $F^p(s)=\mathbb C\otimes_{s,\mathscr 
O_{X^{\log},y}}(F^p\mathscr H_{\mathscr O})_y.$ For a fixed point 
$y\in X^{\log},$ the family $(\mathscr H_{\mathbb 
Z,y},\mathscr Q,F^{\bullet}(s))_{s\in\text{sp}(y)}$ is called 
a \textit{polarized log Hodge structure of weight $n$} on the log point 
$(x,M_{X,x});$ this is the same as a log variation of 
polarized Hodge structures of weight $n$ on the log point $(x,M_{X,x}).$    

\begin{anitem}\label{two-defn-LVPHS}\textbf{Variant.} 
The definition we gave here follows (\cite{Kato-Usui-book}, 0.2.19), 
except that our polarization $\mathscr Q$ is integral. This definition 
differs slightly from the one as in (\cite{Matsubara}, 5.3) and 
(\cite{K-M-N}, 2.3), and is weaker. The main difference is that, in 
\textit{loc. cit.} the locally free $\mathscr O_X$-module $\mathscr E$ 
with its filtration $F^{\bullet}\mathscr E$ is part of the data of the 
definition, and the flat connection is for $\mathscr E$ on $X$ (namely 
$\nabla:\mathscr E\to\Omega_X^1\otimes\mathscr E$) and is required to 
satisfy the Griffiths transversality. 
\end{anitem}

\begin{anitem}\label{LVPHS-geometry}\textbf{LVPHS from geometry.} 
Log variations of polarized Hodge structures arise from 
geometry in the following way. Let $f:X\to Y$ be a projective vertical 
log smooth morphism between log smooth fs log analytic spaces, and we fix 
a line bundle on $\underline{X}$ which is relatively very ample over 
$\underline{Y}.$ Here ``vertical" means $f^{-1}(Y^*)=X^*.$ By a 
theorem of Kajiwara and Nakayama (\cite{Kajiwara-Nakayama}, 0.3), for 
every integer $n,$ the sheaf $R^nf^{\log}_*\underline{\mathbb Z}$ is a 
local system on $Y^{\log}.$ We take $\mathscr H_{\mathbb Z}$ to be 
$R^nf^{\log}_*\underline{\mathbb Z}$ modulo torsion, take $\mathscr Q$ to 
be the pairing induced by the fixed line bundle (which is obtained in 
the same way as (\cite{K-M-N}, 8.2) where we replace all local systems 
of $\mathbb Q$-vector spaces by the integral lattices in them), take 
$\mathscr E$ to be $R^nf_*(\Omega_{X/Y}^{\bullet}),$ with filtration 
$F^p\mathscr E=R^nf_*(\Omega_{X/Y}^{\ge p})\subset\mathscr E,$ and take 
$F^p\mathscr H_{\mathscr O}$ to be $\tau^*F^p\mathscr E.$ Then by a 
theorem of Kato, Matsubara and Nakayama (\cite{K-M-N}, Theorem 8.1), this 
is a log variation of polarized Hodge structures on $Y,$ even in the 
stronger sense (\cite{K-M-N}, 2.3). See (\cite{K-M-N}, Theorem 8.1) for a 
more general version with coefficients. The special case when $Y$ is the 
unit disk $\{|z|<1\}$ in the complex plane with log structure induced by 
the divisor $\{z=0\}$ and $f$ is family of projective manifolds with 
semi-stable reduction over $\{z=0\}$ was proved earlier by Matsubara 
(\cite{Matsubara}, Theorem C), except for the polarization part. 
\end{anitem}

\begin{anitem}\label{can-extn}\textbf{Relation to Deligne's canonical 
extensions.} 
Suppose that $\underline{X}$ and $\underline{Y}$ are smooth, that the 
log structure on $Y$ is induced by a normal crossing divisor $D\subset Y,$ 
and that $f:X\to Y$ is a morphism of semi-stable reduction. Let $f':X^*\to 
Y^*$ be the restriction of $f.$ Then the flat connection 
$(R^nf'_*\Omega_{X^*/Y^*}^{\bullet},\nabla')$ on $Y^*$ has unipotent local 
monodromy around each component of $D,$ and the flat connection 
$(R^nf_*\Omega_{X/Y}^{\bullet},\nabla)$ on $Y$ is its canonical extension 
in the sense of Deligne \cite{LNM163}.
\end{anitem}

\subsection*{Kato-Usui spaces.}

We fix $n,h,H_{\mathbb Z},Q,G_{\mathbb Z},D$ and $D^{\vee}$ 
as in (\ref{subsection1.1}). For a ring $R,$ let $\mathfrak 
g_R=\text{Lie}(G_R).$ A subset $\sigma\subset\mathfrak 
g_{\mathbb R}$ is called a \textit{nilpotent cone} if it is a cone 
$$
\sigma=\sum_{i=1}^n\mathbb R_{\ge0}N_i
$$
generated by mutually commutative nilpotent operators $N_i\in\mathfrak 
g_{\mathbb R}\subset\text{End}(H_{\mathbb R}).$ It is called a 
\textit{rational nilpotent cone} if it can be generated by nilpotent 
operators in $\mathfrak g_{\mathbb Q}.$ 
Let $\Gamma$ be a \textit{neat 
subgroup} of $G_{\mathbb Z}(\mathbb Z),$ i.e. for every element 
$\gamma\in\Gamma,$ its eigenvalues on $H_{\mathbb C}$ generate 
a torsion-free subgroup of $\mathbb C^*.$

\subsubsection*{Nilpotent orbits.}

\begin{definition}\label{nilp-orbit}
Let $\sigma=\sum_i\mathbb R_{\ge0}N_i$ be a nilpotent cone. A 
subset $Z\subset D^{\vee}$ is called a 
\emph{$\sigma$-nilpotent orbit}, if there exists an 
$F_0\in D^{\vee}$ such that 

$\bullet\ Z=\exp(\sum_i\mathbb CN_i)F_0,$

$\bullet\ NF^p_0\subset F^{p-1}_0$ for all $p\in\mathbb Z$ and 
$N\in\sigma,$

$\bullet\ \exp(\sum_iz_iN_i)F_0\in D$ if $\text{Im}(z_i)\gg0$ 
for all $i.$
\end{definition}

We also call the pair $(\sigma,Z)$ a \textit{nilpotent orbit.} 

Let $\Sigma$ be a \textit{fan} in $\mathfrak g_{\mathbb Q},$ 
i.e. $\Sigma$ is a non-empty set of rational nilpotent cones 
in $\mathfrak g_{\mathbb R}$ such that 

$\bullet$ if $\sigma\in\Sigma,$ then all faces of $\sigma$ are 
in $\Sigma,$

$\bullet$ for $\sigma,\sigma'\in\Sigma,$ the intersection 
$\sigma\cap\sigma'$ is a face of both $\sigma$ and $\sigma',$ 

$\bullet$ for every $\sigma\in\Sigma,$ we have 
$\sigma\cap(-\sigma)=0.$

One can then define the set $D_{h,\Sigma}$ (or just 
$D_{\Sigma},$ if there is no confusion) of \textit{nilpotent 
orbits in the directions in $\Sigma$} to be the set of 
nilpotent orbits $(\sigma,Z)$ where $\sigma\in\Sigma.$ There 
is a natural injection 
$$
F\mapsto(0,\{F\}):D\hookrightarrow D_{\Sigma}.
$$

\subsubsection*{The moduli space $M_{\Sigma}$.}

Let $\Sigma$ be a fan in $\mathfrak g_{\mathbb Q}$ and let 
$\Gamma\subset G_{\mathbb Z}(\mathbb Z)$ be a subgroup. Then 
we say that \textit{$\Gamma$ is compatible with $\Sigma$} if 
for every $\gamma\in\Gamma$ and $\sigma\in\Sigma,$ we have 
$Ad(\gamma)(\sigma)\in\Sigma.$ In this case, there is an 
action of $\Gamma$ on $D_{\Sigma}$ given by 
$$
\xymatrix@C=1cm{
(\sigma,Z) \ar@{|->}[r]^-{\gamma} & 
(Ad(\gamma)(\sigma),\gamma Z).}
$$
We say that \textit{$\Gamma$ is strongly compatible with 
$\Sigma$} if every cone $\sigma\in\Sigma$ is generated by 
elements in $\log\Gamma.$ Kato and Usui showed that when 
$\Gamma$ is strongly compatible with $\Sigma$ and the 
arithmetic subgroup $\Gamma$ is neat, the quotient set 
$\Gamma\backslash D_{\Sigma}$ can be given the structure of a 
log locally ringed space over $\mathbb C,$ in fact a 
\textit{log manifold} (see (\cite{Kato-Usui-book}, 3.5.7)). Roughly 
speaking, a log manifold is a log locally ringed space over 
$\mathbb C,$ which is locally isomorphic to the ``zero locus" 
of some log differential forms on a log smooth analytic space.  

Informally speaking, Kato and Usui proved the following. 
First, there is a one-to-one correspondence between 
$D_{\Sigma}$ and the set of polarized log Hodge structures of 
the given type. Second, if $\overline{X}\to\overline{S}$ is a 
log smooth family extending the projective smooth family $X\to 
S,$ where $S\subset\overline{S}$ is the complement of a normal 
crossing divisor, then the period map extends to 
$\overline{S}\to M_{\Sigma}.$ We briefly explain the first 
part in the following. 

We shall show how to get a nilpotent orbit from a polarized 
log Hodge structure on a log point (\cite{Kato-Usui-book}, 0.4.24). Let 
$x$ be an fs log point with log structure $M_x.$ Then 
$\overline{M}_x$ is a sharp fs monoid and 
$\overline{M}_x^{\text{gp}}$ if a free abelian group of finite 
rank, say $r.$ Fix $y\in x^{\log}.$ We have 
$x^{\log}=Hom(\overline{M}_x^{\text{gp}},S^1)\simeq(S^1)^r$ 
and hence $\pi_1(x^{\log})=Hom(\overline{M}_x^{\text{gp}}, 
\mathbb Z)\simeq\mathbb Z^r.$ Let 
$\pi_1^+(x^{\log})\subset\pi_1(x^{\log})$ be the subset 
consisting of those homomorphisms $a:\overline{M}_x^{\text 
{gp}}\to\mathbb Z$ that take $\overline{M}_x$ into $\mathbb 
N;$ this subset is an fs monoid.  

Let $(H_{\mathbb Z},Q,F^{\bullet}H_{\mathscr O})$ be a 
polarized log Hodge structure on $x.$ Let 
$(h_i)_{i=1}^n$ be a family of generators for 
$\pi_1^+(x^{\log})$ and fix an $s_0\in\text{sp}(y).$ Let 
$z_1,\cdots,z_r$ be complex numbers, and let 
$s\in\text{sp}(y)$ be such that 
$$
s\Big(\frac{\log(f)}{2\pi i}\Big)-s_0\Big(\frac{\log(f)}{2\pi 
i}\Big)=\sum_{i=1}^rz_ih_i(f),\quad\text{for 
}f\in\overline{M}_x^{\text{gp}}.
$$
Let $N_i:H_{\mathbb Q,y}\to H_{\mathbb Q,y}$ be the logarithm 
of $h_i.$ Then we have 
$$
F(s)=\exp\Big(\sum_{i=1}^nz_iN_i\Big)F(s_0),
$$
which shows that $(F(s))_{s\in\text{sp}(y)}$ is an orbit of 
filtrations under $\exp(\sigma\otimes\mathbb C)$ for 
$\sigma=\sum_i\mathbb R_{\ge0}N_i.$ Moreover, the condition 2) 
in Definition (\ref{LVPHS}) implies that $F(s)\in D$ if 
$\text{Im}(z_i)\gg 0$ for all $i,$ and the condition 3) in Definition 
(\ref{LVPHS}) implies that $NF(s_0)^p\subset F(s_0)^{p-1}$ for 
all $p\in\mathbb Z$ and $N\in\sigma.$ In other words, the 
family $(F(s))_s$ is a $\sigma$-nilpotent orbit. 

\section{The main component of moduli spaces}\label{MainComp}

\subsection*{Moduli: compactness and main components} 
In Section \ref{Satriano}, we gave an overview of F. Kato's work \cite{FKato} in which 
he uses log geometry to compactify the moduli space $\mgn$ of curves.  Specifically, he shows that the moduli 
space of log smooth curves agrees with the Deligne-Mumford compactification $\mgnb$.  The key philosophic idea in 
that section was that since moduli spaces of log smooth objects already includes degenerate objects, it is reasonable 
to expect that such a moduli space is a compactification of the moduli of objects with trivial log structure.

While the Deligne--Mumford space of stable curves $\mgnb$
turns out to be irreducible, 
it is an unfortunate fact of life that if $\mf{X}$ is a 
moduli space of higher dimensional objects, then moduli-theoretic ``compactifications" $\bar{\mf{X}}$ of $\mf{X}$ tend 
to have many irreducible components.  If $\mf{X}$ is irreducible, then it sits entirely within one of these many 
components of $\bar{\mf{X}}$ and so it is natural then to ask if this ``main component'' can itself be given 
a moduli interpretation.

In Section \ref{Curves-big-picture} we stated a second philosophic principle: log geometry controls 
degenerations; that 
is, moduli of log smooth objects does not incorporate ``too many'' degenerate objects.  This provides a type of 
converse to the aforementioned philosophy that moduli of log smooth objects should be compact.  As explained in 
Section \ref{Curves-big-picture}, the log structure gives us a fighting chance to show that our moduli space 
is irreducible (although it is of course too na\"{\i}ve to expect that moduli of log smooth objects is always 
irreducible).  Combining the two principles, one may hope that if $\mf{X}$ is an irreducible moduli space and $\bar{\mf{X}}$ 
a moduli-theoretic compactification, then by appropriately incorporating log structures into the objects parameterized 
by $\bar{\mf{X}}$, one will isolate the main component.

This technique of using log geometry to isolate the main component of a moduli space has been carried out by 
M. Olsson in several different settings.  In \cite{toricHilb}, Olsson gives a moduli interpretation to the 
normalization of the main component of the toric Hilbert scheme.  In \cite{cancmpt}, he isolates the 
normalization of the main component of V. Alexeev's compactification of the moduli space of principally polarized 
abelian varieties given in \cite{alexeev}; he further constructs a moduli-theoretic irreducible compactification of 
the moduli space of abelian varieties with higher degree polarization.  In \cite{K3}, he gives an irreducible modular 
compactification of the moduli space of polarized $K3$ surfaces.

\subsection*{Example: the toric Hilbert scheme}
Our goal in this section is to explain the technique of isolating the main component of a moduli space by following 
Olsson's work \cite{toricHilb}.  We begin with the definition of the toric Hilbert scheme.  Let $k$ be a field 
and let $P$ and $Q$ be finitely-generated integral monoids with $Q$ sharp and $P^{gp}$ and $Q^{gp}$ torsion-free.  
Fix a surjective morphism $\pi:P\rightarrow Q$.  This yields a closed immersion from $A_Q:=\Spec k[Q]$ to 
$A_P:=\Spec k[Q]$, which is $T_Q$-equivariant, where $T_Q$ (resp. $T_P$) denotes the torus associated to $Q^{gp}$ 
(resp. $P^{gp}$).  Consider the functor $\sH$ whose $S$-valued points are diagrams 
\[
\xymatrix{
Z\ar[r]^{i}\ar[dr]_{g} & A_{P,S}\ar[d]\\
 & S
}
\]
where $i$ is a $T_Q$-invariant closed immersion and for every $q\in Q^{gp}$, the $q$-eigenspace of $g_*\sO_Z$ is a 
finitely-presented projective $\sO_S$-module of rank $1$ if $q\in Q$ and rank $0$ otherwise.  By 
\cite[Thm 1.1]{hs}, this functor is representable by a quasi-projective scheme, which we call the 
\emph{toric Hilbert scheme}.

Given a closed subscheme $Z$ of $A_{P,S}$ as above, we can move $Z$ by the action of $T_P$ on $A_{P,S}$.  This yields 
an action of $T_P$ on $\sH$.  Since $Z$ is $T_Q$-invariant, this action factors through $T_K=T_P/T_Q$, where 
$K$ denotes the kernel of $\pi^{gp}$.  We therefore obtain a map
\[
T_K\longrightarrow \sH
\]
by letting $u\in T_K$ act on the distinguished point of $\sH(k)$ given by the closed immersion $A_Q\rightarrow A_P$.  
By \cite[3.6(2)]{chuv}, this map is an open immersion.  Therefore the normalization $\sS$ of the scheme-theoretic 
closure of its image is a normal toric variety, and hence carries a natural fs log structure $\sM_{\sS}$ which 
makes it log smooth over $\Spec k$ (endowed with the trivial log structure).  The goal of \cite{toricHilb} is to 
give a moduli-theoretic interpretation of $(\sS,\sM_{\sS})$.

Consider the functor $\sH^{log}$ on the category of fs log schemes over $k$ whose $(S,\sM_S)$-valued points are 
given by diagrams 
\[
\xymatrix{
(Z,\sM_Z)\ar[r]^-{i}\ar[dr]_{g} & (A_P,\sM_{A_P})\times (S,\sM_S)\ar[d]\\
 & (S,\sM_S)
}
\]
where the underlying maps on schemes defines a point of $\sH(S)$, where $g$ is log smooth and integral, and 
where the map
\[
P\longrightarrow \sM_{(Z,\sM_Z)/(S,\sM_S)}:=\coker(g^*\sM_S\rightarrow\sM_Z)
\]
induced by $i$ factors through $Q$.  The main theorem of \cite{toricHilb} is then 
\begin{theorem}[{\cite[Thm 1.6]{toricHilb}}]
\label{thm:mainthm}
The functor $\sH^{log}$ is representable by $(\sS,\sM_{\sS})$.
\end{theorem}
We explain briefly how Olsson obtains a natural morphism
\[
F:\sH_{\sS}\longrightarrow \sH^{log}
\]
which he then shows is an equivalence; here $\sH_{\sS}$ denotes the functor of points of $(\sS,\sM_{\sS})$.  
Olsson obtains $F$ by showing that the pullback to $\sS$ of the universal family 
over $\sH$ yields a point of $\sH^{log}(\sS,\sM_{\sS})$.  Explicitly, if $i:\sZ\rightarrow A_P\times\sS$ is 
the pullback of the universal family, he constructs a log structure on $\sZ$ as follows.  Since $\sS$ is a 
toric variety with torus $T_K$, it can be covered by open affines of the form $\Spec k[L]$ with $L$ a submonoid of $K$ 
whose associated group is $K$.  Olsson proves (\cite[2.4]{toricHilb}) that over such an open affine, the closed 
immersion $i$ is given by
\[
\sZ\times_{\sS}\Spec k[L]=\Spec k[E_L]\longrightarrow \Spec k[P\oplus L]=A_P\times\Spec k[L],
\]
where $E_L$ is the image of $P\oplus L$ in $P^{gp}$ under the map $(p,\ell)\mapsto p+\ell$.
\begin{example}
\emph{Let $\pi:\N^2\rightarrow\N$ send $e_1$ to $2$ and $e_2$ to $1$.  Then the kernel $K$ of $\pi^{gp}$ is all 
integer multiples of $2e_2-e_1$.  Let $L$ be the submonoid of $K$ consisting of the non-negative multiples 
of $e_1-2e_2$ and let $M$ be the submonoid consiting of the non-positive multiples.  Then $E_L$ is generated 
by $e_1$, $e_2$, and $2e_2-e_1$; hence,}
\[
\Spec k[E_L]\simeq k[x,y,z]/(xy-z^2).
\]
\emph{We see that $E_M$ is freely generated by $e_2$ and $e_1-2e_2$, so 
$\Spec k[E_M]\simeq\mathbb{A}^2_k$.}
\end{example}
We see then that 
$\sZ\times_{\sS}\Spec k[L]$ carries a natural log structure.  These log structures glue to give a log structure on 
$\sZ$ (by \cite[Lemma 2.8]{toricHilb}) which makes $i$ a closed immersion of log schemes.  This therefore yields 
$F$ above.

\subsection*{Example: moduli of K3 surfaces}
In this section we discuss Olsson's work \cite{K3} on the moduli of $K3$ surfaces.  Recall that a surface 
$\underline{X}/\underline{S}$ is called \emph{K3} if for all geometric points $s$ of $\underline{S}$, the 
canonical divisor $K_{\underline{X}_s}$ is trivial and 
$H^1(\underline{X}_{s},\sO_{\underline{X}_s})=0$.  
In his thesis, R. Friedman constructs partial compactifications of the coarse space of the moduli of polarized 
$K3$ surfaces over $\C$.  The key notion in his construction is that of a combinatorial $K3$ surface, which we now 
recall (see \cite[p.2]{FS}).
\begin{definition}
\label{def:combk3}
If $k$ is an algebraically closed field, then a \emph{combinatorial K3 surface} over $\underline{S}=\Spec k$ is a 
$k$-scheme $\underline{X}$ with normal crossings (see Section \ref{D-SS}) which is $d$-semistable, has trivial 
dualizing sheaf, and satisfies one of the following: $\underline{X}$ is a smooth $K3$ surface, $\underline{X}$ is a 
chain of elliptic ruled surfaces with rational surfaces on either end, or $\underline{X}$ is a union of rational 
surfaces where the double curves on each component form a cycle of rational curves with the dual graph of $\underline{X}$ a 
trianglulation of $S^2$.
\end{definition}
In light of the $d$-semistability condition and our discussion in Section \ref{D-SS}, we should expect to be 
able to ``explain'' Friedman's moduli from the point of view of log geometry.  With this in mind, we give the 
following definition of a log $K3$ surface (see \cite[Def 5.1]{K3}).
\begin{definition}
\label{def:logk3}
A morphism $f:X\rightarrow S$ of log algebraic spaces is a \emph{log K3 surface} if $f$ is log smooth and integral, 
$\underline{f}$ is proper, the cokernel of $f^*\sM_S\rightarrow\sM_X$ in the category of integral sheaves of monoids 
is a sheaf of groups, and for every 
geometric point $s$ of $\underline{S}$, we have $\Omega^2_{X_s/s}=\sO_{\underline{X}_s}$, 
$H^1(\sO_{\underline{X}_s})=0$, and $\underline{X}_s$ is a normal crossing variety.  The log $K3$ surface is 
called \emph{stable} if it has no infinitesimal automorphisms.
\end{definition}
As explained in \cite[Rmk 5.3]{K3}, if $X/S$ is a log $K3$ surface, then for every geometric point $s$ of $\underline{S}$, 
the fiber $\underline{X}_s$ is a combinatorial $K3$ surface.

Now that we have a definition of log $K3$ surfaces, we discuss the logarithmic counterpart to the polarization.  
In order to ease the exposition, in what follows we will always assume that our log $K3$ surfaces satisfy the 
technical \emph{special} condition defined in \cite[Def 2.7]{K3}.  
Olsson introduces a notion of logarithmic Picard functor which generalizes the usual Picard functor in the case of 
trivial log structure (see Definition 4.5 and Corollary 5.6 of \cite{K3}):
\begin{definition}
\label{def:logPic}
If $X/S$ is a log $K3$ surface, then the \emph{log Picard functor} $\underline{Pic}(X/S)$ is the sheafification of the 
presheaf associating to any $\underline{S}$-scheme $\underline{T}$ the isomorphism classes of $\sM_{X_T}^{gp}$-torsors 
on $\underline{X}_{{\underline{T}},\textrm{\'et}}$.
\end{definition}
A polarization on a log $K3$ surface $X/S$ is then defined to be (\cite[Def 5.7]{K3}) a morphism 
$\lambda:\underline{S}\rightarrow\underline{Pic}(X/S)$ such that on each geometric fiber $\underline{X}_s$, there is a 
line bundle $\mathcal{L}$ which lifts $\lambda_s$ to $H^1(\underline{X}_s,\sO^*_{\underline{X}_s})$ and satisfies the 
following.  There is some $N>0$ such that $\mathcal{L}^N$ is generated by global sections and the map defined by 
$\mathcal{L}^N$ only contracts finitely many curves to points.

With these definitions in place, Olsson fixes a positive integer $k$ and considers the stack $\mathbb{M}_{2k}/\mathbb{Q}$ 
whose fiber over $\underline{T}$ is the groupoid of triples $(\sM_T,X/(\underline{T},\sM_T),\lambda)$ where 
$\sM_T$ is a log structure on $\underline{T}$ and $(X/(\underline{T},\sM_T),\lambda)$ is a stable polarized log $K3$ 
surface such that $\lambda^2_t=2k$ for every geometric point $t$ of $\underline{T}$.  Note that $\mathbb{M}_{2k}$ 
carries a natural log structure given by base log structure $\sM_T$ in each fiber.  One of the main results of 
\cite{K3} is then:
\begin{theorem}[{\cite[Thm 6.2]{K3}}]
\label{thm:K3cmpt}
The stack $\mathbb{M}_{2k}$ is smooth, log smooth, Deligne-Mumford, and contains an open substack 
$\mathbb{M}_{2k}^{sm}$ parameterizing polarized $K3$ surfaces in the classical sense.  The compliment of 
$\mathbb{M}_{2k}^{sm}$ in $\mathbb{M}_{2k}$ is a smooth divisor and the induced log structure agrees with the 
natural one on $\mathbb{M}_{2k}$.
\end{theorem}

\section{Twisted curves and log twisted curves}\label{Roots}
Twisted curves are a central object in the theory of {\em twisted stable maps} \cite{Abramovich-Vistoli,Chen-Ruan,Abramovich-Olsson-Vistoli}: in order to have a complete moduli space of stable maps $C \to X$ of type $\Gamma$, where $X$ is a proper tame stack with projective coarse moduli space and $\Gamma=(g,n,\beta)$ are the relevant discrete data, one must allow the curve $C$ itself to be a certain type of stack, called {\em twisted curve}.

The original treatments of twisted  curves relied on ad-hoc methods. The more recent approach of \cite{Abramovich-Olsson-Vistoli} relies on a method introduced in \cite{Olsson-curves}, which uses a construction with logarithmic structures.

\subsection*{Twisted curves} 

For simplicity we will stick with the case of Deligne--Mumford stacks.

First consider the geometric objects: fix an algebraically closed field $k$.

\begin{definition}\label{Def:twisted-curve}
A {\em twisted curve} over $k$ is a tame, purely 1-dimensional Deligne-- Mumford stack $\cC /k$, with at most nodes as singularities, satisfying the following conditions:
\begin{enumerate}
\item Let $\pi: \cC \to C$ be the morphism to the coarse moduli space. Then $\cC^{sm} = \pi^{-1} C ^{sm}$, and $\pi: \cC \to C$ is an isomorphism over a dense open subset of $C$.
\item Consider a node $\bar x \to C$, where the strictly henselian local ring  $\cO_{C,\bar x}$ is the strict henselization of $k[x,y]/(xy)$. Then 
$$\cC \times_C\Spec \cO_{C,\bar x}\ \simeq\
 \left[ \Spec \cO_{C,\bar x}[z,w]/(zw,z^m-x,w^m-y) \ \ \bigg/ \ \ \mu_m\right],$$
 where $\zeta \in \mu_m$ acts by $(z,w) \mapsto (\zeta z, \zeta^{-1}w)$.
\end{enumerate}
\end{definition}
An action such as (2) above is called {\em balanced} - it is crucial to our discussion of log structures below. Note that $\cC$ may have a stack structure at isolated smooth points as well - such points will behave like $[\AA^1/\mu_a]$, where $\mu_a$ acts by multiplication.

Over a general base $S$ twisted curves are detected by their geometric fibers: a {\em twisted curve} $\cC \to S$ is a flat, tame Deligne--Mumford stack locally of finite presentation, all of whose  geometric fibers are twisted curves as in the definition above.

The {\em genus} of $\cC$ is simply the genus of $C$. One typically needs to consider $n$-pointed twisted curves, where the markings are described in families as follows:

\begin{definition}
An {\em $n$-pointed twisted curve} $\cC /S$ marked by disjoint closed substacks $\{\Sigma_i\}_{i=1}^n$ in $\cC$ is assumed to satisfy the following:
\begin{enumerate}
\item the $\Sigma_i$ are contained in the smooth locus $\cC^{sm}$,
\item each  $\Sigma_i$ is a tame \'etale gerbe over $S$, and 
\item $\cC_{gen} := \cC^{sm} \smallsetminus \cup_i\Sigma_i \ \ \longrightarrow\ \  C $ is an open embedding.
\end{enumerate}
\end{definition}
\begin{remark}
When $S=\Spec k$ where $k=\bar k$, then $\Sigma_i = B\mu_{a_i}$, and moreover $a_i$ is locally constant in families.
\end{remark}

\begin{remark}
When $(\cC/S, \{\Sigma_i\})$ is an $n$-pointed twisted curve, then the coarse moduli space of $\Sigma_i$ is isomorphic to $S$. This means that the composite morphism $\Sigma_i\to \cC \to C$ factors through a section $p_i: S \to C$. It follows that $(C,\{p_i\})$ is an $n$-pointed curve in the usual sense. This gives a functor 
$$(\cC/S, \{\Sigma_i\}) \mapsto (C,\{p_i\})$$
\end{remark}

One can ask oneself: what does one need in order to recover a twisted $n$-pointed curve  $(\cC/S, \{\Sigma_i\})$ from a usual $n$-pointed curve $(C,\{p_i\})$? In other words, can we enrich the functor above to something like
$$(\cC/S, \{\Sigma_i\}) \mapsto (C,\{p_i\})+\text{?}$$
which is nice and explicit and actually an equivalence of categories?

The stack structure at the marking definitely needs the data of the integers $a_i$, but in fact this is all that is necessary for the markings: near $p_i$, the curve $\cC$ is canonically isomorphic to the root stack $C(\sqrt[a_i]{p_i})$. If $x$ is a local generator of the ideal of $p_i$, then Zariski locally we have
$$\cC\ \  \simeq \ \ \left[ \Spec \cO_C[z]/(z^{a_i} -x)\ \  \bigg/\ \  \mu_{a_i}\right].$$

The story is a bit more interesting at a node. It has to be - a twisted curve $\cC$ with a node of index $m> 1$ has ``ghost" automorphisms in $\mu_m$ which are not detectible on the coarse 
curve $C$: using the local coordinates given in Definition \ref{Def:twisted-curve} (2), the $\mu_m$ action 
$$(z,w) \mapsto (\zeta z, w)$$
on $\Spec \cO_{C,\bar x}[z,w]/(zw,z^m-x,w^m-y)$ commutes with the action defining the quotient stack. It  therefore descends to a nontrivial action on $\cC$ which becomes trivial on the coarse moduli space $C$. 

\subsection*{Log twisted curves} 
Let $X$ be a Deligne--Mumford stack. Recall from Definition \ref{Def:fine} that a fine log structure $M$ on $X$ is said to be {\em locally free} if for every geometric point $\bar x \to X$ we have that the characteristic sheaf $\overline M_{\bar x}$ is isomorphic to $\NN^r$ for some $r$.

In this situation we say that a morphism of sheaves of monoids $M \to M'$ is {\em simple} if for   every geometric point $\bar x \to X$  we can identify the map as the diagonal map 
$$\xymatrix{ \overline M_{\bar x} \ar[rr]\ar^\simeq[d]&& \overline M'_{\bar x}\ar^\simeq[d]\\
\NN^r\ar[rr]^{(m_1,\ldots,m_r)}&& \NN^r}$$
where all $m_i$ are prime to the characteristic of the field.

\begin{definition}
An {\em $n$-pointed log twisted curve over $S$} is the data
$$(C/S, \{\sigma_i,a_i\},\ell:M_S \to M_S')$$ where
\begin{itemize}
\item $(C,\{\sigma_i\})/S$ is an $n$-pointed nodal curve.
\item $M_S$ is the canonical log structure coming from the family $(C,M_C) \to (S,M_S)$ (see section \ref{Satriano}).
\item $a_i: S \to \ZZ_{>0}$ are locally constant, with $a_i(s)$ invertible in the residue field $k(s)$. 
\item $\ell: M_S\to M_S'$ is a simple morphism.
\end{itemize}
\end{definition}

And we have the following:
\begin{theorem}[\protect{\cite[Theorem 1.8]{Olsson-curves}}]
The fibered category of $n$-pointed twisted curves is naturally equivalent to the stack of $n$-pointed log twisted curves.
\end{theorem}

The picture is as follows: we have already noted that we can replace a marking $p_i$ by a stacky marking just using the data $a_i$.  Now the $j$-th node which looks like $\Spec \cO_S[x,y]/(xy-t)$ needs to be replaced by $\left[\Spec \cO_S[z,w]/(zw-t^{1/m_j})\ \big/\ \mu_{m_j}\right]$, and the data is encoding by deviding the $j$-th generator of $\NN^r$ by $m_j$.


\begin{remark}
We can decompose the stack according to $a_i$:
$$\xymatrix{\mathop{\coprod}\limits_{\underline a}^{\vphantom{a}} \cM^{tw}_{g,n,\underline a} \ \ \ \ \ = 
& \cM^{tw}_{g,n}\ar[d]&= \ \ \ \text{ stack of $n$-pointed twisted curves}\\
&\cS_{g,n}&= \ \ \ \text{ stack of $n$-pointed nodal curves}\ \ 
} $$
and it can be deduced form the theorem that $\cM^{tw}_{g,n,\underline a} $ is obtained from $\cS_{g,n}$ using a root construction applied to the boundary divisor of $\cS_{g,n}$, see \cite[Remark 1.10]{Olsson-curves}. 

In fact we have to apply all possible roots, accounting for all possible twisting of nodes, and glue together, so $\cM^{tw}_{g,n,\underline a} $ is highly non-separated.
\end{remark}

Below we sketch the main ideas in proving this. We stress that the assumption that our twisted curves are {\em balanced} is crucial - the case of unbalanced curves has not been treated.

\subsection*{From twisted curves to log twisted curves}

Fix a twisted curve $f:\cC \to S$.

We can follow F. Kato \cite{FKato}, giving log structures on nodal curves: consider all possible triples $$(M_S,M_{\cC}, f^\flat: f^{-1}M_S \to M_\cC)$$ such that 
\begin{enumerate}
\item $(\cC,M_\cC) \to (S,M_S)$ is log smooth;
\item $M_\cC,M_S$ are locally free; and
\item for all geometric points $\bar x \to \cC$ mapping to nodes we have
$$\xymatrix{\overline M_{\cC,\bar x}\ar[r]^\sim & \NN^{r-1}\oplus \NN^2 \\
\overline M_{S,f(\bar x)}\ar[r]^\sim\ar[u] & \NN^{r-1}\oplus \NN\ar[u]_{id\oplus \Delta}
}$$
\end{enumerate}

If $S = \Spec k$ we get a natural map $$\NN^\text{number of nodes} \to \overline M_S. $$ We say that $(M_\cC,M_S,f^\flat)$ is {\em special} if for every geometric point $\bar s \to S$ this map is an isomorphism. 

A result \cite[Theorem 3.6]{Olsson-curves} analogous to F. Kato's \cite[Theorem 2.1]{FKato} says that there is a unique special triple $(M_\cC,M'_S,f^\flat)$ associated to $f: \cC \to S$. Analyzing the coarse moduli space we obtain a unique diagram

$$\xymatrix{
(\cC, M_\cC)\ar[r]\ar[d] & (C, M_C)\ar[d]\\
(S, M'_S) \ar[r]^{(id,\ell)} & (S, M_S)
}$$
where $\ell: M_S \to M'_S$ is simple. In particular we obtain a log twisted curve
$$(C/S, \{\sigma_i,a_i\}_{i=1}^n, \ell: M_S \to M'_S)$$
\subsection*{From log twisted curves to twisted curves}
Now we fix a log twisted curve $(C/S, \{\sigma_i,a_i\}_{i=1}^n, \ell: M_S \to M'_S)$.
In particular we have the log smooth curve 
$$(C, M_C) \to (S, M_S)$$ which is the coarse moduli space of a putative twisted curve. 
We want to describe $\cC/C$ as the stack parametrizing natural objects over $T \to C$. Here it is! If we denote the relevant maps as follows
$$\xymatrix{T\ar[r]^s\ar[dr]_h& C \ar[d]\\& S,}$$ then $\cC$ is the groupoid of diagrams

$$\xymatrix{h^*M_S\ar[r]^\ell\ar[d] & H^* M'_S\ar[d]\\
    s^*M_C\ar[r]^k&M'_C,
}$$
where
\begin{enumerate}
\item $k$ is simple and for every geometric point $\bar t \to T$ the map 
$\overline M'_{S,\bar t} \to \overline M'_{C,\bar t}$ is either an isomorphism (at a general point), or of the form
$$\xymatrix{\NN^r\ \ \ar@{^(->}[rr]^{(id,0)}&&\ \  \NN^r\oplus \NN}$$ (at a marked point), or 
$$\xymatrix{\NN^{r-1} \oplus \NN\ \ \ar@{^(->}[rr]^{(id,\Delta)}&&\ \  \NN^{r-1}\oplus \NN^2}$$ (at a node).
\item for all $i$ and geometric point  $\bar t \to T$ with $s(\bar t)\subset \sigma_i(S)\subset C$ the group 
$$\text{Coker} ( \overline{M}^{gp}_{C,\bar t} \to \overline{M'}^{gp}_{C,\bar t})$$ is cyclic of order $a_i$. 
\end{enumerate}

\section{Log stable maps}\label{Stablemaps}
\subsection*{From curves to maps and expansions}
\subsubsection*{Curves}
In section \ref{Satriano} we discussed how prestable curves can be encoded as log smooth curves, and how in particular the stack of Deligne--Mumford stable curves can be interpreted as a logarithmic stack, representing log smooth stable curves over the category of fine  and saturated log schemes. Stability in this situation just means that $\Omega^1_{(C,\cM_C)}$ is an ample line bundle. This is the same as saying that $\omega_C(D)$ is ample, where $D$ is the divisor of markings.

\subsubsection*{Maps}
Kontsevich \cite{Kontsevich} introduced the moduli  of {\em stable maps} of prestable curves into a projective target variety $X$. This is a proper Deligne--Mumford stack $\ocM_\Gamma(X)$ having projective coarse moduli space, where $\Gamma = (g,n,\beta)$ is the relevant numerical data: genus, number of markings and homology class of the image curve.
It parametrizes maps $f: C \to X$, where this time stability means that $\omega_C$ is $f$-ample. 

Kontsevich's moduli space has the property that it carries a {\em perfect relative obstruction theory} giving rise to a {\em virtual fundamental class} $[\ocM_\Gamma(X)]^{vir}$, see \cite{Li-Tian, Behrend-Fantechi}. This is a key ingredient in defining Gromov--Witten invariants, with their applications in enumerative geometry and theoretical physics. The simplest GW invariants are 
$$\langle\gamma_1\cdots \gamma_n\rangle := 
\int_{[\ocM_\Gamma(X)]^{vir}} e_1^*\gamma_1\cdots e_n^*\gamma_n,$$ 
where $\gamma_i\in H^*(X,\QQ)$ and $e_i: \ocM_\Gamma(X) \to X$ are the natural {\em evaluation maps} at the $n$ markings.

So far no logarithmic structures are necessary.

\subsubsection*{Degenerations}

Among the methods of computing GW invariants, the {\em degeneration formula} is among the most powerful ones. It was introduced by A.-M. Li and Y. Ruan in symplectic geometry \cite{Li-Ruan}, see also Ionel-Parker \cite{Ionel-Parker1,Ionel-Parker2}. For our purposes, Jun Li's treatment in algebraic geometry \cite{Li,Li-degeneration} is most relevant.

We are interested in the invariants of the smooth projective variety  $X$. Since these are deformation invariant, it is natural to consider a degeneration of  $X$, with smooth total space, into a union $X_0 = Y_1\sqcup_D Y_2$ of two smooth varieties $Y_i$ meeting transversally along a smooth divisor $D$. One wants to show that invariants of $X$ coincide with suitably defined invariants of $X_0$, and these in turn can be computed in terms of appropriately defined invariants of $Y_i$ relative to $D$.

This is where logarithmic structures begin to show up, but there is still some way to go.

\subsubsection*{Perfect obstruction and Li's approach}
The difficulty with the degeneration is precisely the fact that the variety $X_0$ is singular, and therefore the natural obstruction theory on the moduli space of stable maps is not perfect in general.  A similar situation occurs when considering the pair $(Y_i,D)$, but we will not get into this discussion.

The problem occurs when a component of the source curve $C$ maps entirely into $D$.
 
 Jun Li's approach uses 
  {\em expanded degenerations}. There are similar ideas in  \cite{Li-Ruan}, but the symplectic approach builds in deformations of Cauchy-Riemann equations and has, at least on the surface, a significantly different flavor.

The idea is, that just as in stable pointed curves, if a marking travels towards a node one sprouts a new component of the curve, Jun Li says that when a component of $C$ travels into  $D$ we can let $X_0$ sprout a new component. 

\subsubsection*{Expansions}
Here is how the new component looks like. 
Denote by $N_{D\subset Y_{i}}$ the normal bundle of $D$ in $Y_{i}$. Since the total space is smooth, we have $N_{D\subset Y_{1}}\cong N_{D\subset Y_{2}}^{\vee}$. Let $\PP=\PP roj_{D}(1_{D}\oplus N_{D\subset Y_{1}})$. We have $\PP\cong \PP roj_{D}(N_{D\subset Y_{2}}\oplus 1_{D})$ so we can denote by $D^{+}$ and $D^{-}$ the smooth divisors in $\PP$ which correspond to the normal bundle $N_{D\subset Y_{1}}$ and $N_{D\subset Y_{2}}$ respectively. Note that the divisor $D^{+}$ and $D^{-}$ are canonically isomorphic to $D$.

We now glue things up.
Let $\PP_{i}$ for $i\in \N$ be copies of $\PP$. We can glue $Y_{1}$ and $\PP_{1}$ along $D$ and $D^{-}$ respectively, $\PP_{i}$ and $\PP_{i+1}$ along $D^{+}$ and $D^{-}$ respectively, and $\PP_{n}$ and $Y_{2}$ along $D^{+}$ and $D$ respectively. We denote the resulting gluing by 
\begin{equation}
X_{0}[n]= Y_{1} \bigsqcup_{D_{1}} \PP_{1} \bigsqcup_{D_{2}} \cdots \bigsqcup_{D_{n}} \PP_{n} \bigsqcup_{D_{n+1}}Y_{2}, 
\end{equation}
where $D_{1},\cdots, D_{n+1}$ are the disjoint singular loci of $X_{0}[n]$.

Such a beast is known as an {\em expanded degeneration}, or by the more folksy name {\em an $n$-accordion}.

This led Jun Li to define {\em degeneration stable maps with target $X_0$}  as {\em nondegenerate maps $C \to X_0[n]$:} a map is nondegenerate if no component of $C$ maps into any of the $D_i$.

\subsubsection*{Predeformability}

But here Jun Li meets another phenomenon, already present in the space of admissible covers of Harris--Mumford \cite{Harris-Mumford}: nondegenerate maps to $C \to X_0[n]$ have many redundant components which have nothing to do with maps to the generic fiber $X$. Near a singular point of $X_0[n]$ which looks like $\{xy=0\}$, a curve map locally deforms to a smoothing of $X_0$ if and only if the curve looks like $\{uv=0\}$, and the map given by $x=u^m, y=v^m$. 
Such nice maps are called {\em predeformable}. But predeformable maps are clearly not open in the space of maps - they are actually closed among nondegenerate maps. This means that the restricted obstruction theory on them is ``wrong" - definitely not perfect

Of course the virtual fundamental class of Gromov--Witten theory has no problem dealing with the total moduli space with its extra components, but these extra components do get in the way of decomposing invariants of $X_0$ in terms of $(Y_i,D)$. So we really do want to stick by predeformable maps. 

\subsection*{Logarithmic methods: from Jun Li to Bumsig Kim}
\subsubsection*{Predeformable deformations}

At this point Jun Li's reasoning arrives at a point where a new obstruction theory on predeformable maps is needed. Having read this chapter, the reader will immediately recognize that 
\begin{enumerate} 
\item nodal curves are log-smooth, 
\item $n$-accordions are d-semistable and admit a canonical log-smooth structure, and
\item predeformble maps can be viewed as log maps from the log-smooth curve to the d-semistable target.
\end{enumerate}

Jun Li also recognized this fact, as did Shin Mochizuki before him \cite{Mochizuki} when he in his turn revisited the space of admissible covers of Harris--Mumford. What he lacked at the time was a formalism for logarithmic deformation theory of singular spaces, such as the moduli space itself: a perfect obstruction theory is a two-term complex mapping to the cotangent complex of the moduli space, but the moduli space is, as usual highly singular, even taking its logarithmic structure into account.

So Jun Li resorted to an ad-hoc construction of his perfect obstruction theory. This is the most difficult part of his work.

In section \ref{LogDef} we saw how log deformation theory works in the necessary generality. This is where Bumsig Kim's paper \cite{Kim} comes in: he provides a correct formalism for nondegenerate logarithmic stable maps into expanded degenerations, and shows that it carries a perfect obstruction theory. The degeneration formula in this formalism has been worked out by one of us (Q. Chen), and should appear as part of a larger project indicated below.

One aspect that deserves mention is Kim's notion of {\em minimal log strutures} on maps.  Recall that the log stack $\mgnb$ can be constructed as a stack over the category of fine  saturated log schemes, whose objects over $S = (\underline S,M)$ are log smooth stable pointed curves over $S$. But in order to exhibit the underlying stack, one needs to use the {\em canonical} log structure on $\underline S$, which is initial among all possible ones carrying the log smooth curve. 

Kim describes his stack similarly - given a predeformable map $\underline C \to \underline X$ of underlying schemes over $\underline S$,  it amounts to describing what he calls {\em minimal} log structure $S$ on $\underline S$ carrying a log map $C \to X$. Kim does this by explicitly describing the combinatorial structure of such log structures.

\subsection*{Unexpanded log maps: from Siebert into the future} 
A very different approach was proposed in a 2001 lecture by Bernd Siebert, but lay dormant for almost a decade.

The point is this. If one embraces logarithmic structures, and logarithmic maps from log smooth curves to some logarithmic scheme, then expansions are no longer necessary. Defined correctly, the space of log maps automatically has a perfect log-obstruction theory, which in view of sections \ref{LogStacks} and \ref{LogDef} can be viewed as an obstruction theory relative to the stack $LOG$. This automatically results in a virtual fundamental class.  

With this way of thinking, one can consider much more general logarithmic stable maps, gaining access to invariants of much more general degenerations of varieties. This has been a desired goal for a number of years.

So what is the correct definition? Consider a fine and saturated log scheme $X$. Following the work of F. Kato \cite{FKato} as discussed in section \ref{Satriano}, one  comes up with a definition of a category $\ocM_\Gamma(X)$ fibered over the category $LSch_{fs}$ of fine and saturated log schemes: an object over a fine saturated log scheme $S$ is a log smooth curve $C \to S$ and a log map $f:C \to X$. We further require it to be {\em stable}: the line bundle  of logarithmic differentials $\Omega^1_{C/S}$ is required to be $f$-ample. This is tantamount to requiring  the map of underlying schemes to be a stable map. 

So the main claim is: $\ocM_\Gamma(X)$ is represented by a logarithmic Deligne--Mumford stack with projective coarse moduli space. In fact this stack is proper and quasi-finite over the usual stack of stable maps $\ocM_\Gamma(\underline X)$ of the underlying scheme. As in the discussion of Kim's work, the underlying stack can be viewed as a moduli of log maps with {\em minimal} log structure.

This is the subject of current work - of Gross and Siebert on the one hand and of two of us (mainly Chen, and to a lesser extent Abramovich) on the other, so it would not be appropriate to get into details until definite results appear.   

Let us instead  put this in a larger context. Consider logarithmic schemes $Z \to B$ and $X$,  and assume we are given a morphism of underlying schemes $\underline f: \underline Z \to \underline X$. We can define a category $\Lift_{\underline f}$ fibered over $LSch_{fs}$ whose objects over a fine saturated log $B$-scheme $S\to B$ are lifts $f_S:Z_S \to X$ of the morphism of underlying schemes $\underline f_{\underline S}:\underline{Z}_{\underline S} \to \underline X$. 

One can ask the following general questions: 
\begin{question}
\begin{enumerate}
\item Under what conditions is $\Lift_{\underline f}$ a log stack  locally of finite type over $B$?
\item What natural numerical data cut out a substack of finite type?
\item Under what conditions is the result proper?
\end{enumerate}
\end{question}

We want to stress our belief that this question is natural, important  and quite tractable. For instance, the case where $B = \Spec \CC$ with trivial structure and $Z = X = \Spec (\NN \to \CC)$ the result is a countable union of components. It  is similar in nature to an inertia stack. 
The more general case where 
$B =X$, $Z = X\times \Spec (\NN \to \CC)$ and $\underline f$ is the diagonal, is the relevant analogue of the inertia stack of $X$. It is  important for Gromov--Witten theory - up to $\CC^*$ action the result is the natural target for evaluation maps associated to log smooth curves. Its components account for the contact orders of relative stable maps. Further examples of a similar nature govern gluing of nodes of log-smooth curves.


\end{document}